\documentclass[numerical,a4paper,nofootinbib]{article}
\pdfoutput=1
\usepackage[top=3cm, bottom=3cm, left=2.5cm, right=2.5cm]{geometry}

\usepackage{graphicx}
\usepackage{calc}
\usepackage{amsmath}
\usepackage{bbm}
\usepackage{amsfonts}
\usepackage{amssymb}
\usepackage{amsbsy}
\usepackage[T1]{fontenc}
\usepackage{theorem}
\usepackage{stmaryrd}
\usepackage{tipa}
\usepackage{textcomp}
\usepackage{subfigure}
\usepackage{float}
\usepackage{epsfig}
\usepackage{lmodern}
\usepackage{framed}
\usepackage[subfigure]{tocloft}
\usepackage{subeqnarray}
\usepackage{alphalph}
\usepackage[refpage,cfg,prefix]{nomencl}
\usepackage{multirow}
\usepackage{xifthen}
\usepackage{enumerate}
\usepackage{color}
\usepackage{natbib}  
\usepackage{wasysym}
\usepackage{bibentry}
\usepackage{verbatim}
\usepackage{todo}
\usepackage{url}
\usepackage[symbol]{footmisc}
\usepackage{setspace}
\usepackage{xr}
\usepackage{prettyref}

\usepackage{color,colortbl}
\usepackage{relsize}

\makenomenclature

\newcommand\norm[1]{\left\lVert#1\right\rVert}

\newcommand{\listofalgorithms}{\textbf{\Huge{List of Algorithms}}}

\newlistof{Algorithm}{exp}{\listofalgorithms}

\newif\ifletter
\def\ack{\ifletter\bigskip\noindent\ignorespaces\else
    \section*{Acknowledgments}\fi}

\definecolor{light-gray}{gray}{0.4}
\usepackage{epstopdf}

\begin{document}

\title{\Huge{Pair correlation functions for identifying spatial correlation in discrete domains}}
\author{\textbf{Enrico Gavagnin\footnote{These authors contributed equally to this work. Corresponding authors: EG: e.gavagnin@bath.ac.uk, JO: j.owen@bath.ac.uk }   \,, Jennifer P. Owen$^*$ and Christian A. Yates} \\ \small{\textit{Department of Mathematical Sciences},}\\ \small{\textit{University of Bath, Claverton Down, Bath, BA2 7AY,  UK}}}%
\date{}

\maketitle

\begin{abstract}

Identifying and quantifying spatial correlation are important aspects of studying the collective behaviour of multi-agent systems. Pair correlation functions (PCFs) are powerful statistical tools which can provide qualitative and quantitative information about correlation between pairs of agents. Despite the numerous PCFs defined for off-lattice domains, only a few recent studies have considered a PCF for discrete domains. Our work extends the study of spatial correlation in discrete domains by defining a new set of PCFs using two natural and intuitive definitions of distance for a square lattice: the taxicab and uniform metric. We show how these PCFs improve upon previous attempts and compare between the quantitative data acquired. We also extend our definitions of the PCF to other types of regular tessellation which have not been studied before, including hexagonal, triangular and cuboidal. Finally, we provide a comprehensive PCF for any tessellation and metric allowing investigation of spatial correlation in irregular lattices for which recognising correlation is less intuitive.
\\

\end{abstract}

\section{Introduction}
\label{sec:introduction}
A system of agents is considered in a state of spatial correlation if, given any agent in this system, the likelihood that there are other agents at a certain close distance is either increased or decreased with respect to the situation in which the agents are distributed uniformly at random. 
Spatial correlation is a dominant feature of many biological and physical systems \cite{steinberg1996adh,cavagna2010sfc,green2010nlm,thomas2005hsc,chalub2006mhc,othmer1997abc,stevens2000sca,murakawa2015cmc,simpson2013emi,treloar2013mtd,hinde2010vpc,binder2015qtf}. 
For example, in cell biology, spatial correlation can be seen in the form of patterns on animal fur or fish skin \cite{steinberg1996adh,green2010nlm,deutsch2007cam}. In a clinical setting, cell aggregation is a characteristic feature of melanoma and its identification is essential for early diagnosis and effective therapy \cite{friedman1985edm,weinstock2000edm}. Resource competition in ecology can lead to spatial correlation in the form of segregation, for example, in ant nest displacement in a competitive environment \cite{bourke1995sea}. In epidemiology, spatial correlation can be observed in the occurrence of disease across different geographical regions \cite{keeling1999els}.    

The same spatial configuration can have different origins. For example, spatial aggregation in cell biology can be caused by a result of cell-to-cell adhesion \cite{murakawa2015cmc}, external signals, as in \textit{chemotaxis} \cite{chalub2006mhc,othmer1997abc}, or even slime following \cite{stevens2000sca}. Alternatively, cells may form clusters during development due to a combination of a high proliferation rate and a low movement rate \cite{simpson2013emi}.
Given a system exhibiting spatial correlation, one may hypothesise an underlying mechanism responsible for these properties. These assumptions may form the basis of a mathematical model that can be simulated for the purpose of testing. Quantifying spatial correlation in both the simulation and observed experimental data can be a way to connect these studies and to validate or disprove such a theory. As a result, a great number of statistical tools have been developed in the past decade to analyse and measure spatial correlation \cite{diggle1976sas,illian2008sam,binder2013qss,agnew2014dmc,binder2015sap,binder2015qtf,hackett2012gis}. Among the most popular are pair correlation functions (PCF) \cite{bahcall1983scf,donev2005pcf,young2001rpc,hinde2010vpc,binny2016ccb,raghib2011mme,binder2013qss,agnew2014dmc,binder2015sap,binder2015qtf} and Fast Fourier Transform (FFT) \cite{bone1986fpa,takeda1982ftm}. In this paper we focus our attention on the study of spatial correlation using PCFs.

Given a system of agents, a PCF determines whether pairs of agents are more or less likely to be found with a given separation than in the situation in which the agents are positioned uniformly at random in the domain. 
A PCF is considered effective if it fits two main criteria. Firstly, the PCF distance metric should be well-defined, but most importantly be readily interpretable in the context of the system considered. This criteria is essential so that in the case of correlation (aggregation or segregation), the PCF can be used to obtain more details about the spatial configuration. For example, if the system exhibits aggregation, the PCF should be able to provide a measure for the average size of the clusters and their pairwise separation.  
Secondly, the PCF should be correctly calibrated. The PCF should be able to distinguish between three basic types of configurations: spatial randomness, aggregation and segregation. For this, the PCF should be normalised correctly \textit{i.e.} the PCF should return the value unity at all pairwise distances (no correlation) when applied to a uniformly distributed set of agents. If the PCF is not normalised correctly, a spatially random set of agents may be incorrectly identified as a correlated system. This inconsistency makes PCF profiles hard to interpret. 

Depending on the type of investigation, the mathematical framework can either be continuous (off-lattice) or discrete (on-lattice). The corresponding PCF has to be defined in accordance with the given framework. Despite the abundance of PCFs defined for off-lattice domains \cite{young2001rpc,hinde2010vpc,binny2016ccb,raghib2011mme}, only a few recent studies have defined a PCF for domains partitioned with a lattice \cite{binder2013qss,agnew2014dmc,binder2015sap,binder2015qtf}.
On-lattice PCFs often assume exclusion properties, that is, that each lattice site in a domain can be occupied by at most one agent at any given time. This is consistent with typical on-lattice correlation studies, such as those designed to quantify correlation in binary pixelated images, or to determine spatial correlation in exclusion processes simulated using a discrete domain. 

Currently, there are two PCFs defined on-lattice. The first is a naive approach consisting of applying the classic off-lattice PCF to lattice-based systems. We refer to this from now on as the \textit{Annular} PCF. In the Annular PCF, given some small positive $\delta$, the number of agents at distance $m$ from a focal agent is defined as the number of agents whose centres lie in the annulus ($m-\delta$,$m$], where distance is defined by the Euclidean metric. A limitation of this method is that, whilst the normalisation is a good approximation for a continuous domain (see Section \ref{sec:previous} for more details), it is poor in the case of a discrete domain, thus the PCF is not correctly calibrated.
In other, more recent work, \citet{binder2013qss,binder2015sap} defined a PCF specifically designed for a two-dimensional square on-lattice exclusion process which we will refer to as the \textit{Rectilinear} PCF (see Section \ref{sec:previous} for more details). Whilst their approach correctly identifies the spatial correlation in many examples, due to an anisotropic definition of distance, spatial structures which are biased in either Cartesian directions can remain unidentified by this PCF. 
To summarise, to the best of our knowledge, a discrete isotropic PCF with correct normalisation does not currently exist in the literature. 

In this paper, we extend the study of pairwise spatial correlation for on-lattice exclusion processes which tackles the flaws of previous PCFs.
We define new isotropic PCFs for a square lattice on which distance is defined using two of the most natural and intuitive metrics for a discrete domain: the taxicab and uniform metric.  
We call these the \textit{Square Taxicab} PCF and \textit{Square Uniform} PCF, respectively, after the square lattice set up and metric type. We define them in both the non-periodic and periodic boundary cases. 
Using synthetically generated data, we demonstrate that our PCF can correctly distinguish between spatial randomness, aggregation and segregation. Furthermore, we show that it can also provide quantitative information about the structure of the system, such as approximate aggregate size or segregation distance both in the short and long scales. Moreover, we investigate how the choice of metric, uniform or taxicab, can affect this quantitative information. 
We demonstrate that our PCFs represent a significant improvement on previous on-lattice PCFs by showing that, firstly, our method is correctly calibrated (unlike the Annular PCF) and secondly, that it can identify anisotropic patterns of the type that are routinely missed by the Rectilinear PCF.
As a natural extension, we define PCFs for higher dimensions and other types of tessellations (cubic, triangular and hexagonal) which have not been considered previously. We name these the \textit{Triangle} PCF, \textit{Hexagon} PCF, the \textit{Cube Taxicab} and \textit{Cube Uniform} PCFs after the lattice set up and metric type. 

Finally, we extend the concept of a PCF by introducing the \textit{General} PCF. This PCF can be defined using any metric, 
on any discrete domain type, with the caveat that it is more computationally expensive. We give an example of how we can use this PCF on a discrete irregular lattice (both tessellation and domain shape) where we define adjacent sites to be at unit distance from one another. We show how our PCF can identify aggregation and segregation on an irregular domain using some synthetic examples.  All the MATLAB codes which accompany the paper can be found as online supplementary material (see the online published version \citep{gavagnin2018pcf}).

The paper is organised as follows. In Section \ref{sec:previous} we discuss the successes and limitations of previous on-lattice PCFs. In Section \ref{sec:PCF} we introduce our novel Square Taxicab and Square Uniform PCF.  We apply our Square Taxicab and Square Uniform PCF to some relevant examples and make comparisons with previous on-lattice PCFs from the literature in Section \ref{sec:results}. In Section \ref{sec:tria-hex} we define the Triangle PCF, Hexagon PCF, the Cube Taxicab and Cube Uniform PCF. We extend our PCF to more generic and possibly irregular lattices by defining the General PCF in Section \ref{sec:graph}\footnote{For reference, in Section \ref{SUPP-sec:summary_norm} of the Supplementary Materials we supply a table summarising all the formulae for the normalisations of our PCFs.}. Finally, we conclude, in Section \ref{sec:conclusions}, by summarising the relevance of our results and discussing potential avenues for future work.
\section{Existing on-lattice Pair Correlation Functions}
\label{sec:previous}
 In this section we provide a summary of the only two existing PCFs defined for discrete domains: the Annular PCF and the Rectilinear PCF. For each, we describe their strengths and limitations.

First, consider a system of agents on a two-dimensional square lattice of size $L_x \times L_y$, lattice step $\Delta$ and with the exclusion property that, at any given time, each lattice site can be occupied by at most one agent.  
If $N$ agents occupy the domain, then the occupancy of the lattice can be represented by a matrix $M$: 
\begin{align}
M_{xy}=
\begin{cases}
&  0 \text{  if lattice site $(x,y)$ is vacant,}\\
&1 \text{  if lattice site $(x,y)$ is occupied,}
\end{cases}
\end{align}
where 
\begin{align}
N=\sum_{x=1}^{L_x}\sum_{y=1}^{L_y} M_{xy} \leq L_xL_y .
\end{align}
Let $\psi^M$ be the set of all agent pairs in the lattice defined by matrix $M$, \textit{i.e.}:
\begin{align}
\begin{split}
\psi^M=&\{(\boldsymbol{a},\boldsymbol{b})\in \mathbb{L}\times\mathbb{L}\hspace{0.6em} |\hspace{0.6em} \boldsymbol{a}=(x_a,y_a),\   \boldsymbol{b}=(x_b,y_b),\  a\neq b,\  M_{x_a,y_a}=M_{x_b,y_b}=1 \} 
\end{split} ,
\end{align}
where $\mathbb{L}=\{1,\dots, L_x\}\times \{1, \dots , L_y\}$ is the set of all sites in the lattice. With agents in configuration $M$, let us define the subset of agent pairs separated by distance $m$ according to some (as yet unspecified) definition of distance, denoted by $d$,  as
\begin{align}
\begin{split}
C_d(m)=&\{(\boldsymbol{a},\boldsymbol{b}) \in \psi^M \ \big|\  \norm{\boldsymbol{a}-\boldsymbol{b}}_{d}=m\} , \text{  for }  m \in \mathcal{D}_d,
\end{split} 
\end{align}
where $\mathcal{D}_d$ is the set of possible distances under the metric $d$.
We define the total number of pairs of agents for each value of distance $m \in \mathcal{D}_d$ as:
\begin{align}
c_d(m)=\left|C_{d}(m)\right|. 
\end{align}
Similarly, we define the set of pairs of sites (regardless of their occupancy) which are separated by distance $m$ according to the metric $d$ as
\begin{align}
\label{eq:S_def_2}
\begin{split}
S_d(m)=&\{(\boldsymbol{a},\boldsymbol{b}) \in \mathbb{L}\times \mathbb{L} \ \big|\  \norm{\boldsymbol{a}-\boldsymbol{b}}_{d}=m\} , \text{  for }  m \in \mathcal{D}_d,
\end{split} 
\end{align}
hence the total number of pairs of sites at distance $m$ is given by
\begin{equation}
	\label{eq:def_s}
	s_d(m)=|S_d(m)|\, .
\end{equation}

To produce a PCF, we aim to normalise the index $c_d(m)$ with the number of pairs we would expect at distance $m$ if the system had no spatial correlation. That is, we consider the case in which the same number of agents are displaced uniformly at random on the lattice and compute the expected number of pairs at distance $m$.
Let $U$ be a random matrix such that $U_{xy}=0$ for all sites $(x,y) \in \mathbb{L}$, apart from $N$ sites chosen uniformly at random without replacement, for which $U_{xy}=1$. Then we define 
\begin{align}
\begin{split}
\bar{C}_d(m)=&\{(\boldsymbol{a},\boldsymbol{b}) \in \psi^U \ \big|\  \norm{\boldsymbol{a}-\boldsymbol{b}}_{d}=m\} , \text{  for }  m \in \mathcal{D}_d.
\end{split} 
\end{align}
Hence, for each value of $m \in \mathcal{D}_d$, the PCF at distance $m$ is defined as:
\begin{equation}\label{eq:PCF_definition}
f_d(m):=\frac{c_d (m)}{\mathbb{E}[\bar{c}_d (m)]}  ,
\end{equation} 
where $\bar{c}_d (m)=|\bar{C}_d(m)|$ and $\mathbb{E}$ represents the expectation operator. 

\subsection{The Annular PCF}
\begin{figure}[h!!]
	\begin{center}
	\includegraphics[width=0.4 \columnwidth]{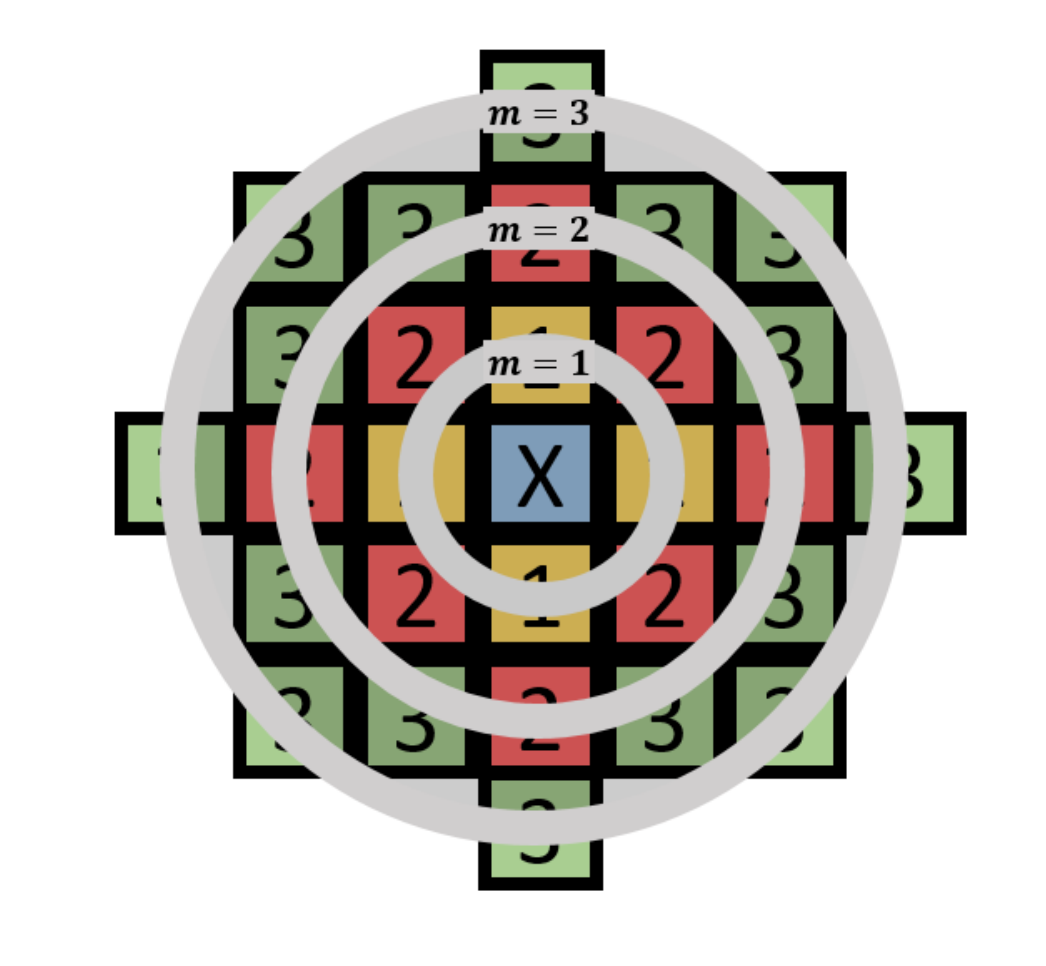}\label{fig:euc1}  
	\end{center}
	\caption{Example agent pairs in the sets $C_A(m)$ with $m=1,2,3$ and bandwidth $\delta=1$. Concentric annuli, $(m-\delta, m]$, are superposed on top of the square lattice and the sites whose centres fall into each annulus are coloured differently. Sites in yellow, red and green (labelled 1, 2 and 3 respectively) are defined to be distance one, two and three from the blue site respectively labelled X.}
	\label{fig:euc} 
\end{figure}

The Annular PCF was originally designed for two-dimensional off-lattice systems \cite{illian2008sam}, but can be extended to on-lattice systems with periodic boundary conditions (BC). For the Annular PCF, the set $C_A$, where $A$ denotes the annular metric, is defined as follows:
\begin{align}
\begin{split}
C_A(m)=&\left\lbrace (\boldsymbol{a},\boldsymbol{b}) \in \psi^M  \big|\  \sqrt{(x_a-x_b)^2+(y_a-y_b)^2} \in \left(m-\delta, m\right]\right\rbrace ,
\end{split} 
\end{align}
where $ m \in \left\lbrace \delta k  \,| \, k \in \mathbb{N}^+ \right\rbrace$ and $\delta$ is a small real number which determines the bandwidth of the PCF. The schematics in Fig. \ref{fig:euc} show a representation of some elements in the sets  $C_A(1)$, $C_A(2)$ and $C_A(3)$ with $\delta =1$.
The normalisation factor is given by
\begin{align}
\mathbb{E}[\bar{c}_A(m)] \approx \frac{ N(N-1) (2 \pi m \delta)} {L_x L_y}\; , \label{eq:norm_ring}
\end{align} 
where $2\pi m \delta $ approximates the area of the $m^{th}$ annulus (assuming small $\delta$) from any given agent. The Annular PCF, $f_A$, follows from definition \eqref{eq:PCF_definition}.

The normalisation in expression \eqref{eq:norm_ring} is a good approximation for a continuous domain with a small $\delta$. However, when the agents are positioned on a lattice, this approach is no longer appropriate. The main issue is that the counts of agents in each annulus vary in an unpredictable manner with the distance, $m$, and the annular width, $\delta$. 
For example, consider a square lattice with spacing $\Delta$. The only possible distances two agents can be separated by are in the countable set
  \begin{align} 
  \mathcal{D}_A = \left\lbrace \Delta\sqrt{x^2+y^2} \, |\, (x,y) \in \mathbb{N}^2 \backslash \lbrace 0,0 \rbrace \right\rbrace= &\lbrace \Delta,\sqrt{2} \Delta,2 \Delta, \sqrt{5} \Delta, \dots \rbrace.
  \end{align}
Partitioning these distances into regularly spaced intervals, as is required by the Euclidean distance metric, we can see that the number of agent pairs does not increase smoothly with the distance, $m$. Depending on the value of $\delta$ it may not even increase monotonically.  However, the definition of the normalisation factor \eqref{eq:norm_ring} suggests that the expected number of pairs increases smoothly and monotonically with both $m$ and $\delta$. This disparity means the on-lattice Annular PCF will not be properly normalised and will either be an over- or under- approximation, making results hard to interpret (see Fig. \ref{fig:ex_unif} \subref{fig:ex_unif_b} as an example).

\subsection{The Rectilinear PCF}
In more recent work, \citet{binder2013qss,binder2015sap} define the Rectilinear PCF specifically for two-dimensional, on-lattice exclusion processes with non-periodic BC. Their definition is easily extendible to periodic BC. They define two PCFs for the two Cartesian directions. In each case the distance is defined by the number of columns (or rows) separating two agents. 

\begin{figure}[h!!]
	\begin{center}
		\subfigure[][]{\includegraphics[width=0.2 \columnwidth]{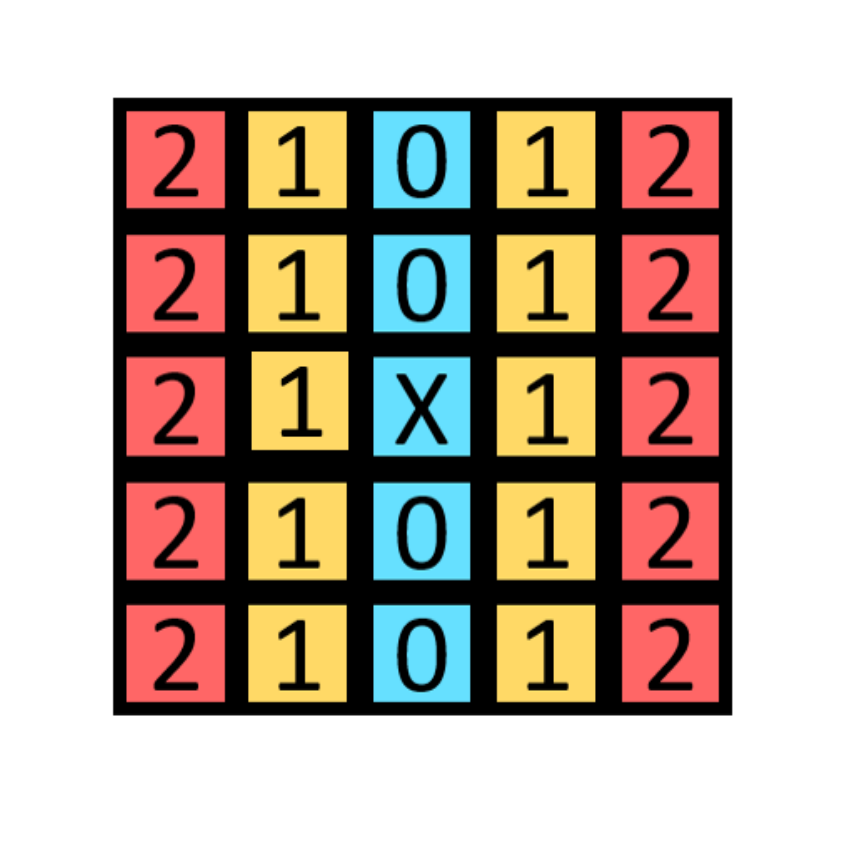}\label{fig:rect1} } \quad
		\subfigure[][]{\includegraphics[width=0.2 \columnwidth]{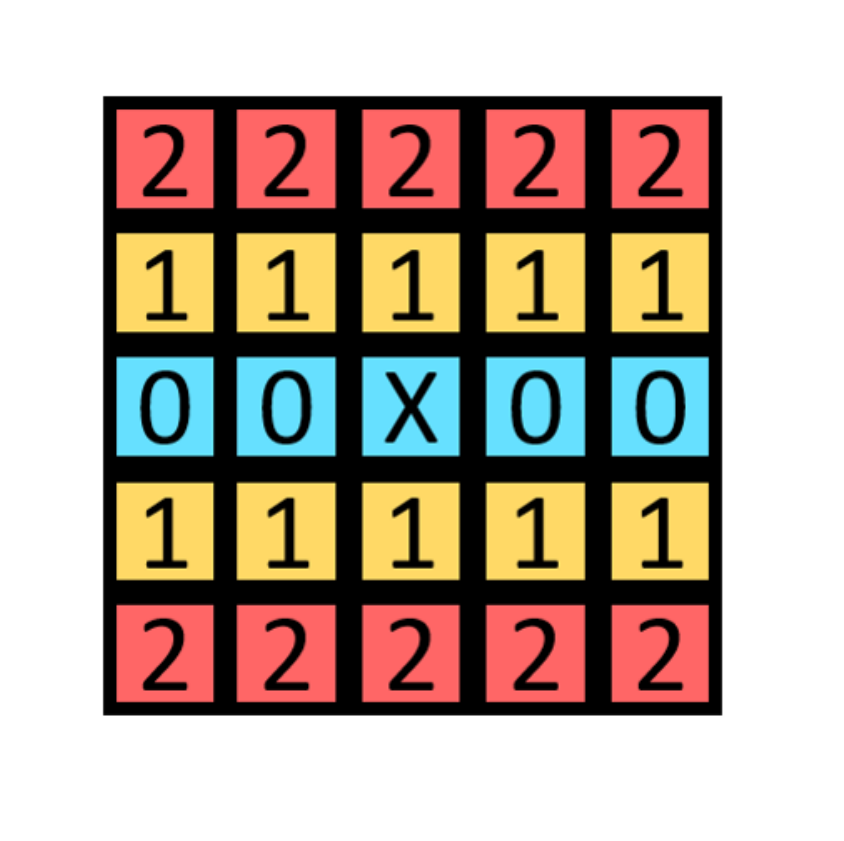}\label{fig:rect2} }
	\end{center}
	\caption{Schematic of agent pairs using (a) the $R_x$ metric, (b) the $R_y$ metric. Sites in yellow and red (labelled 1, 2 respectively) are defined to be distance-one and -two neighbours respectively from the blue site labelled X.}
	\label{fig:rect} 
\end{figure}

Thus, the set of pairs of agents separated by integer distance $m \in \mathbb{N}^+$ are defined in the $x$ direction and $y$ direction respectively as
\begin{subequations}
\begin{align}
C_{R_x}(m)=&\{(\boldsymbol{a},\boldsymbol{b}) \in \psi^M \  \big| \  |x_a-x_b|=m \},
\\C^M_{R_y}(m)=&\{(\boldsymbol{a},\boldsymbol{b}) \in \psi^M \  \big| \  |y_a-y_b|=m \},
\end{align}
\end{subequations}
where subscripts $R_x, R_y$ refer to the metrics defined by the Rectilinear PCFs.
The schematics in Fig. \ref{fig:rect} represent examples of sites separated by distances $m=0$, $m=1$ and $m=2$ for metrics $R_x$ and $R_y$. 

The counts are then normalised by the expected number of pairs of agents at distance $m$ assuming $N$ uniformly distributed agents:
\begin{subequations}
\begin{align}
\mathbb{E}[\bar{c}_{R_x} (m)] = &\frac{N-1}{L_xL_y-1}\frac{N}{L_xL_y}L_y^2(L_x-m),
\\\mathbb{E}[\bar{c}_{R_y} (m)] = &\frac{N-1}{L_xL_y-1}\frac{N}{L_xL_y}L_x^2(L_y-m),
\end{align}
\end{subequations}
respectively. For details of the derivation of these factors see \citep{binder2013qss}.
 The final Rectilinear PCF is defined as the arithmetic average of the two orthogonal PCFs, \textit{i.e.}
\begin{align}
f^M_R(m)=\frac{1}{2} \bigg[f_{R_x}(m)+f_{R_y}(m)\bigg] \hspace{1em} \text{where} \hspace{1em} f^M_{R_x}(m)=& \frac{c_{R_x}}{\mathbb{E}[\bar{c}_{R_x} (m)]}, \, f^M_{R_y}=\frac{c_{R_y}}{\mathbb{E}[\bar{c}_{R_y} (m)]}.
\end{align} 
 The Rectilinear PCF correctly identifies spatial correlation in many examples and, unlike the Annular PCF, is normalised correctly. However, one major issue that the Rectilinear PCF suffers from is that, due to the inherent anisotropy of its definition, spatial structures which are biased in either Cartesian direction may be missed. For such patterns, the PCF given by the $R_x$ and $R_y$ metrics are approximately constant functions of distance, because the averaged row and column densities are constant along the axes, despite the fact that clustering can still be present. Examples of these spatial patterns include many biologically and chemically relevant cases such as diagonal stripes and chessboard patterns \cite{ouyang1991tuh,bhide1999dsc} (see Fig. \ref{fig:ex_scs}). We note that when the pattern structure is biased in only one Cartesian direction, the preaveraging Rectilinear PCFs $f_{R_x}$ and $f_{R_y}$ will identify further information about the direction of the spatial pattern. Another limitation of the Rectilinear PCF is that it applies only to regular square lattices and a generalisation to other forms of tessellations would be challenging. 

\section{The Square Taxicab and Square Uniform PCFs}
\label{sec:PCF}
\begin{figure}[h!!]
\begin{center}
\subfigure[][]{\includegraphics[width=0.25 \columnwidth]{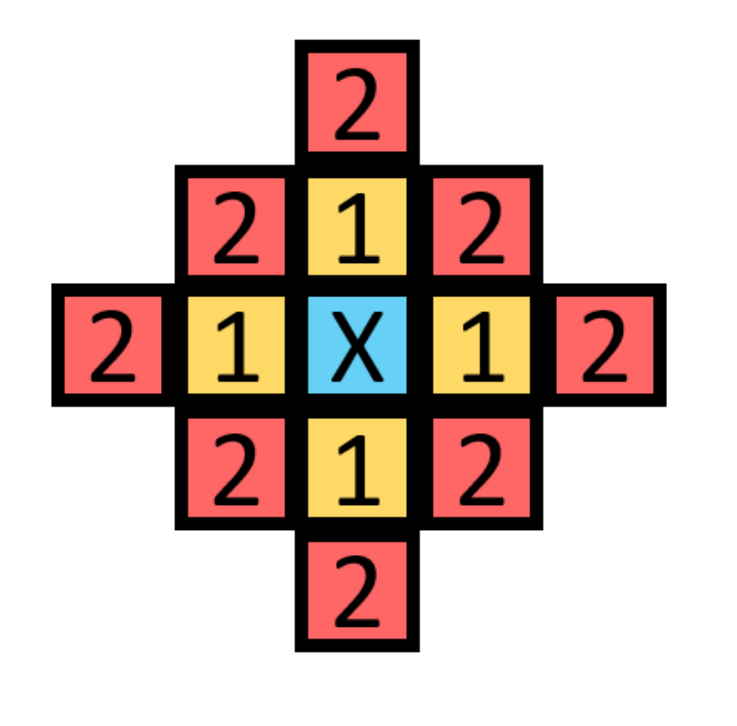}\label{fig:normsvis_TA} }
\subfigure[][]{\includegraphics[width=0.25 \columnwidth]{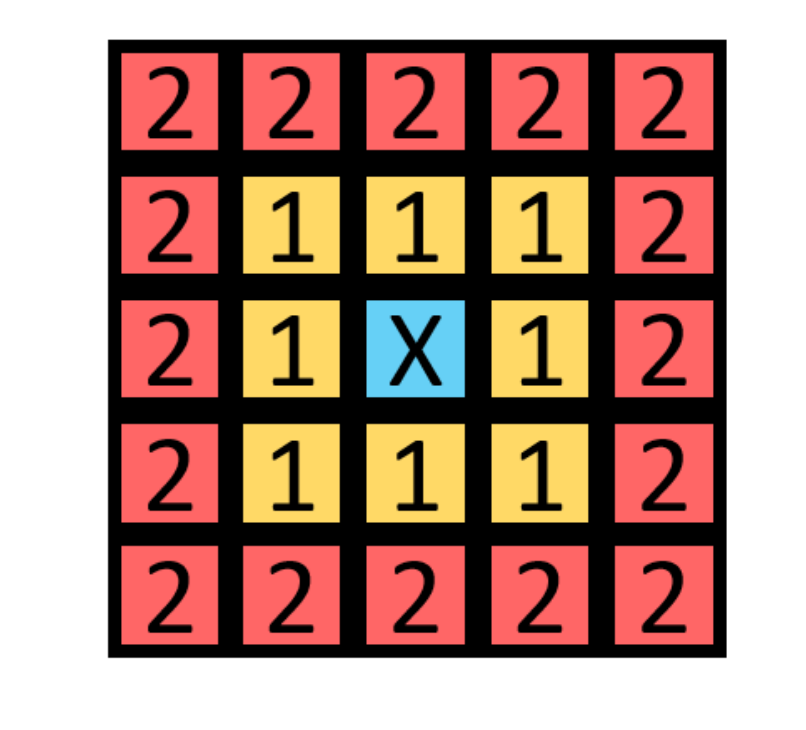}\label{fig:normsvis_UN} }
\end{center}
\caption{Schematic of agent pairs using (a) the taxicab metric (b) the uniform metric. Sites in yellow and red, labelled with numbers 1 and 2 respectively, are defined to be distance-one and -two neighbours from the site marked in blue (labelled with X).}
\label{fig:normsvis} 
\end{figure}

In this section we define two new discrete PCFs for a square lattice: the Square Taxicab and Square Uniform PCF, using the taxicab and uniform metric, respectively, under both periodic and non-periodic BC.
Using the same notation as in Section \ref{sec:previous}, we define the subsets of agent pairs separated by distance $m$ under non-periodic BC as
\begin{subequations}
\begin{align}
\label{eq:def_bigS_nonper}
C^{n}_1(m)=&\{(\boldsymbol{a},\boldsymbol{b}) \in \psi^M \ \big|\  \norm{\boldsymbol{a}-\boldsymbol{b}}_{1}=m\}, \quad  \, m \in \mathcal{D}^n_1,
\\C_{\infty}^{n}(m)=&\{(\boldsymbol{a},\boldsymbol{b}) \in \psi^M \  \big|\  \norm{\boldsymbol{a}-\boldsymbol{b}}_\infty=m \}, \quad m \in \mathcal{D}^n_\infty,
\end{align}
\end{subequations}
for the taxicab and uniform metric respectively. Here, $\mathcal{D}^n_1 = \mathcal{D}^n_\infty= \big\lbrace 1,2,...\min \lbrace L_x,L_y \rbrace -1 \big\rbrace$ and the superscript $n$ refers to the fact we are considering \textit{non-periodic} BC. Using the definitions of the uniform and taxicab metrics, we can express these sets as: 
\begin{subequations}
\begin{align}
C^{n}_1(m)=&\{(\boldsymbol{a},\boldsymbol{b}) \in \psi^M \  \big|\  |x_a-x_b|+|y_a-y_b|=m\} ,
\\C_{\infty}^{n}(m)=&\{(\boldsymbol{a},\boldsymbol{b}) \in \psi^M \  \big |\  \max \lbrace |x_a-x_b|\text{, }|y_a-y_b| \rbrace=m  \} .
\end{align}
\end{subequations}
Similarly we define the subsets of agent 
pairs separated by distance $m$ under periodic BC as

\begin{subequations}	
\begin{align}
\label{eq:def_bigS_per}	
&\begin{aligned}
C^{p}_1(m)= \Big\lbrace(\boldsymbol{a},\boldsymbol{b})\in  \psi^M \mid \min \lbrace & |x_a-x_b|, L_x-|x_a-x_b| \rbrace \\&+ \min \lbrace (|y_a-y_b|, L_y-|y_a-y_b|\rbrace =m\Big\rbrace , \quad  \, m \in \mathcal{D}^p_1,
\end{aligned}\\
\label{eq:def_S_p_unif}
&\begin{aligned}
C^{p}_\infty(m)=\Big\lbrace (\boldsymbol{a},\boldsymbol{b})\in \psi^M \mid \max \big\lbrace & \min \lbrace |x_a-x_b|, L_x-|x_a-x_b| \rbrace, \\& \min\lbrace (|y_a-y_b|, L_y-|y_a-y_b|\rbrace \big\rbrace =m \Big\rbrace , \quad  \, m \in \mathcal{D}^p_\infty,
\end{aligned}
\end{align}
\end{subequations}
where $\mathcal{D}^p_1 = \mathcal{D}^p_\infty= \big\lbrace 1,2,...\min \big\lbrace \big\lfloor\frac{L_x}{2}\big\rfloor,\big\lfloor\frac{L_y}{2}\big\rfloor \big\rbrace \big\rbrace$. The corresponding definitions of $S_1^n$, $S_\infty^n$, $S_1^p$ and $S_\infty^p$ can be obtained similarly by using equation \eqref{eq:S_def_2}. Here the superscript $p$ refers to the fact we are considering \textit{periodic} BC. Notice that we restrict the largest $m$ to be $\min \big\lbrace \big\lfloor\frac{L_x}{2}\big\rfloor,\big\lfloor\frac{L_y}{2}\big\rfloor \big\rbrace$ in order to simplify the computation of the normalisation factor (see Section \ref{sec:norm_periodic2}). However, with some work this restriction could be relaxed.
The schematics in Fig.  \ref{fig:normsvis} represent examples of sites separated by distance $m=1$ and $m=2$ using the taxicab (a) and uniform (b) metrics. 

For the normalisation factors, \citet{binder2013qss} use:
\begin{align}
\label{normalisation_factor_1}
&\mathbb{E}[\bar{c}_{d} (m)]=\left(\frac{N}{L_x L_y}\right)\left(\frac{N-1}{L_x L_y-1}\right) s_d (m),
\end{align} where $s_d(m)$ is defined as in equation \eqref{eq:def_s}
and $d$ refers to the metric used.
In other words, the expected number of pairs of agents at distance $m$ on a lattice with $N$ uniformly distributed agents can be written as the probability that two different sites at distance $m$ are simultaneously occupied, multiplied by the total number of pairs of sites at distance $m$ in the domain.

To complete the definitions of our square lattice PCFs we need to provide an expression for $s_1$ and $s_\infty$. We address this in the next two sections, distinguishing between the cases of periodic and non-periodic BC.

\subsection{Normalisation of the Square Taxicab and Square Uniform PCF under periodic boundary conditions}
\label{sec:norm_periodic2}
As the derivation for the normalisation is simple under periodic BC and more complicated under non-periodic BC, we first consider the system with periodic BC and determine $s^{p}_{1}$, $s^{p}_\infty$, where $p$ denotes periodic BC.
Let us define the number of sites separated by a distance $m$ from any given reference site on a lattice as $t_1(m)$, $t_\infty(m)$ under the taxicab and uniform metric respectively. These read
 \begin{subequations}
 \label{eq:t_2D}
 \begin{align}
 \label{eq:t_2Da}
 &t_1(m)=4m,
 \\&t_\infty(m)=8m. \label{eq:t_2Db}
 \end{align}
 \end{subequations}
 The proofs of the expressions \eqref{eq:t_2D} are omitted, but they can be obtained easily by induction on $m$. Examples for $m=1,2$ can be seen in Fig. \ref{fig:normsvis}.
Notice that for $m\le\min \big\lbrace \big\lfloor\frac{L_x}{2}\big\rfloor,\big\lfloor\frac{L_y}{2}\big\rfloor \big\rbrace$, given any site on the lattice, the number of sites at distance $m$ from this reference site in the case of periodic BC is exactly $t(m)$. Consider the lattice of size $L_x \times L_y$ with $L_x, L_y >2$. If we multiply the total number of lattice sites by $t(m)$, we count each pair of sites separated by distance $m$ exactly twice. Hence we conclude that: 
  \begin{align}
   \label{eq:c^A(m)}
 s_{d=1,\infty}^{p}(m)=\frac{t(m)L_xL_y}{2},
 \end{align}
 using the taxicab and uniform metrics. Substituting values for $t_1(m)$ and $t_\infty(m)$ from equations \eqref{eq:t_2D} we deduce that
  	
 \begin{subequations} 
 \label{cA}
  \begin{align}
 s^{p}_1(m)=2mL_xL_y,
 \\s^{p}_\infty(m)=4mL_xL_y.
 \label{eq:def_s_unif}
 \end{align}
 \end{subequations}
 Therefore, by substituting expressions \eqref{cA} into equation \eqref{normalisation_factor_1}, the normalisation factors under periodic BC are:
 
 \begin{subequations}
 \begin{align}
 \mathbb{E}\left[ \bar{c}^{p}_1(m)\right]=\frac{2mN(N-1)}{L_xL_y-1},
 \\\mathbb{E}\left[ \bar{c}^{p}_\infty(m)\right]=\frac{4mN(N-1)}{L_xL_y-1}.
 \end{align}
 \end{subequations}
 
\subsection{Normalisation of the Square Taxicab and Square Uniform PCF under non-periodic BC}\label{sec:norm_non_per}
In this section we derive expressions for $s_1^{n}(m)$ and $s_{\infty}^{n}(m)$, where $n$ denotes non-periodic BC. Notice that, for all $m \in \mathcal{D}$, we have that $s^{p}(m) > s^{n}(m)$ since $s^{p}(m)$ includes pairs that cross the domain boundary whereas $s^{n}(m)$ does not. Therefore, in order to find a formula for $s^{n}(m)$, it is enough to determine a formula for the remainders defined by:
  \begin{subequations}
  \begin{align}
\label{eq:def_rem}
r_1(m)=&s^{p}_1(m)-s^{n}_1(m),
\\r_\infty(m)=&s^{p}_\infty(m)-s^{n}_\infty(m).
\end{align} 	
  \end{subequations}
The remainders count the number of pairs of sites which cross a boundary under the periodic BC. For simplicity, throughout this section, we will only derive the normalisation for the taxicab metric, however, the derivation for the uniform metric is similar and can be found in SM Section \ref{SUPP-sec:norm_unif} for reference.
Let us define the set of pairs of sites separated by distance $m \in \mathcal{D}^n_1$ that cross the $x$ boundary (horizontal axis) or $y$ boundary (vertical axis), respectively, as
\begin{subequations}
\begin{align}
P_1^x(m)=&\big\lbrace (\boldsymbol{a}, \boldsymbol{b}) \in S_1^p(m)\, \  \big| \ |y_a-y_b|>L_y-|y_a-y_b| \big\rbrace, 
\\ P_1^y(m)=&\big\lbrace (\boldsymbol{a}, \boldsymbol{b}) \in S_1^p(m)\, \  \big| \ |x_a-x_b|>L_x-|x_a-x_b| \big\rbrace. 
\end{align}
\end{subequations}
Of these pairs, let us consider those pairs that are at distance $k \in \lbrace 1,... m \rbrace$ rows or columns from each other respectively. We define these subsets as:
\begin{subequations}
\begin{align}
P^x_1(m,k)=&\big\lbrace (\boldsymbol{a}, \boldsymbol{b}) \in \, P_1^x(m) \big| \quad L_y-|y_a-y_b|=k \big\rbrace,
\\P^y_1(m,k)=&\big\lbrace (\boldsymbol{a}, \boldsymbol{b}) \in \, P_1^y(m) \big| \quad L_x-|x_a-x_b|=k \big\rbrace . 
\end{align}
\end{subequations}
\begin{figure}[h!!!!!!!!]
	\begin{center}
	\subfigure[][]{\includegraphics[height=0.34 \columnwidth]{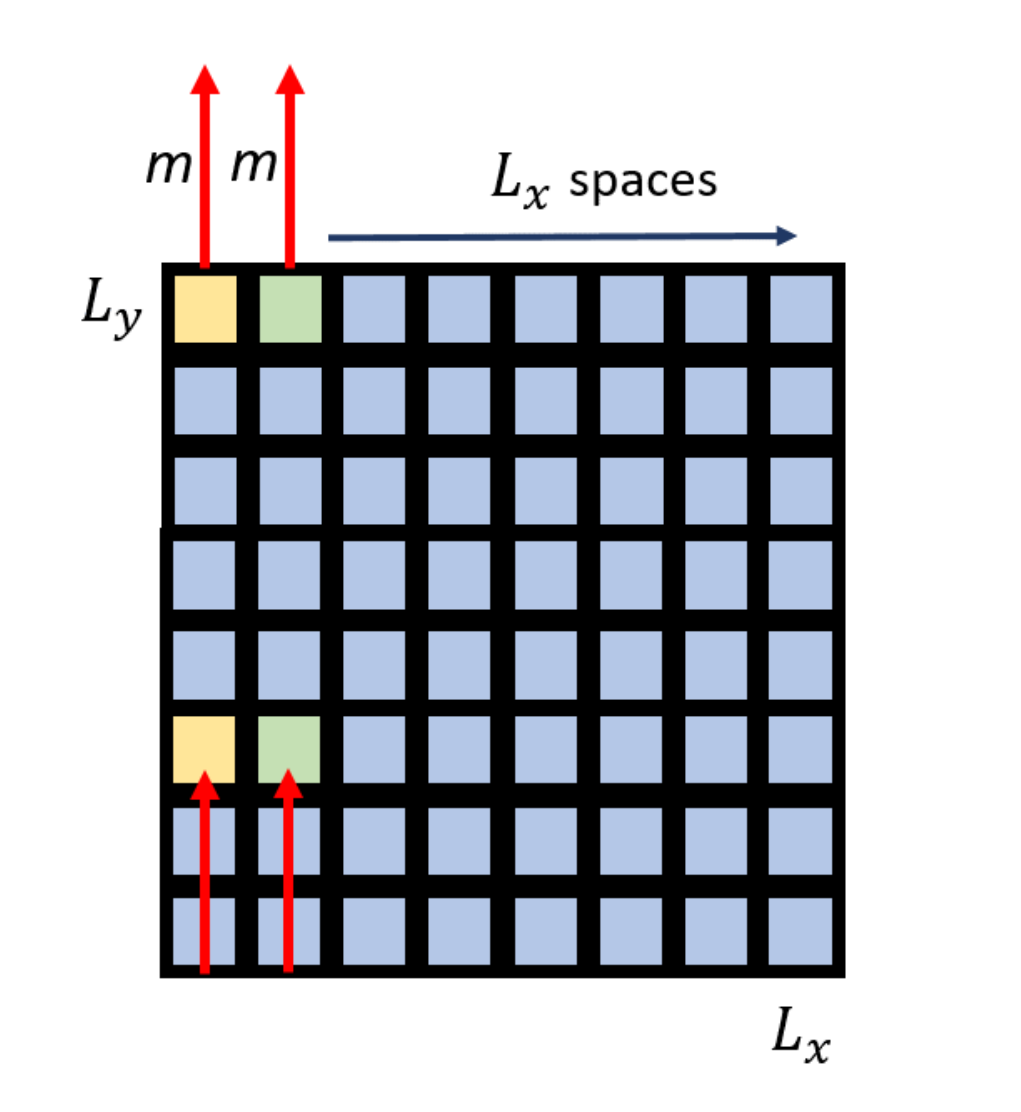}\label{figure:pairs_a}}
	\subfigure[][]{\includegraphics[height=0.34 \columnwidth]{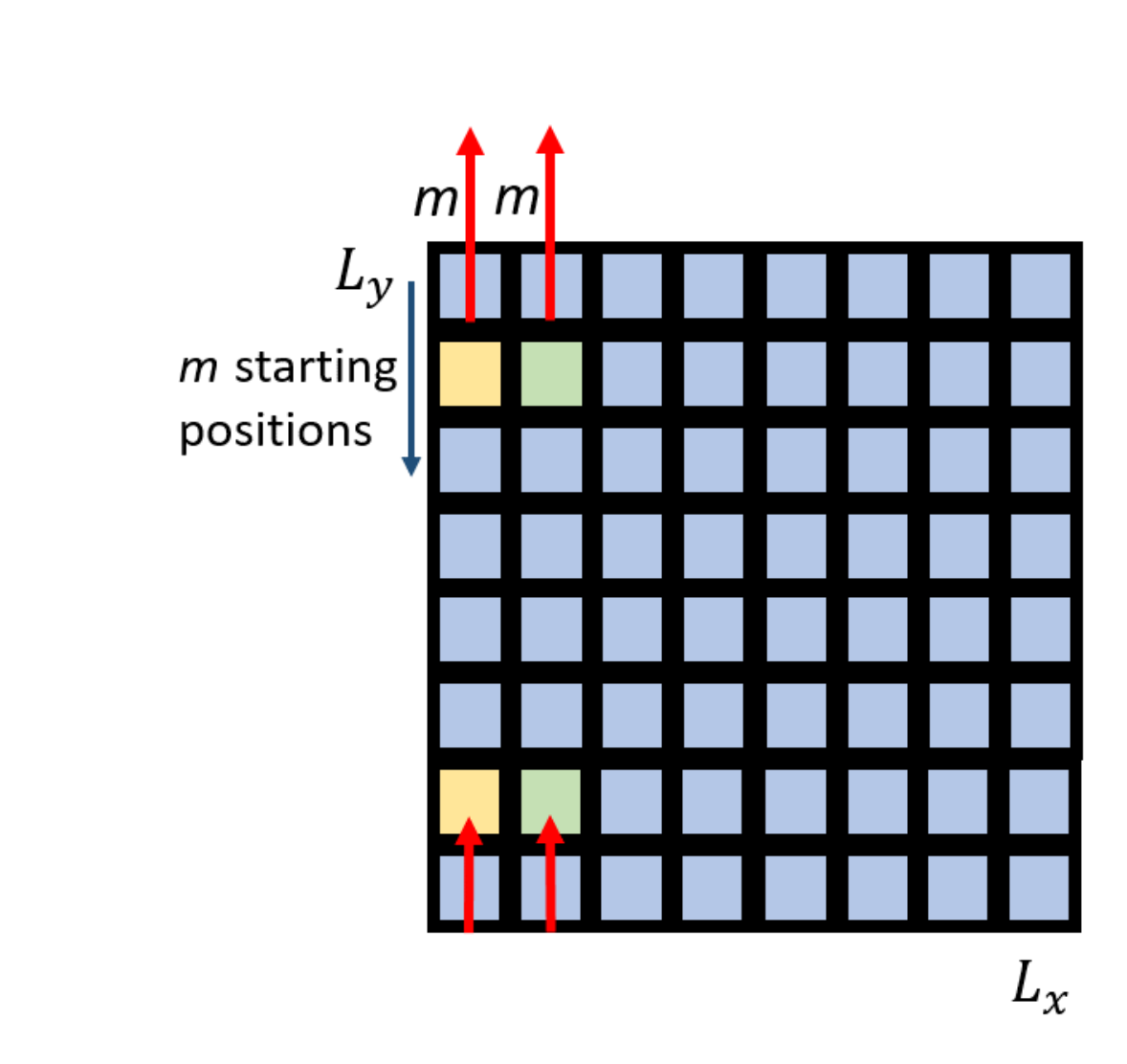}\label{figure:pairs_b}}\\
	\subfigure[][]{\includegraphics[height=0.2 \columnwidth]{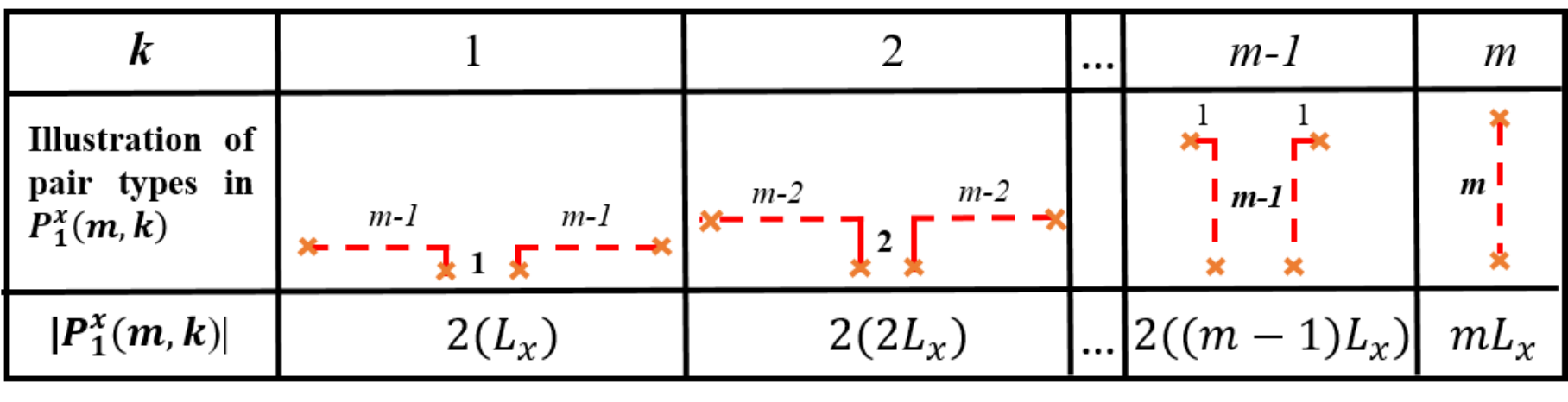} \label{figure:pairs_list}}
\end{center}
	\caption{A visualisation of the pairs of sites in $P^x_1(m)$. Panels (a) and (b) show two different site pairs in $P^x_1(m,m)$. For each of the $L_x$ columns, each site in the rows $\lbrace L_y-m+1, L_y-m+2, \dots , L_y \rbrace$ has a single corresponding site at distance $m$ separated by $m$ rows and 0 columns and reached by crossing the horizontal boundary. Therefore $|P^x_1(m,m)|=mL_x$. Panel (c) shows all the different types of pairs in $P^x_1(m,k)$ for $k=1, \dots , m$ with the corresponding value of $|P^x_1(m,k)|$.} \label{fig:periodic_pairs}
\end{figure}
Notice that $P_1^x(m)=\bigcup_{k =1}^m P_1^x(m,k)$ and $P_1^y(m)=\bigcup_{k =1}^m P_1^y(m,k)$.
Fig. \ref{fig:periodic_pairs} \subref{figure:pairs_a} and \subref{figure:pairs_b} give visualisations of pairs of sites within $P_1^x(m,m)$. Fig. \ref{fig:periodic_pairs} \subref{figure:pairs_list} gives examples of distances between pairs of sites in $P_1^x(m,k)$, for $k=1,\dots,m$.
By definitions \eqref{eq:def_bigS_nonper} and \eqref{eq:def_bigS_per} we have that
\begin{align}
S_1^p(m) \, \backslash \, S_1^n(m)=P^{x}_1(m) \cup P^{y}_1(m).
\end{align}
Hence, by combining equations \eqref{eq:def_rem} and \eqref{eq:def_s}, we obtain
\begin{align}
r_1(m)&=|P^{x}_1(m) \cup P^{y}_1(m)| \nonumber
\\&=\sum_{k=1}^m | P_1^x(m,k) | + \sum_{k=1}^m | P_1^y(m,k) | - |P_1^x(m) \cap P_1^y(m)| \label{eq:rem}	\, .
\end{align}
To conclude the computation we derive an expression for the two sums in equation \eqref{eq:rem} and the corresponding equation for the size of the intersection. 
By counting the contribution of each type of pair (see Fig. \ref{fig:periodic_pairs} for a visualisation), one can write down the following expressions for the two sums in equation \eqref{eq:rem}: 
\begin{subequations}
\begin{align}
	\sum_{k=1}^m | P_1^x(m,k) | &= 2(L_x + 2L_x +\dots +L_x(m-1)) + L_x m \;,
	\\\sum_{k=1}^m | P_1^y(m,k) | &= 2(L_y + 2L_y +\dots +L_y(m-1)) + L_y m\; .
\end{align}
\end{subequations}
Hence
 \begin{align}
\sum_{k=1}^m | P_1^x(m,k) | + \sum_{k=1}^m | P_1^y(m,k) |&=2(L_x + 2L_x +\dots + L_x(m-1)) + L_x m \nonumber \\&+2(L_y + 2L_y +\dots + L_y(m-1)) + L_y m\nonumber
\\&=(L_x+L_y)\bigg(m+2\sum_{i=1}^{m-1}i\bigg) \nonumber
\\&=(L_x+L_y)\bigg(m+2\frac{(m-1)m}{2}\bigg) \nonumber
\\&=(L_x+L_y)m^2 \, .
\label{eq:diff_tc}
\end{align}
We now focus on deriving an expression for the size of intersection, $|P_1^x(m) \cap P_1^y(m)| $, in equation \eqref{eq:rem}. 
The set $P_1^x(m) \cap P_1^y(m) $ consists of pairs of sites separated by distance $m$ that cross both the $x$ and $y$ boundaries simultaneously. 
\begin{figure}[t]
	\centering
	\includegraphics[width=0.6 \columnwidth]{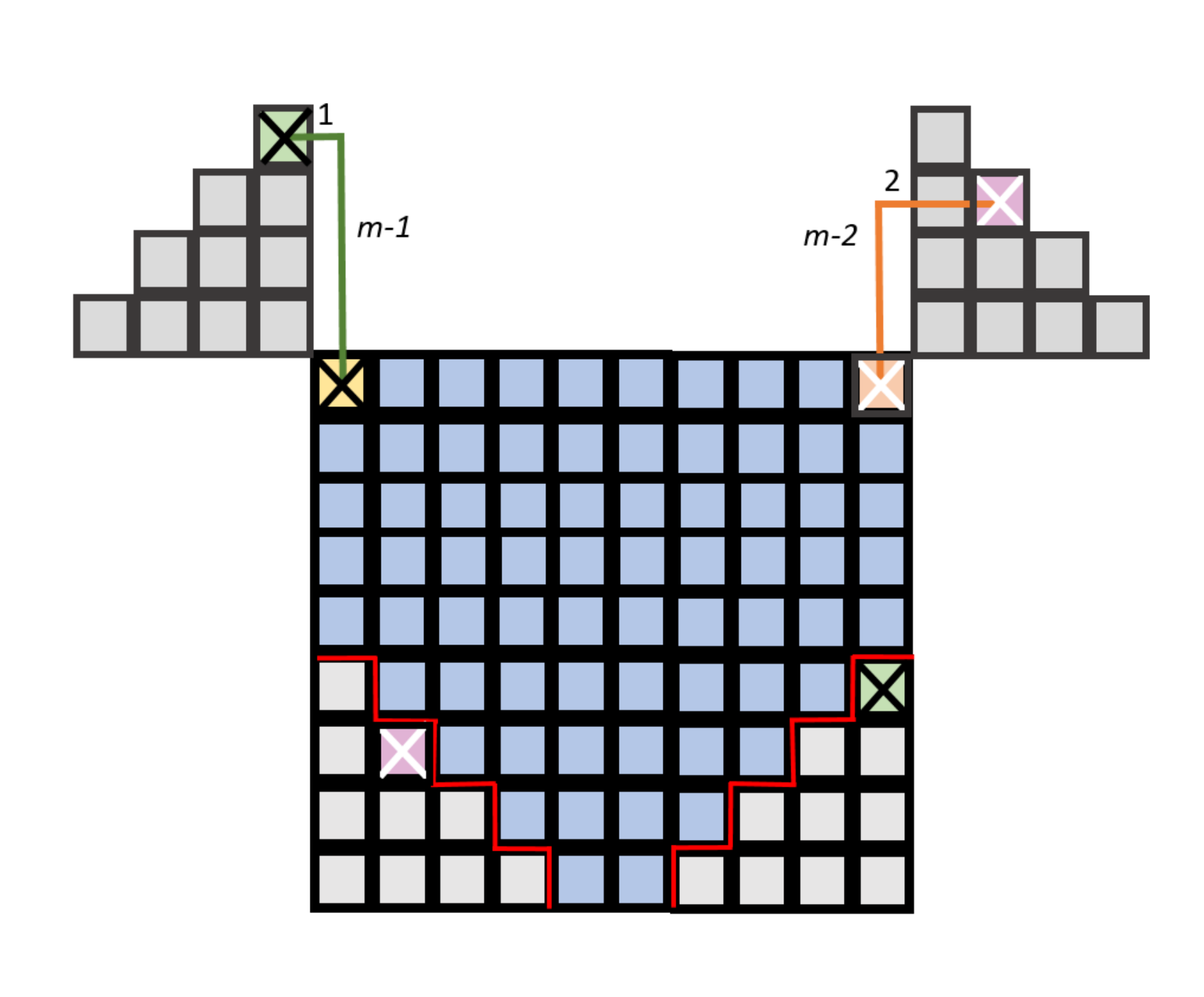}
	\caption{Examples of pairs of sites separated by distance $m=5$ that cross both the $x$ and $y$ boundaries, \textit{i.e.} pairs of sites in $P_1^x(m) \cap P_1^y(m)$. The grey sites outside the domain correspond to the grey sites inside the domain in the diametrically opposite corner. As illustrations of site pairs at a distance $m$ which cross both boundaries, the orange site containing a white cross is distance $m$ from the pink site containing a white cross. Similarly, the yellow site containing a black cross is distance $m$ from the green site containing a black cross.} \label{fig:example_corner}
\end{figure}

\begin{figure}
\begin{center}
		\includegraphics[width=0.4 \columnwidth]{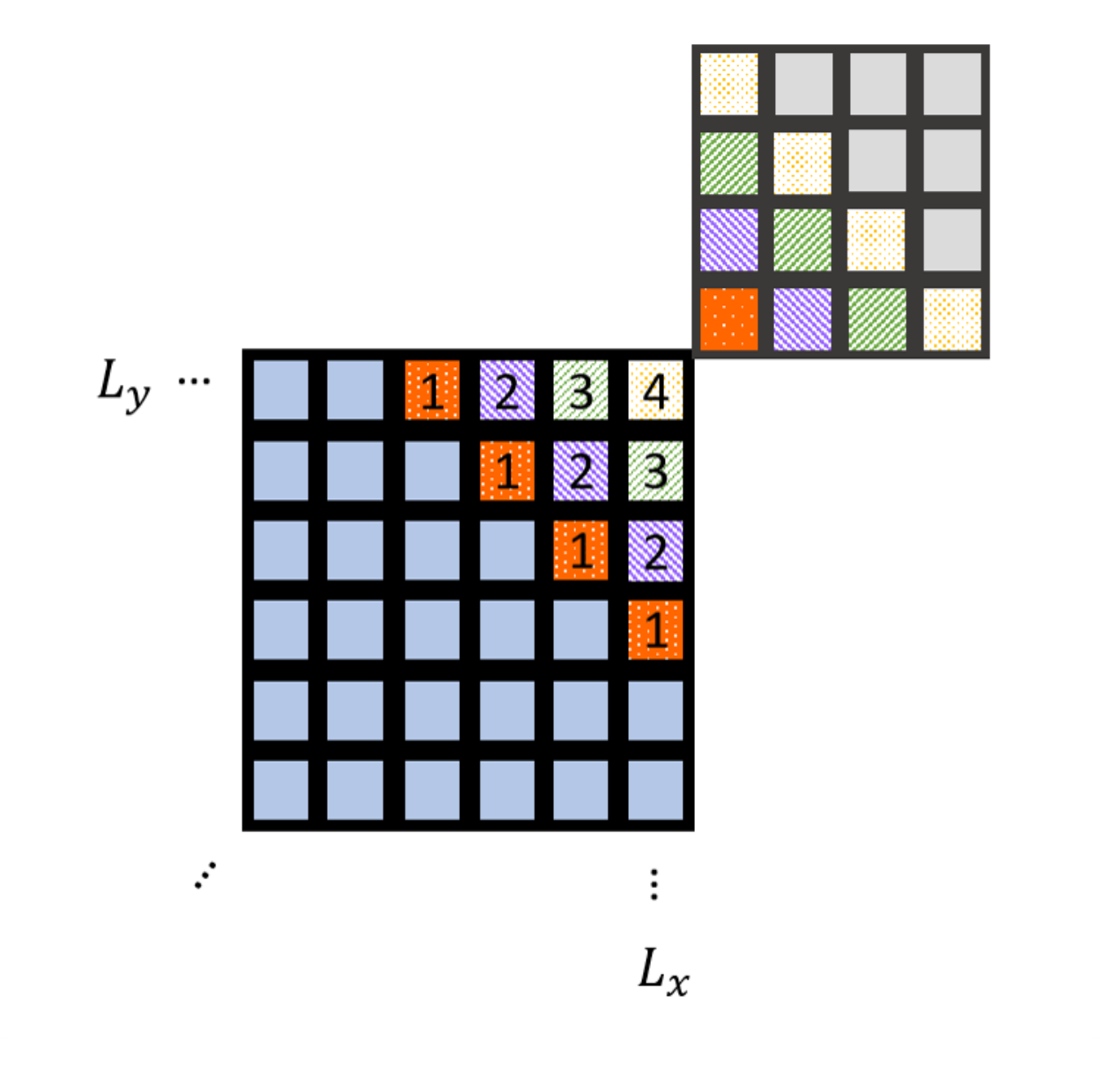}
		\caption{Examples of pairs in $P_1^x(5) \cap P_1^y(5)$ on a zoomed-in corner of a larger domain. Sites with a given pattern (not plain) are distance $m=5$ away from sites with the same pattern that can be reached by crossing the $x$ and $y$ boundaries. The number in each given site corresponds to the number of sites at distance $m=5$, reached by crossing the $x$ and $y$ boundaries.}
		\label{fig:cross_two_axis}
\end{center}
\end{figure}

There are two regions of the domain where site pairs cross two boundaries. These are any two consecutive corners of the four corners of the domain. Examples of these regions and site pairs within these regions are visualised in Fig. \ref{fig:example_corner}. 
In Fig. \ref{fig:cross_two_axis} we give an illustrative example in which we count the number of these pairs for $m=5$. All sites inside the boundaries of the domain coloured in orange, purple, green and yellow are distance $m=5$ from other sites of the same colour outside the boundaries of the domain. Notice that, the yellow site in the corner at $(L_x,L_y)$ is distance five from a total of four sites, reached by crossing the $x$ and $y$ boundaries, denoted by a 4 in the site. Similarly, the two green sites at $(L_x-1,L_y)$ and $(L_x,L_y-1)$ are distance $m$ from three sites, reached by crossing the $x$ and $y$ boundary, denoted by a 3 in the two sites. 
 $|P_1^x(5) \cap P_1^y(5)|$ is the sum of all the numbers in the coloured sites multiplied by two to account for the second corner region.
Extrapolating, for any value of $m$, the number of pairs of sites that cross the two boundaries is exactly
\begin{align}
|P_1^x(m) \cap P_1^y(m)|&=2\big((m-1)+2(m-2)+3(m-3)+\dots +m-1 \big) \; \nonumber
\\&=2\sum_{i=1}^{m-1} (m-i) i\nonumber
\\&=\frac{m^3-m}{3}.
\label{eq:taxicab_rem}
\end{align}
By substituting equations \eqref{eq:taxicab_rem} and \eqref{eq:diff_tc} into equation \eqref{eq:rem} we gain an expression for the remainder $r_1(m)$. By rearranging equation \eqref{eq:def_rem} we determine $s_1^n(m)$, which we then substitute into equation \eqref{normalisation_factor_1} to obtain the exact expression for the normalisation in the non-periodic case. This is given by
\begin{align}
\mathbb{E}\left[\bar{c}^{n}_1(m)\right]=&\bigg(\frac{N}{L_xL_y}\bigg)\bigg(\frac{N-1}{L_xL_y-1}\bigg)\bigg(2mL_xL_y-(L_x+L_y)m^2+\frac{m^3-m}{3}\bigg).
\end{align}

A similar approach can be used to obtain the normalisation factor for the uniform metric under non-periodic BC:
\begin{align}
\mathbb{E}\left[\bar{c}^{n}_\infty(m)\right] =&\bigg(\frac{N}{L_xL_y}\bigg)\bigg(\frac{N-1}{L_xL_y-1}\bigg)\bigg(4mL_xL_y-3(L_x+L_y)m^2+2m^3\bigg). \label{eq:unif_non}
\end{align}
For more details on the derivation of expression \eqref{eq:unif_non} see SM Section \ref{SUPP-sec:norm_unif}.

\section{Results}
\label{sec:results}
In this section we use the Square Taxicab and Square Uniform PCFs defined in Section \ref{sec:PCF} to analyse the spatial correlation in some examples. We compare our results with previously suggested on-lattice PCFs. 

\begin{figure}[h!!]
	\begin{center}
		\subfigure[][]{\includegraphics[ width=0.32 \columnwidth]{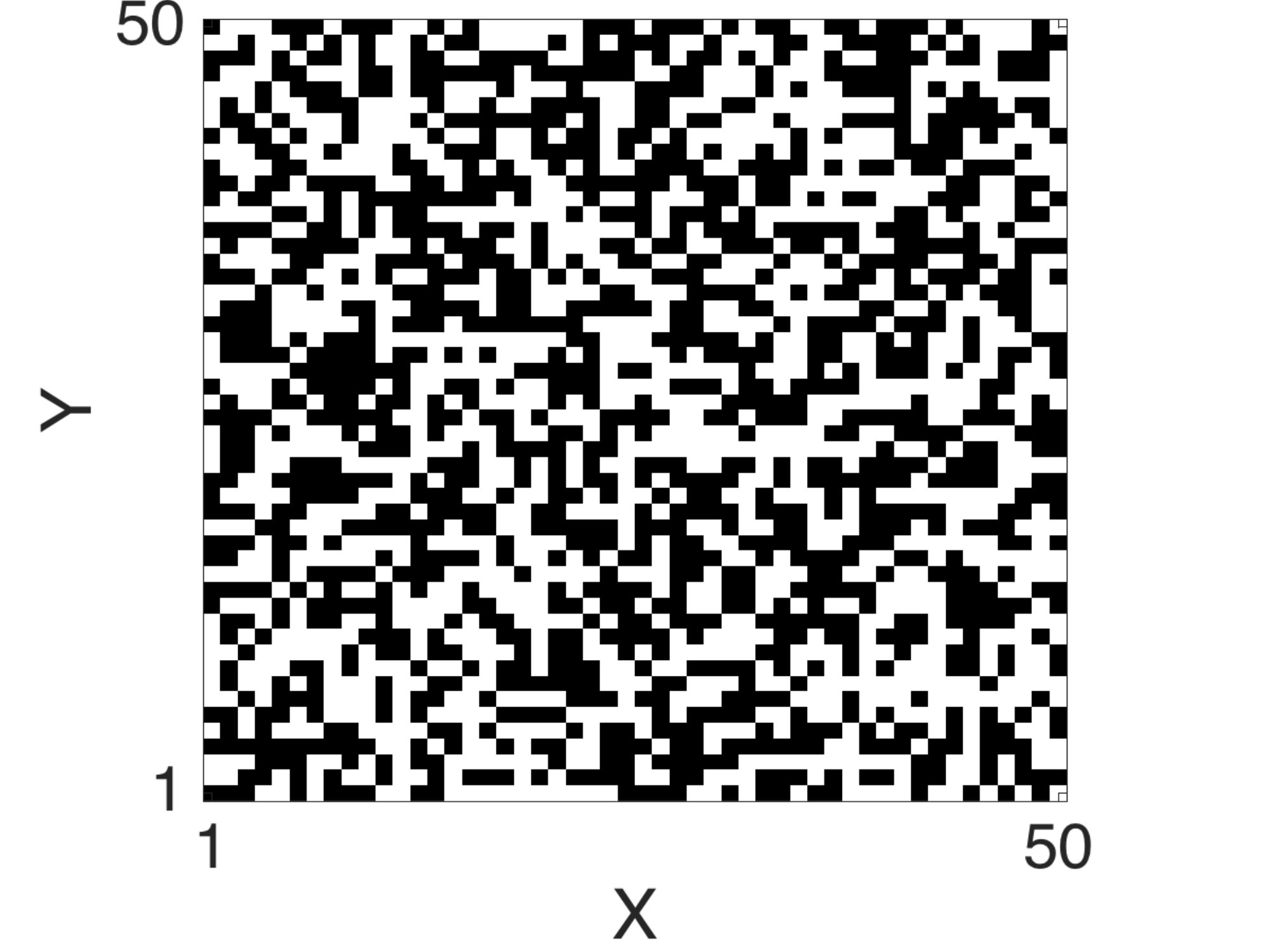}
			\label{fig:ex_unif_a} }
		\subfigure[][]{\includegraphics[width=0.32 \columnwidth]{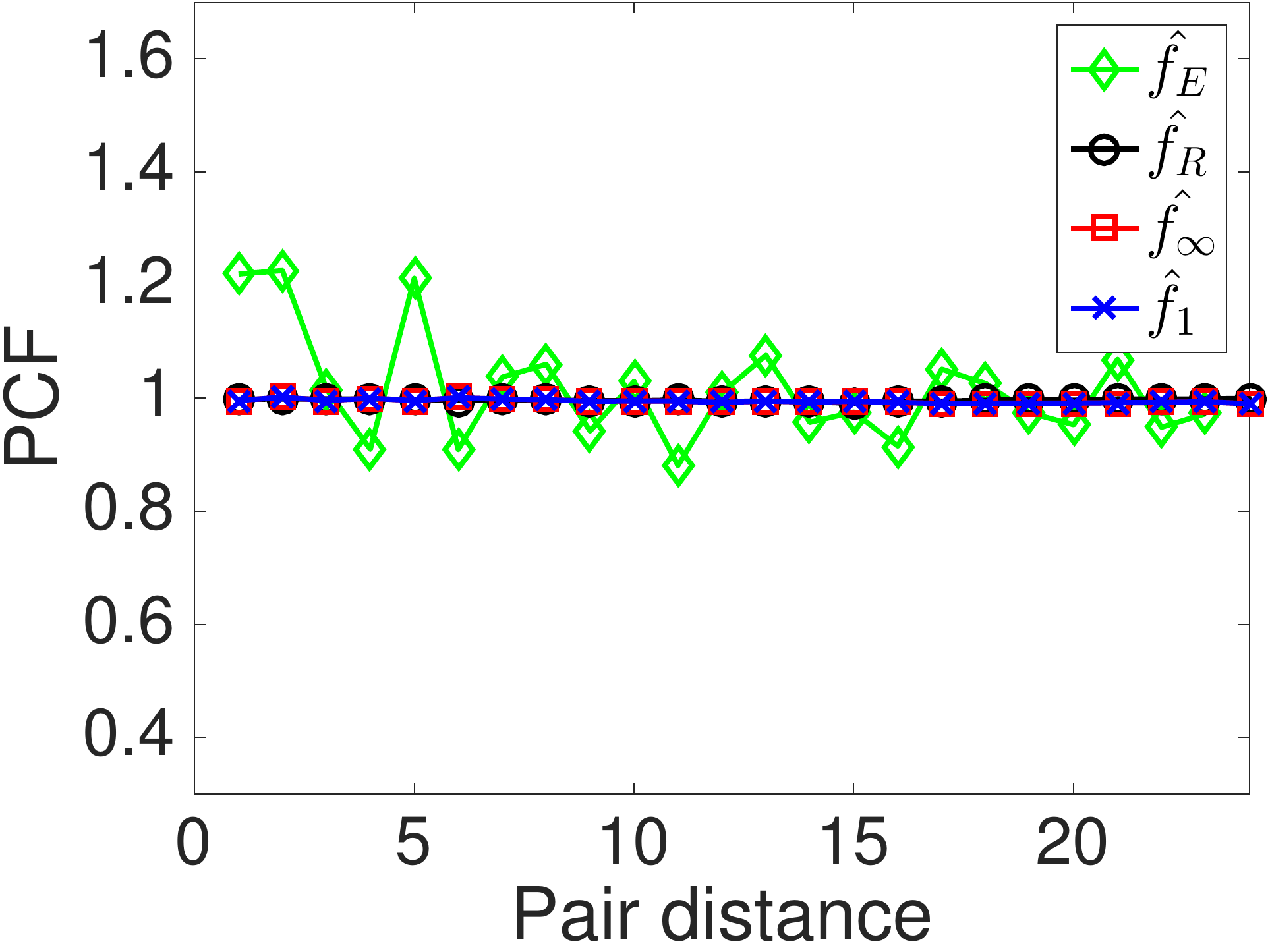}
			\label{fig:ex_unif_b} }\\
		\subfigure[][]{\includegraphics[ width=0.32 \columnwidth]{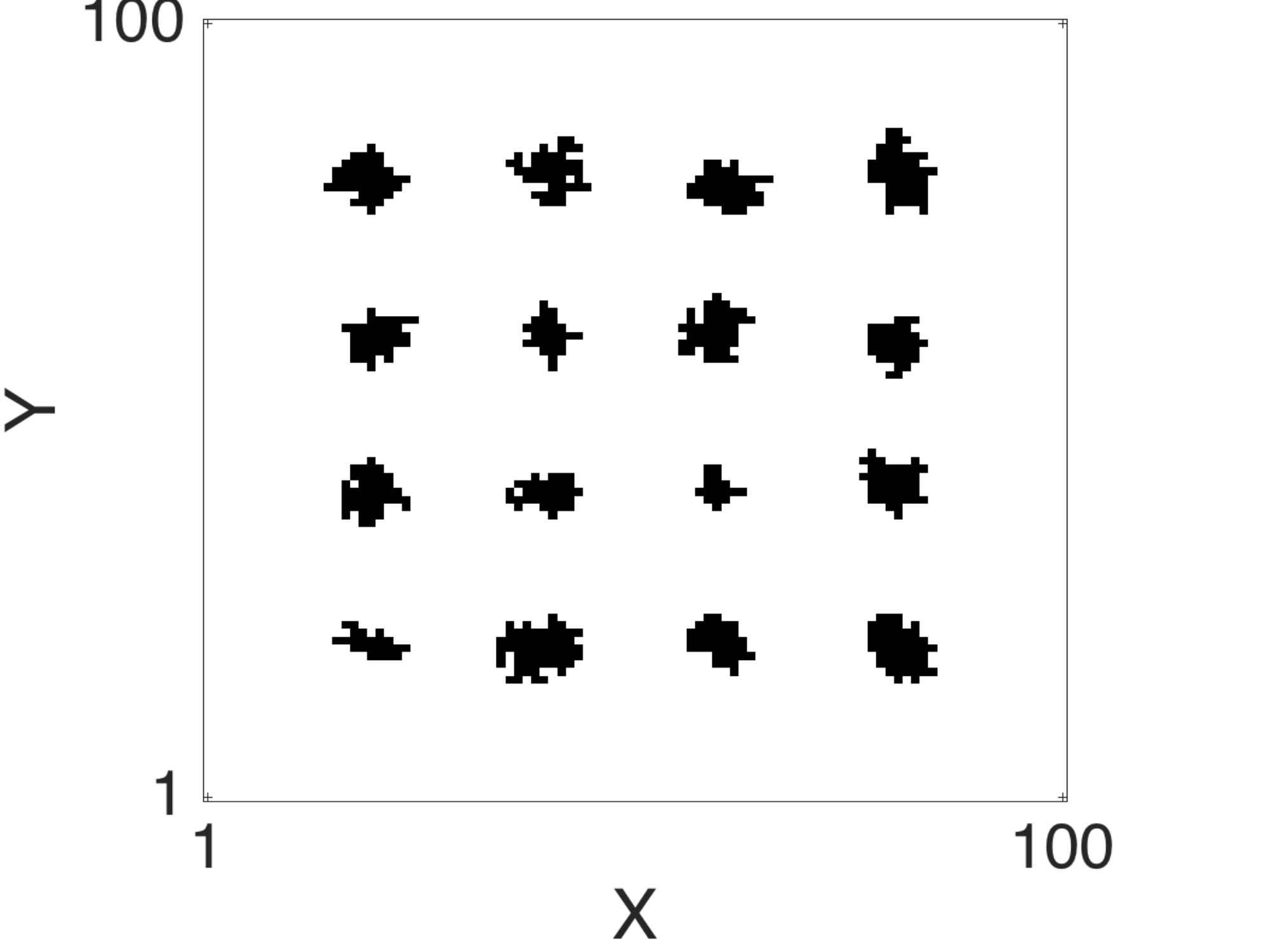}
			\label{fig:ex_binder_a} }
		\subfigure[][]{\includegraphics[width=0.32 \columnwidth]{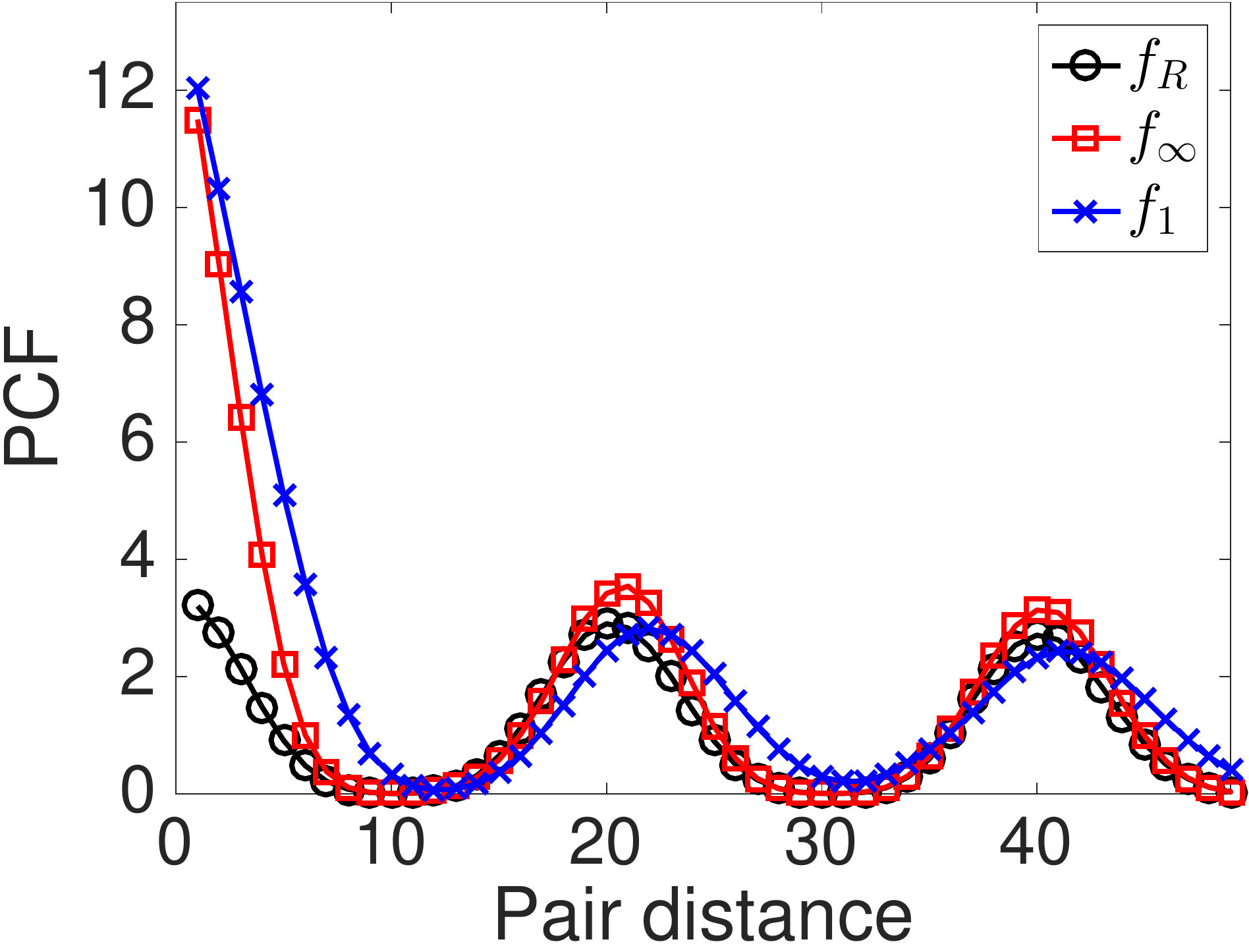}
			\label{fig:ex_binder_b} }\\
	\end{center}
	\caption{Examples of spatial structure analysis. In panels \subref{fig:ex_unif_a} and \subref{fig:ex_unif_b} a case with no spatial correlation is considered. Panel \subref{fig:ex_unif_a} is an example visualisation of an occupancy matrix uniformly populated with density $0.5$. Occupied sites are coloured in black and are white otherwise. Panel \subref{fig:ex_unif_b} shows the four PCFs averaged over $50$ uniformly populated matrices. Panels \subref{fig:ex_binder_a} and \subref{fig:ex_binder_b} refer to a discrete simulation of the agent-based model described in the text at time $t=10$. Panel \subref{fig:ex_binder_a} is a visualisation of the occupancy matrix and in panel \subref{fig:ex_binder_b} our PCFs are compared with the Rectilinear PCF}
	\label{fig:ex_unif} 
\end{figure}

We start by computing the PCFs for a system without any spatial correlation. We consider 50 independent occupancy matrices, $U_i$, $i=1,\dots, 50$, populated uniformly at random with density $0.5$ (see Fig. \ref{fig:ex_unif_a}). For each realisation, $U_i$, we compute the corresponding PCF, $f_d^{U_i}$, and then we average the results over the 50 realisations which we denote $\hat{f}_d$. If the normalisation is correct, $\hat{f}_d(m)$ should return the value unity for every pair distance, $m$, meaning that no spatial correlation is found. In Fig. \ref{fig:ex_unif_b} all four aforementioned averaged PCFs are plotted: $\hat{f}_A, \hat{f}_R, \hat{f}_1$ and $\hat{f}_{\infty}$. The results show that both the averaged Square Uniform PCF, $\hat{f}_1$, and  Square Taxicab PCF, $\hat{f}_\infty$, correctly predict that there is no spatial correlation. The averaged Rectilinear PCF, $\hat{f}_R$, also correctly predicts no spatial correlation. However, the averaged Annular PCF, $\hat{f}_A$, has clear peaks, suggesting, incorrectly, the presence of spatial correlation. Since the results are averaged over multiple repeats, such a discrepancy can not be attributed to stochasticity, but due to incorrect normalisation as explained in Section \ref{sec:previous}. Note that the Annular PCF can still correctly identify spatial correlation in many examples, however, the incorrect normalisation often makes the results hard to interpret. This is because it makes it difficult to distinguish between genuine correlation and systematic error. For this reason, for the rest of this Section, we omit the results of the Annular PCF and continue to compare between our PCFs and the Rectilinear PCF using non-periodic BC.

Next we consider examples of strong spatial correlation. Fig. \ref{fig:ex_binder_a} shows an example of aggregation driven by a proliferation mechanism. The occupancy matrix is obtained by simulating an on-lattice agent-based model with periodic BC as described in \citet{binder2013qss}, which we summarise as follows. The model is initialised with 16 agents, located at coordinates given by $\lbrace (x,y) \,|\, x,y \in \lbrace 20,40,60,80 \rbrace \rbrace$ on a regular square lattice with $L_x= 100$, $L_y = 100$. Time is discretised with a time step $\tau=1$ and the number of agents at time $t$ is denoted by $n(t)$. At each time step the configuration at time $t+\tau$ is obtained from the configuration at time $t$, by repeating the following steps $n(t)$ times. 1) An agent is chosen uniformly at random from the $n(t)$ agents present at the end of the previous time-step; 2) one of its four von Neumann neighbours is selected at random with equal probability; 3) if the selected site is empty, a new agent is placed in this site and $n(t+\tau)=n(t)+1$, otherwise the configuration is left unchanged. 

Fig. \ref{fig:ex_unif} \subref{fig:ex_binder_a} shows a single realisation after 10 time steps and
Fig. \ref{fig:ex_unif} \subref{fig:ex_binder_b} shows the corresponding PCFs: $f_R$, $f_1$ and $f_\infty$. The results indicate that all of the  PCFs considered correctly identify aggregation. However, the quantitative information about aggregate sizes at different length scales provide by each PCF varies. For example, we see that all PCFs in Fig. \ref{fig:ex_unif} \subref{fig:ex_binder_b} exhibit three peaks; at $m=1$, $m \approx 20$ and $m \approx 40$. The different peaks and troughs of the PCF profiles have different qualitative meanings related to the correlation type.  
Due to the local approach of the Square Taxicab PCF and Square Uniform PCF, the first peak at $m=1$ is three times higher than the peaks at larger values of distance. These differences in amplitude highlight the different peak origins. Specifically, the first and highest peak distinguishes the individual cluster aggregate and the later peaks indicate correlation between different clusters. In contrast, all three peaks in the Rectilinear PCF are the same amplitude. 
Note that, in the case of aggregation, the average diameter of the aggregate corresponds to the first value of distance which achieves the minimum of the PCF. The Rectilinear, Square Uniform and Square Taxicab PCFs estimate the aggregate diameter to be 9, 9 and 11 respectively.
Importantly, the PCFs capture the fact that this diameter depends on the metric used. In particular the distance between two sites measured using the uniform metric is always less than or equal to the taxicab distance. This phenomenon is seen more clearly in later examples. 

\begin{figure}[h!!]
	\begin{center}
		\subfigure[][]{\includegraphics[ width=0.32 \columnwidth]{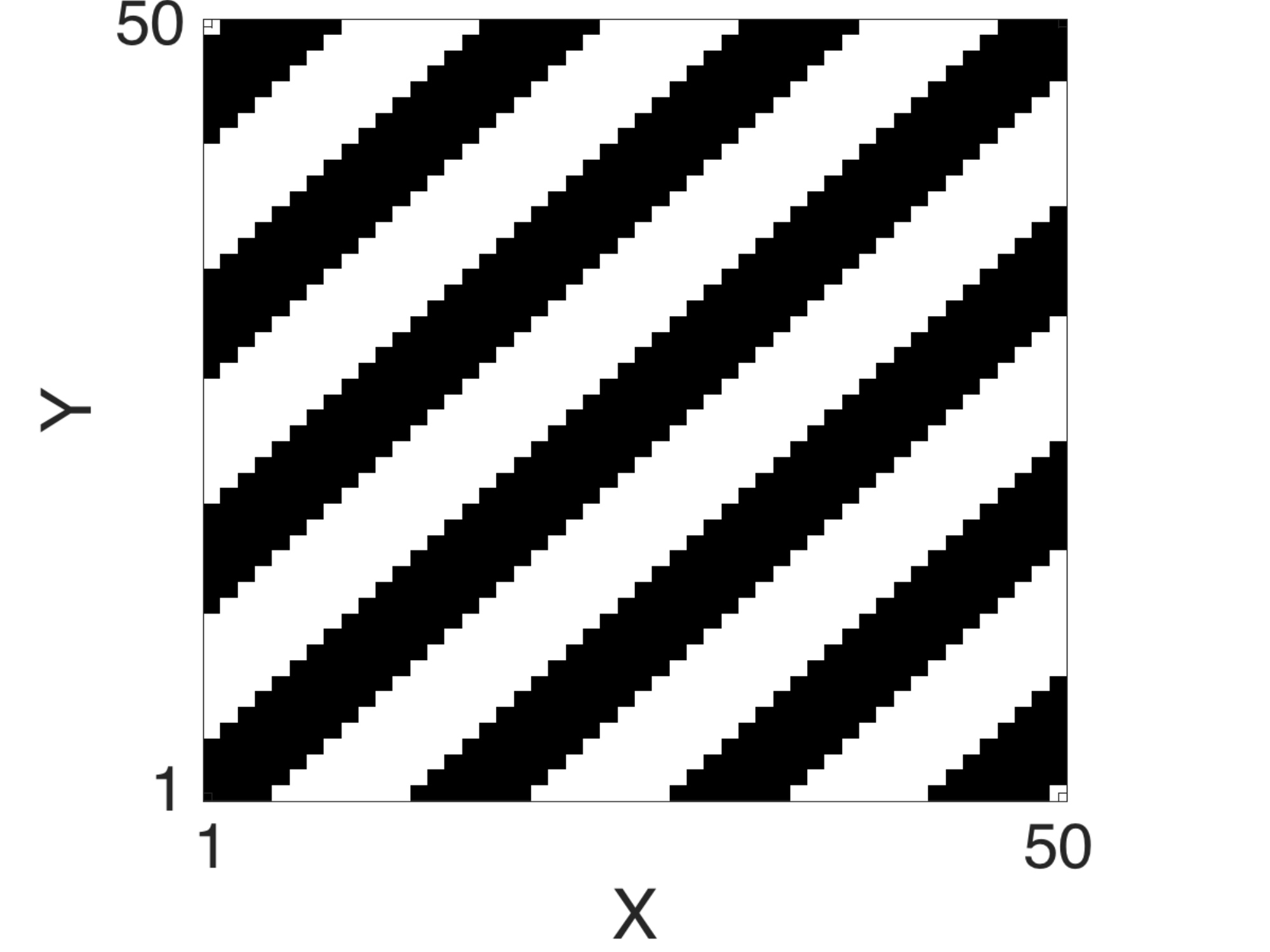}\label{fig:ex_stripes_matrix}\label{fig:ex_stripes_a}}
		\subfigure[][]{\includegraphics[width=0.32 \columnwidth]{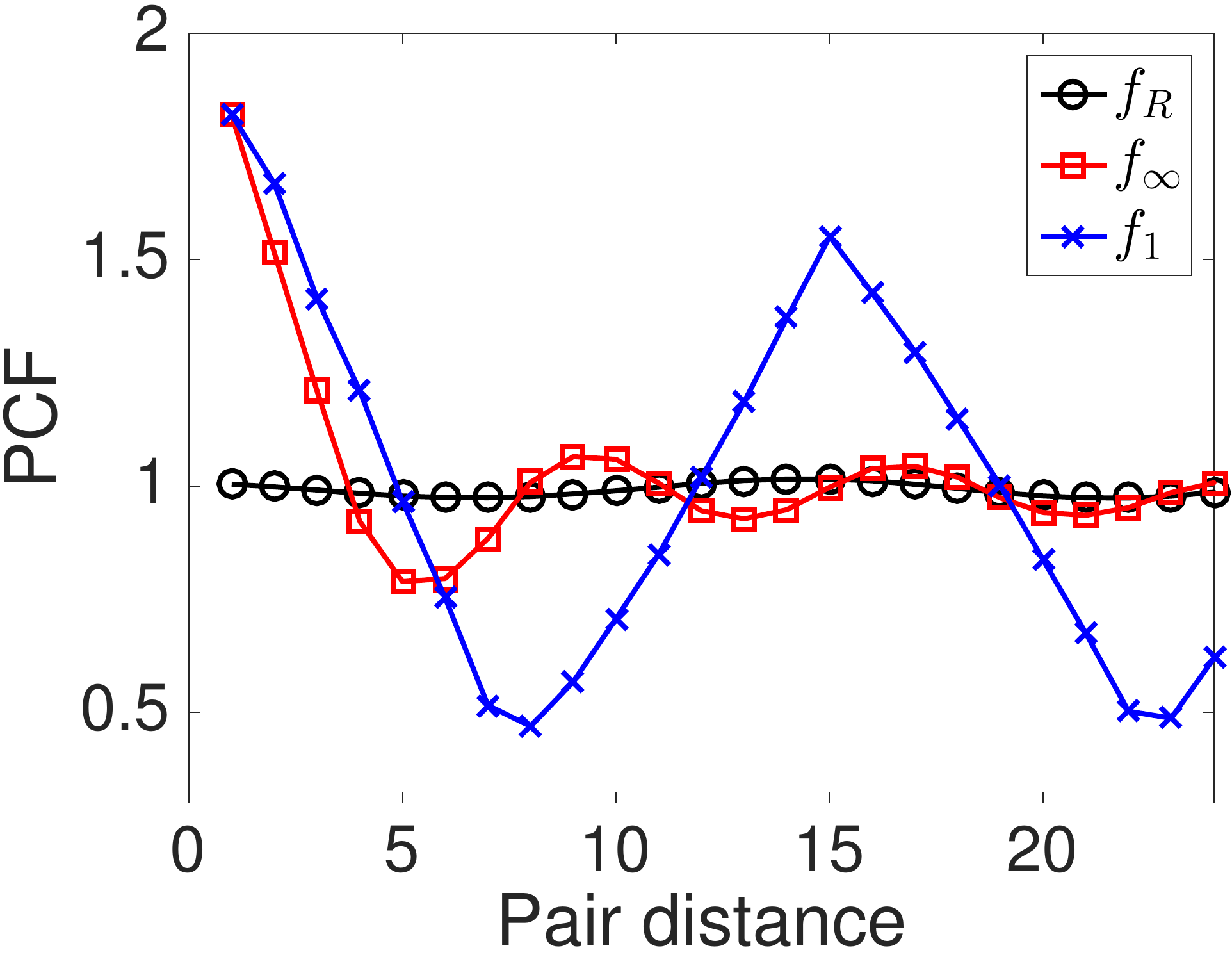}\label{fig:ex_stripes_b}}\\
		\subfigure[][]{\includegraphics[ width=0.32 \columnwidth]{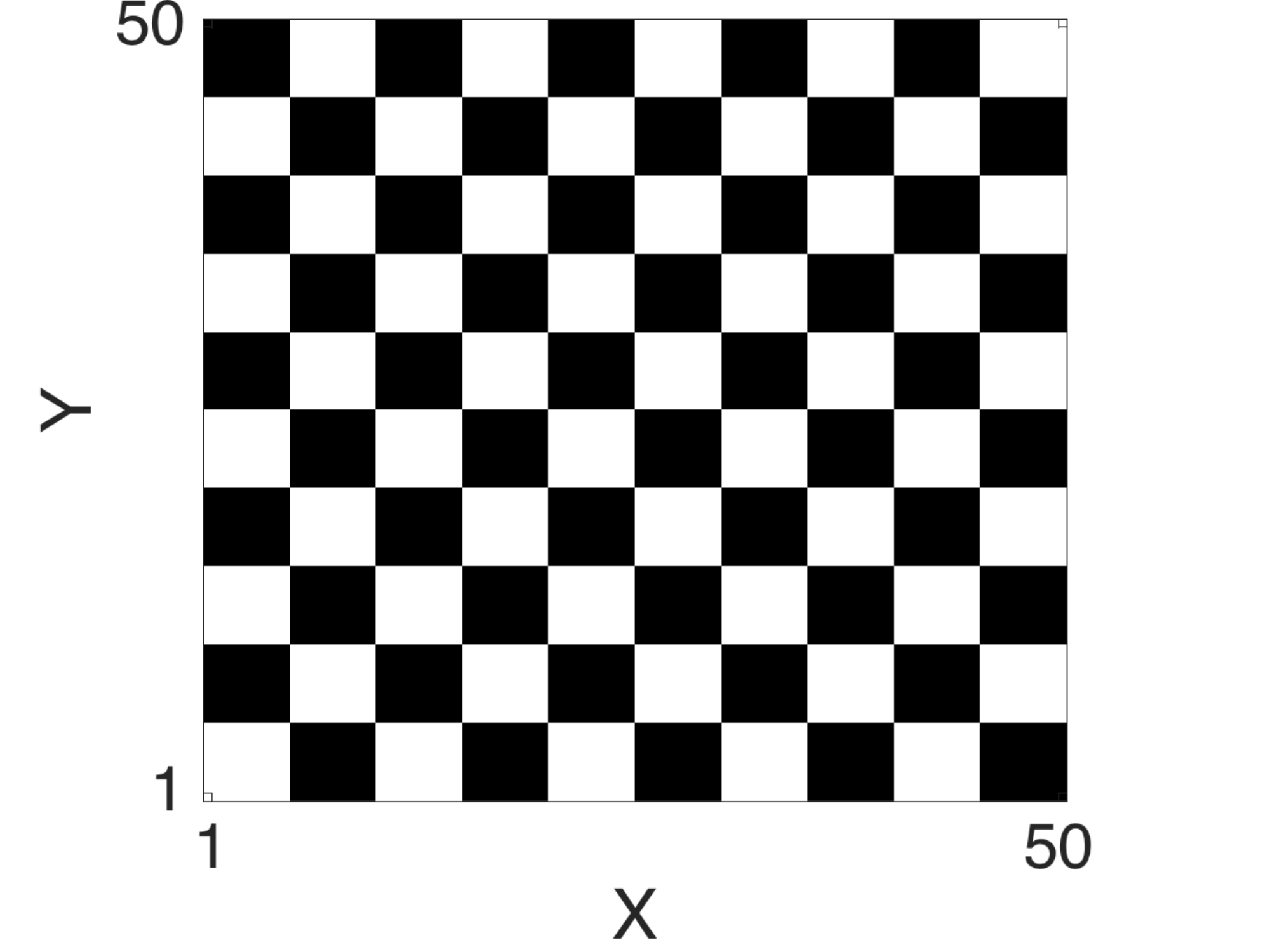}\label{fig:ex_chess_a}}
		\subfigure[][]{\includegraphics[width=0.32 \columnwidth]{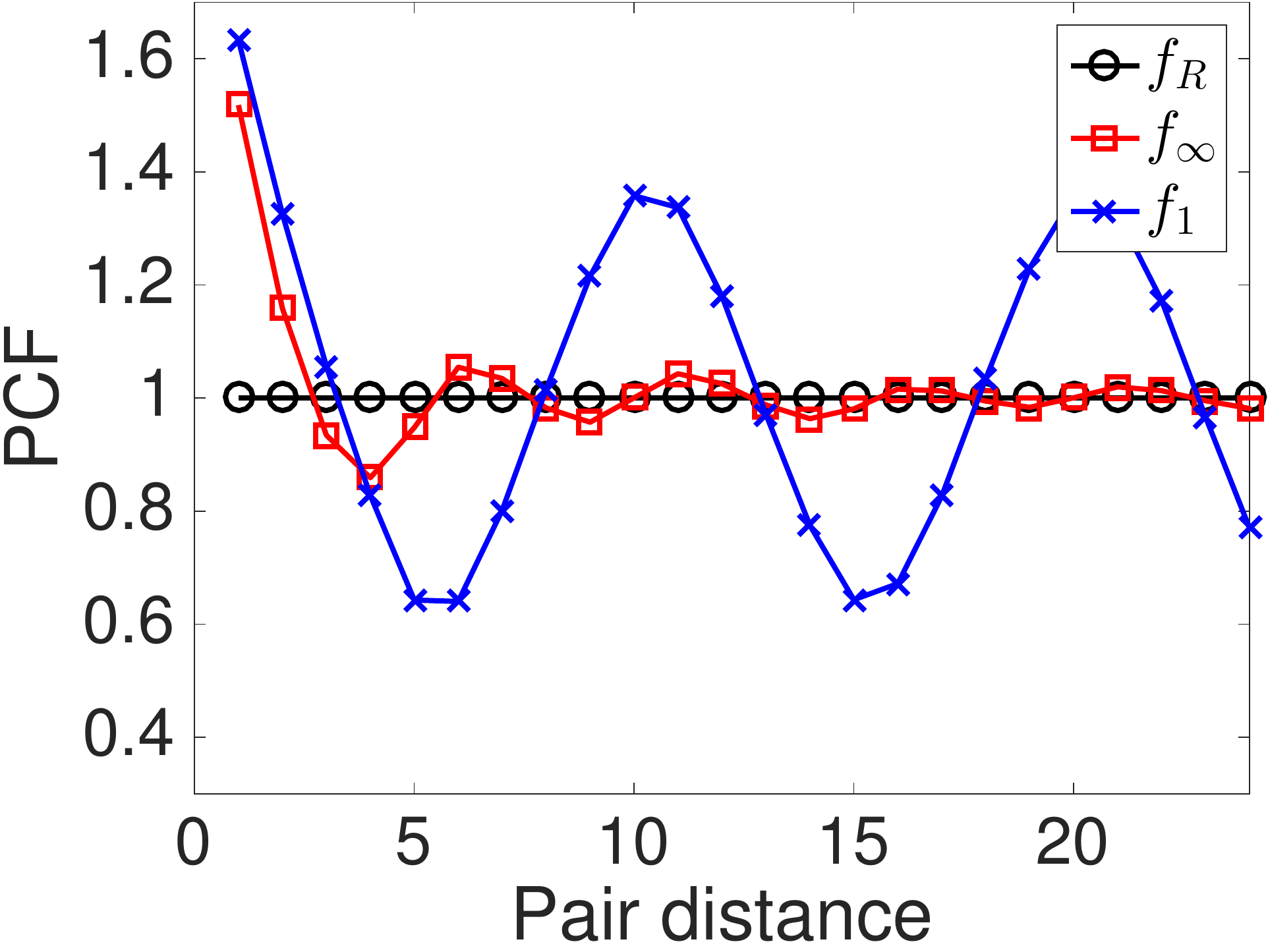}\label{fig:ex_chess_b}}\\
		\subfigure[][]{\includegraphics[ width=0.32 \columnwidth]{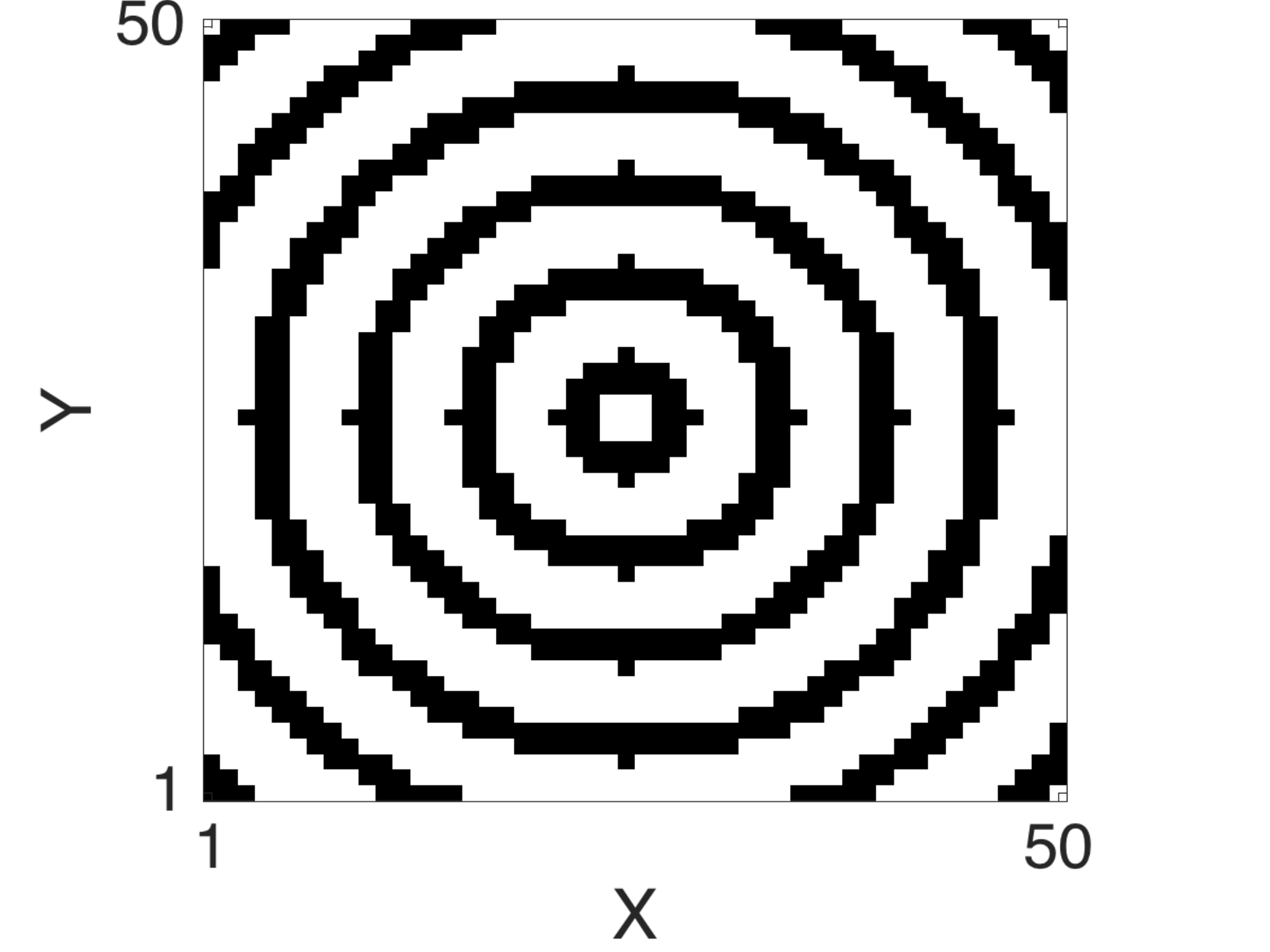}\label{fig:ex_circles_a}}
		\subfigure[][]{\includegraphics[width=0.3 \columnwidth]{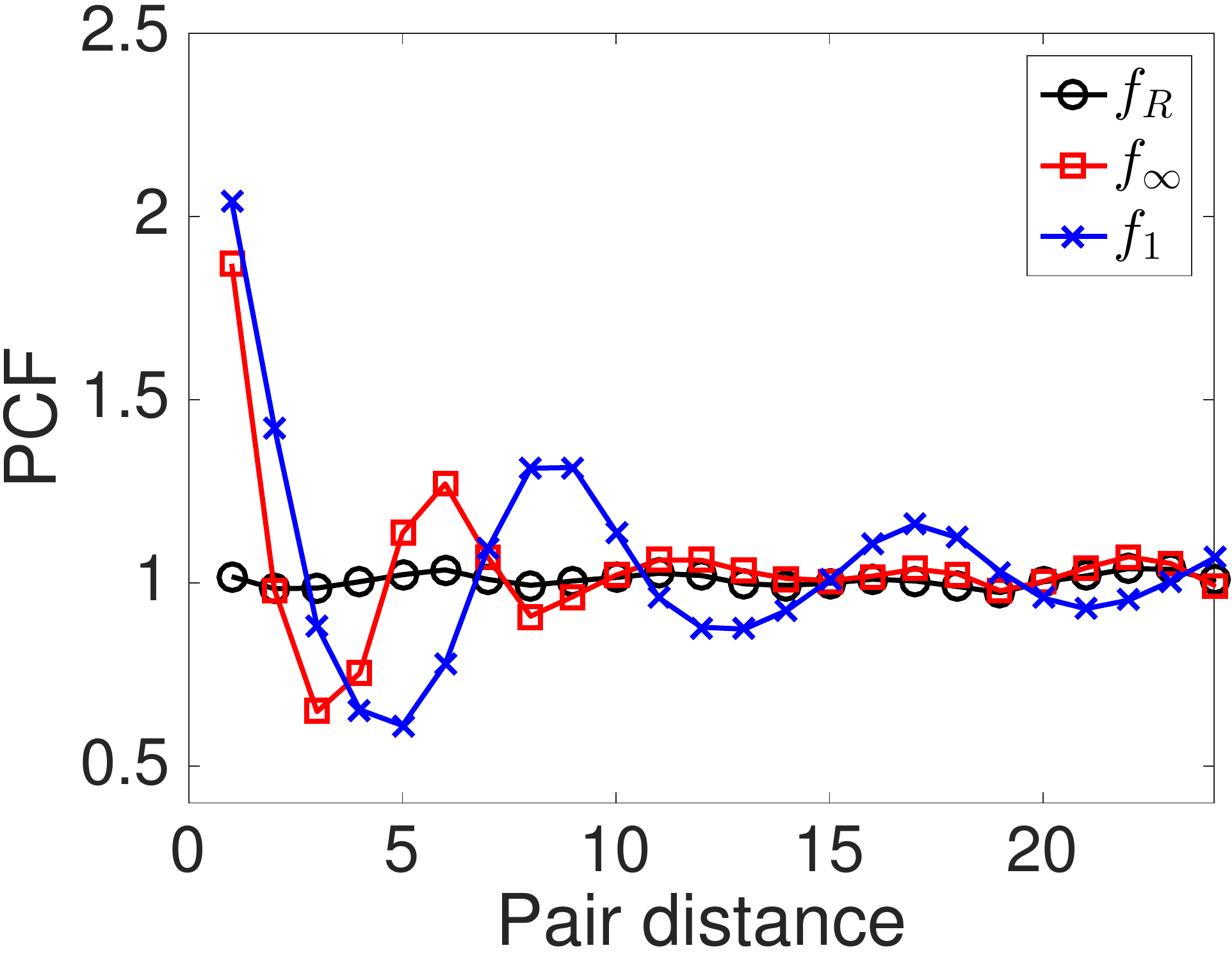}\label{fig:ex_circles_b}}
	\end{center} 
	\caption{Examples of pattern analysis. Panels \subref{fig:ex_stripes_a}, \subref{fig:ex_chess_a} and \subref{fig:ex_circles_a} visualise three constructed spatial patterns. Panels \subref{fig:ex_stripes_b}, \subref{fig:ex_chess_b} and \subref{fig:ex_circles_b} displays the corresponding Square Taxicab, Square Uniform and Rectilinear PCF using the occupancy matrices for \subref{fig:ex_stripes_a}, \subref{fig:ex_chess_a} and \subref{fig:ex_circles_a}, respectively.}
	\label{fig:ex_scs} 
\end{figure}

We now consider a series of examples with spatial correlation constructed artificially in order to compare and evaluate the different PCFs. In Fig. \ref{fig:ex_scs} we compare our Square Uniform and Square Taxicab PCF with the Rectilinear PCF for three different patterns with strong spatial correlation. All three examples (Fig. \ref{fig:ex_scs} \subref{fig:ex_stripes_a}: diagonal stripes, Fig. \ref{fig:ex_scs} \subref{fig:ex_chess_a}: chessboard pattern and Fig. \ref{fig:ex_scs} \subref{fig:ex_circles_a}: concentric circles) are chosen so that the column- and row-averaged densities are constant and hence the spatial structure is not recognised by the Rectilinear PCF, as shown in Figs. \ref{fig:ex_scs} \subref{fig:ex_stripes_b}, \subref{fig:ex_chess_b} and \subref{fig:ex_circles_b}. This is in contrast to the approach of our PCFs (both Square Uniform and Square Taxicab) which successfully recognise the spatial structure in all three examples. In addition, these examples uncover other interesting differences between the taxicab and the uniform approaches. Consider the PCF for the case of diagonal stripes and the chessboard pattern (Fig. \ref{fig:ex_scs} \subref{fig:ex_stripes_b} and Fig. \ref{fig:ex_scs} \subref{fig:ex_chess_b}). Here the Square Uniform PCF quickly converges to unity (no spatial correlation) for large distance $m$, whilst in both cases, the Square Taxicab PCF still shows a strong oscillatory behaviour for large distance $m$ suggesting spatial correlation. To give an intuitive explanation of this phenomenon, let us consider the shapes of balls of size $m$ centred at a given site $\boldsymbol{a}$ under the two metrics. These balls are defined as $B_m(\boldsymbol{a})=\lbrace \boldsymbol{b} \in \mathbb{L} \, |\, \norm{\boldsymbol{a}-\boldsymbol{b}}_d \le m\rbrace$ with $d=1,\infty$, respectively (see Fig. \ref{fig:normsvis}). The ball corresponding to the uniform metric (Fig. \ref{fig:normsvis} \subref{fig:normsvis_UN}) has a square shape with the sides aligned with the directions of the axis. This implies that when distance, $m$, becomes close to either $L_x$ or $L_y$ in size, the ball corresponding to the uniform metric of distance $m$ contains most of the sites in the corresponding row or column at distance $m$. For large $m$, therefore the uniform metric begins to work in a similar way to the Rectilinear PCF and thus fails to recognise anisotropic patterns biased in the axial directions. 
The ball of the taxicab metric (see Fig \ref{fig:normsvis} \subref{fig:normsvis_TA}), however, has a diamond shape. Consequently, the long-distance correlations appear clear even for patterns in which both the average column and row densities are constant, as in Fig. \ref{fig:ex_scs}.

\begin{figure}[h!!]
	\begin{center}
		\subfigure[][]{\includegraphics[width=0.32 \columnwidth]{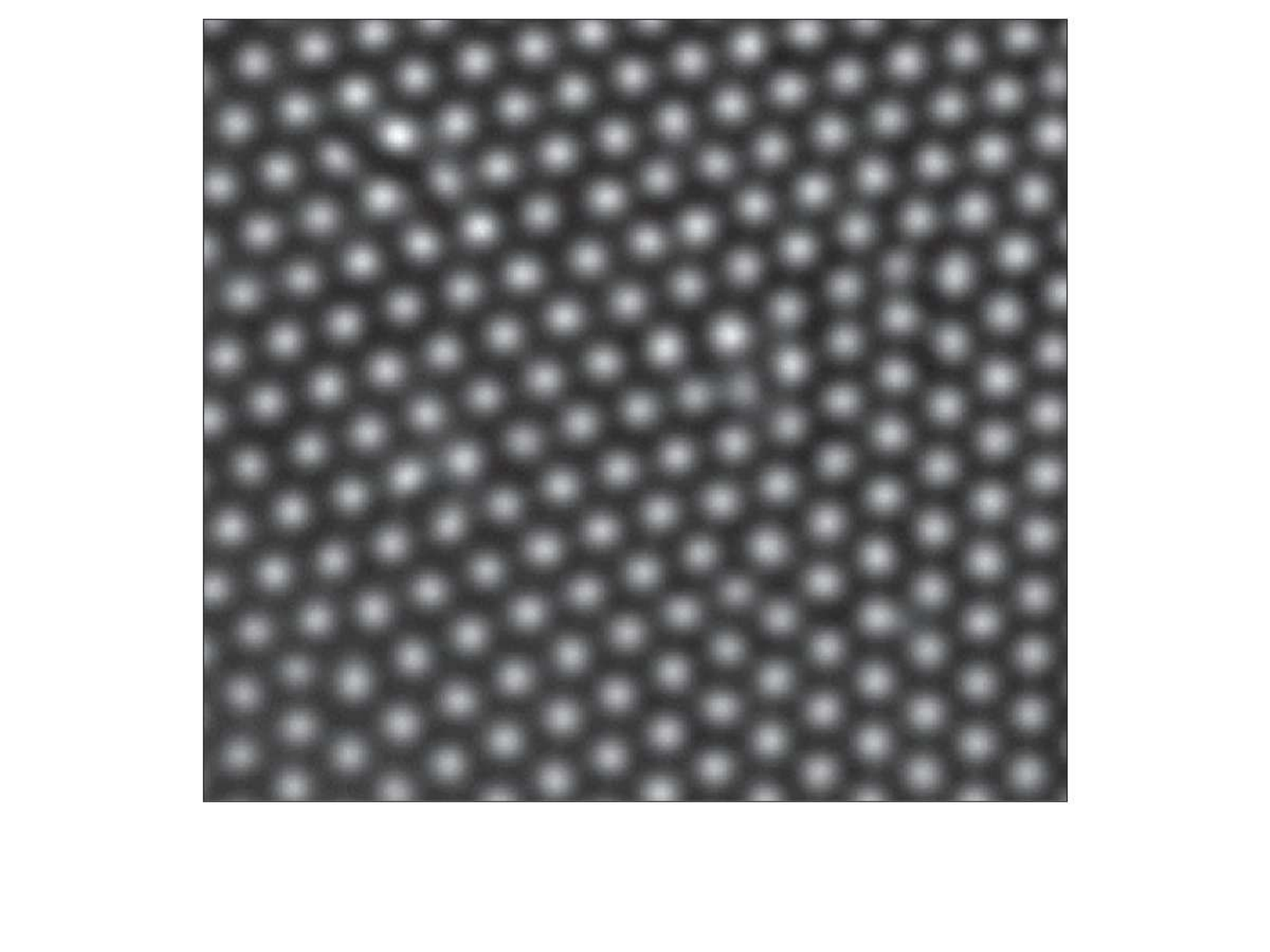}
			\label{fig:ex_a}} 
		\subfigure[][]{\includegraphics[ width=0.32\columnwidth]{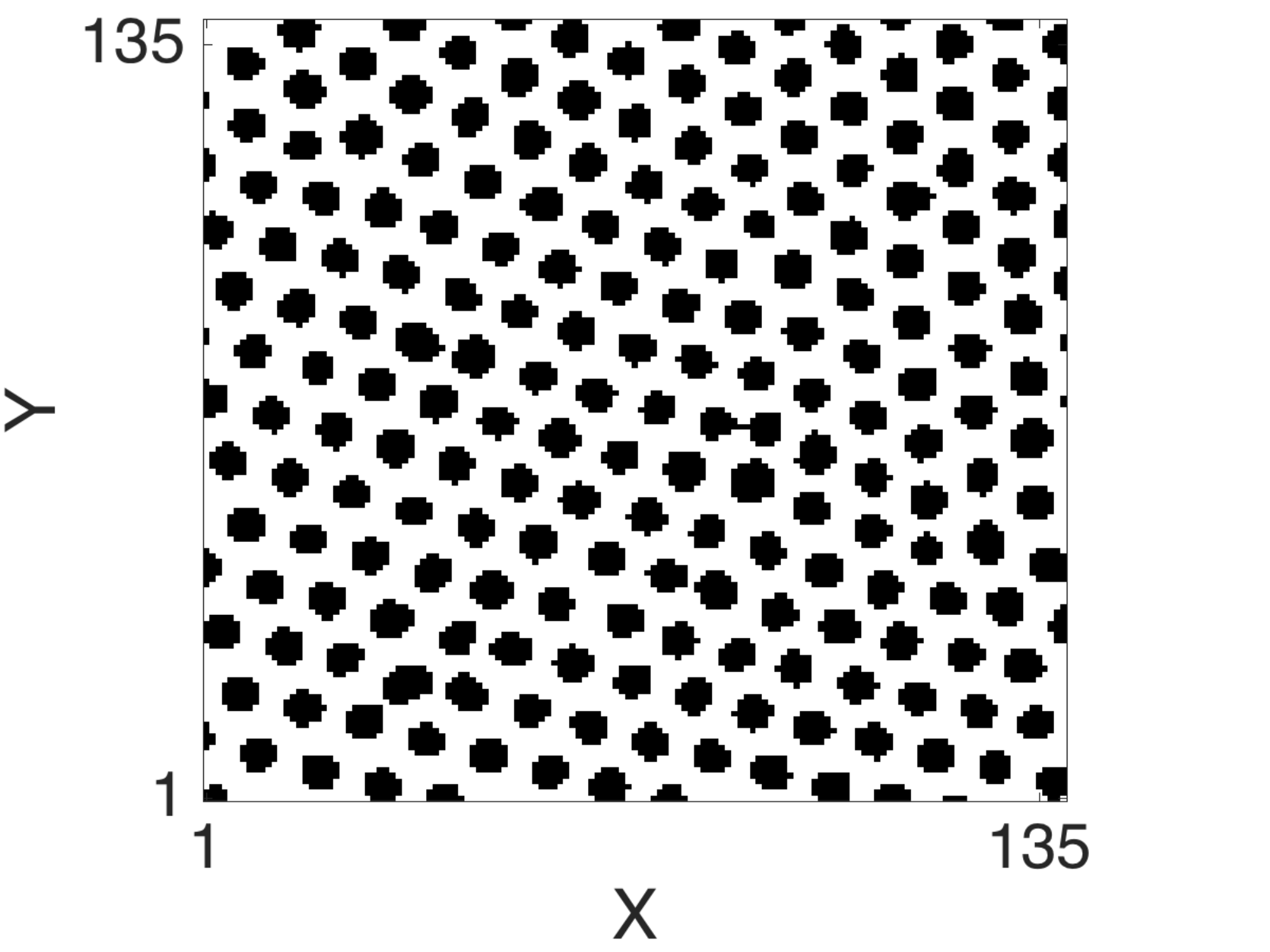}\label{fig:ex_a_matrix}}
		\subfigure[][]{\includegraphics[width=0.32 \columnwidth]{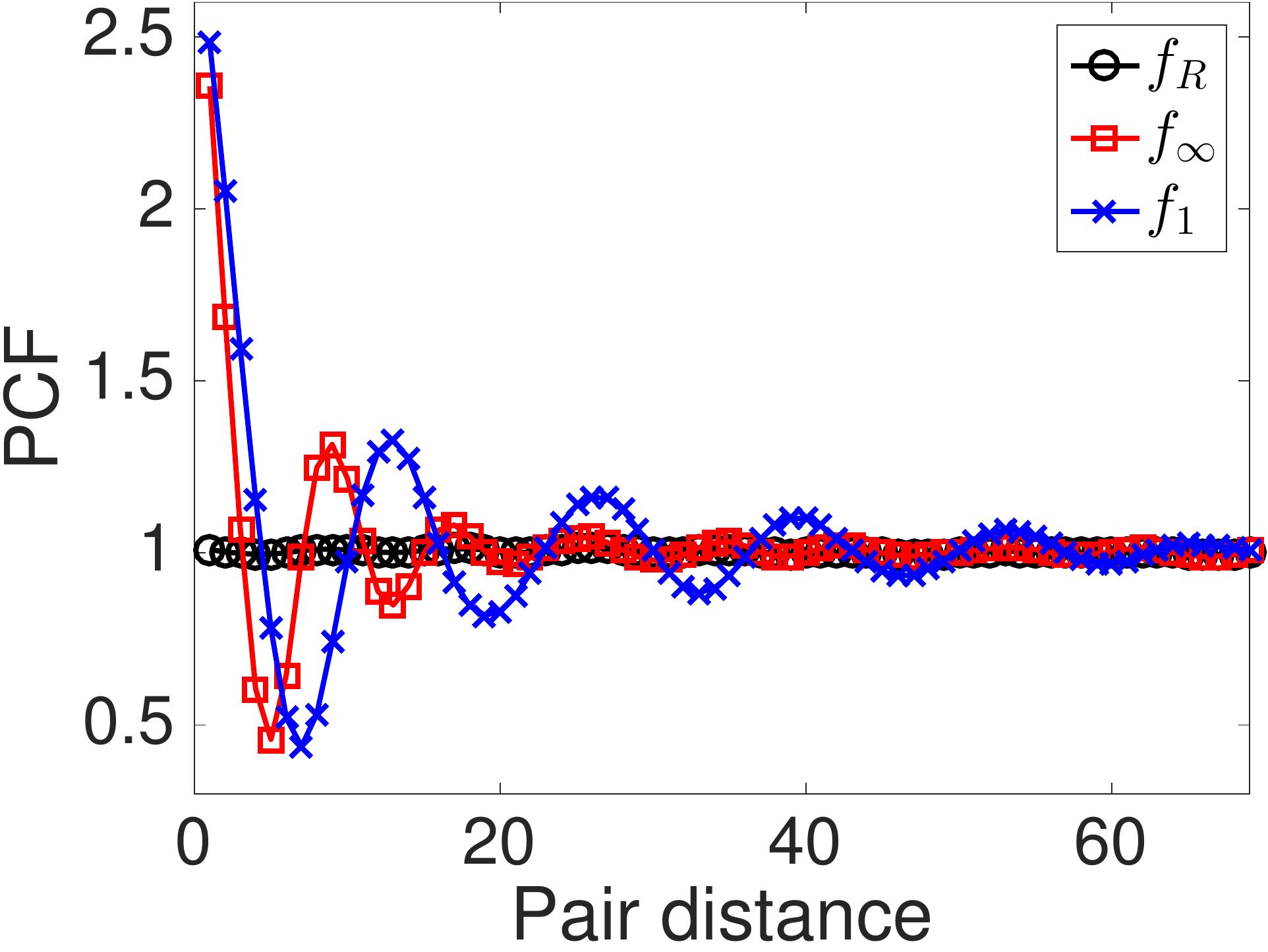}\label{fig:ex_a_PCF}} \\
		\subfigure[][]{\includegraphics[width=0.33 \columnwidth]{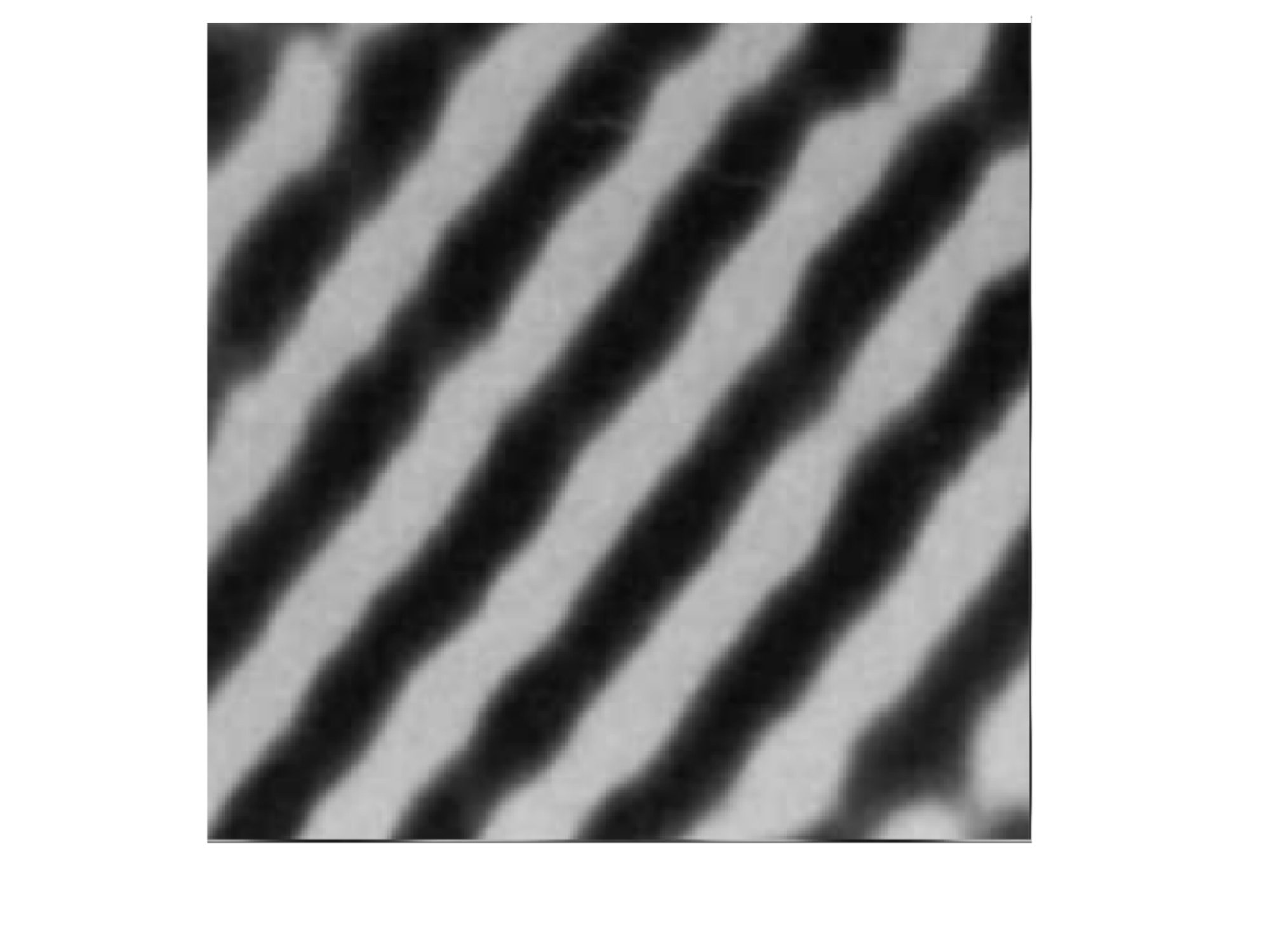}
			\label{fig:ex_b} } 
		\subfigure[][]{\includegraphics[ width=0.32\columnwidth]{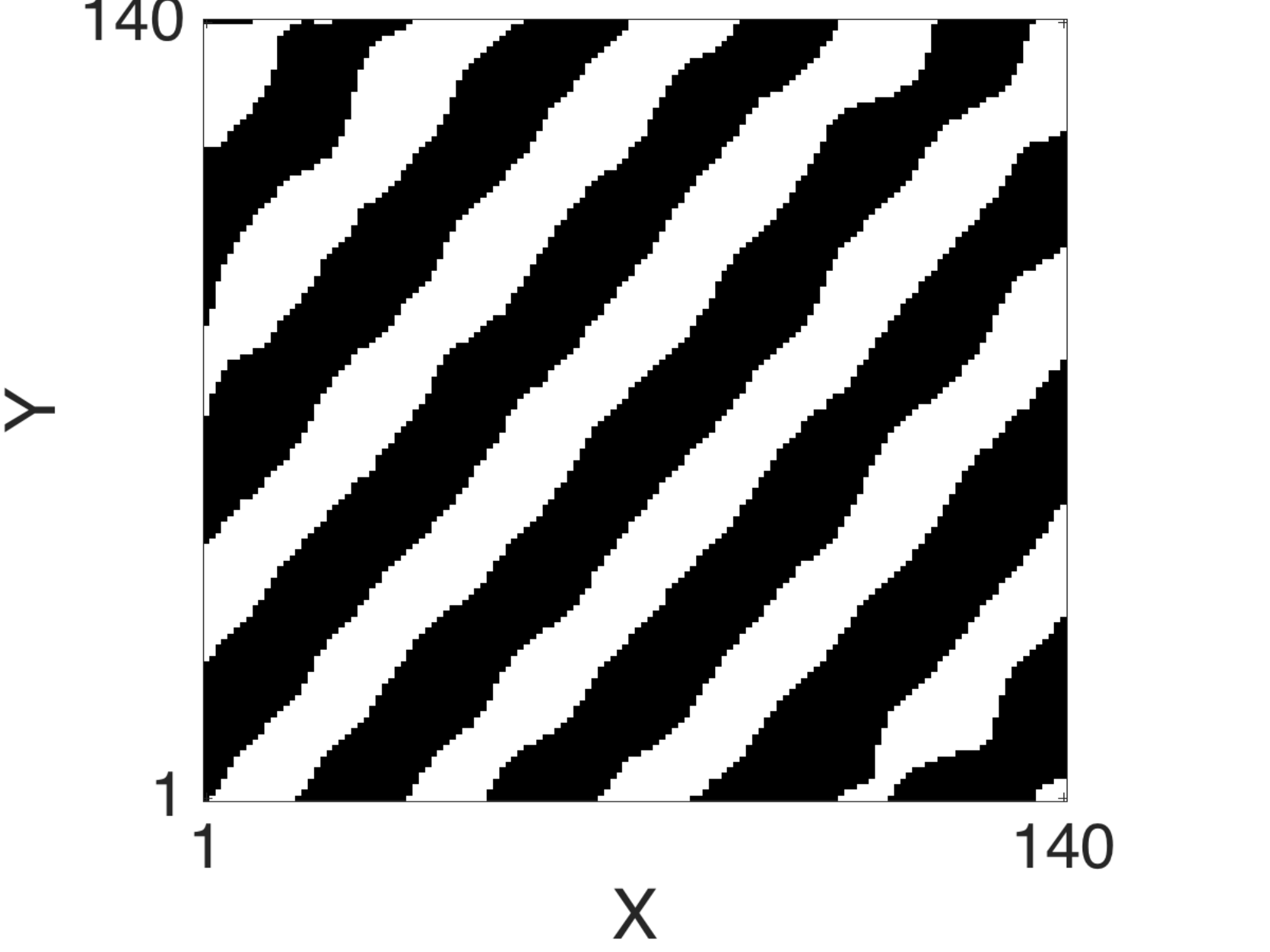}\label{fig:ex_b_matrix}}
		\subfigure[][]{\includegraphics[width=0.32 \columnwidth]{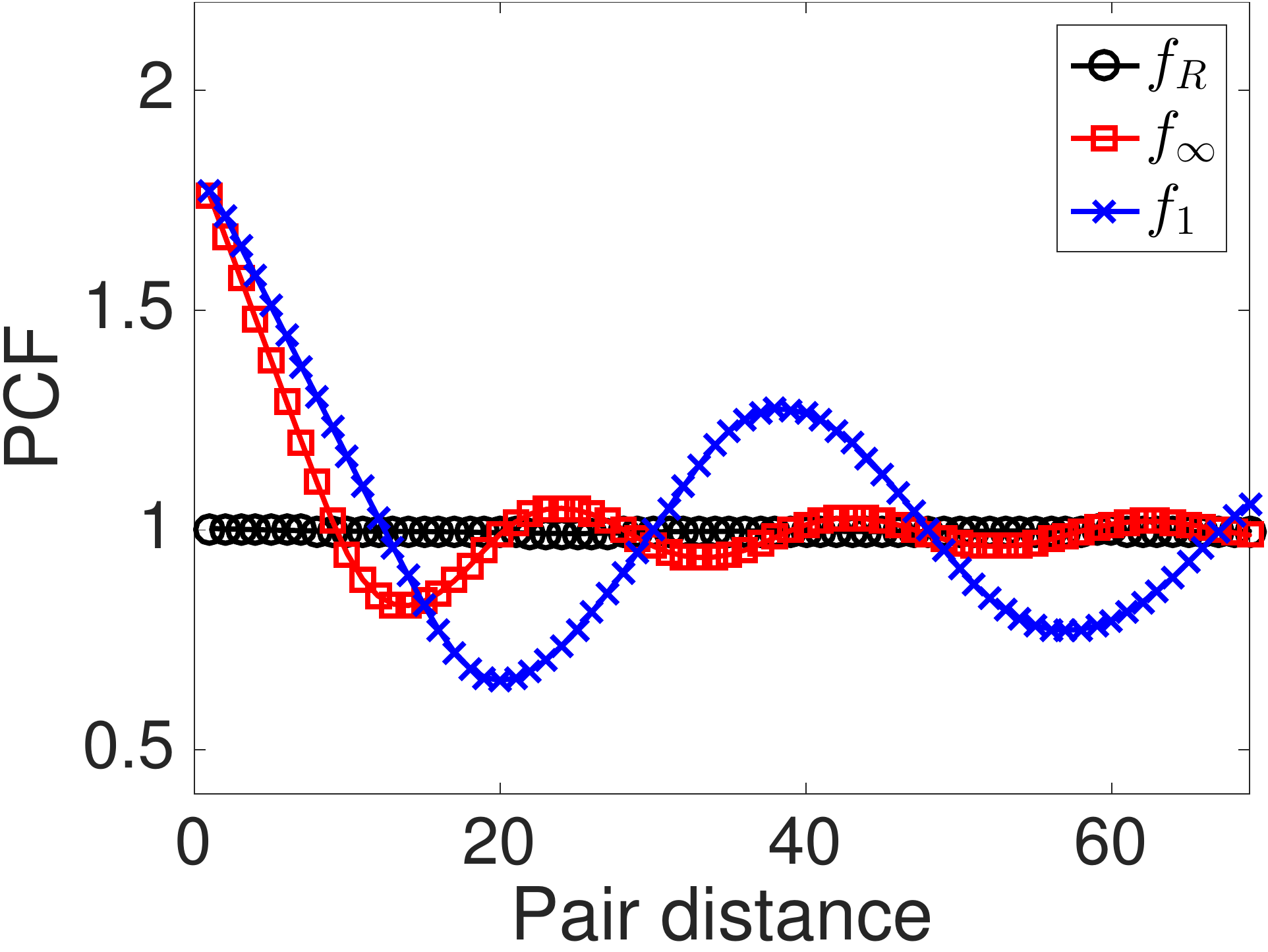}\label{fig:ex_b_PCF}}  \\
		\subfigure[][]{\includegraphics[width=0.32 \columnwidth]{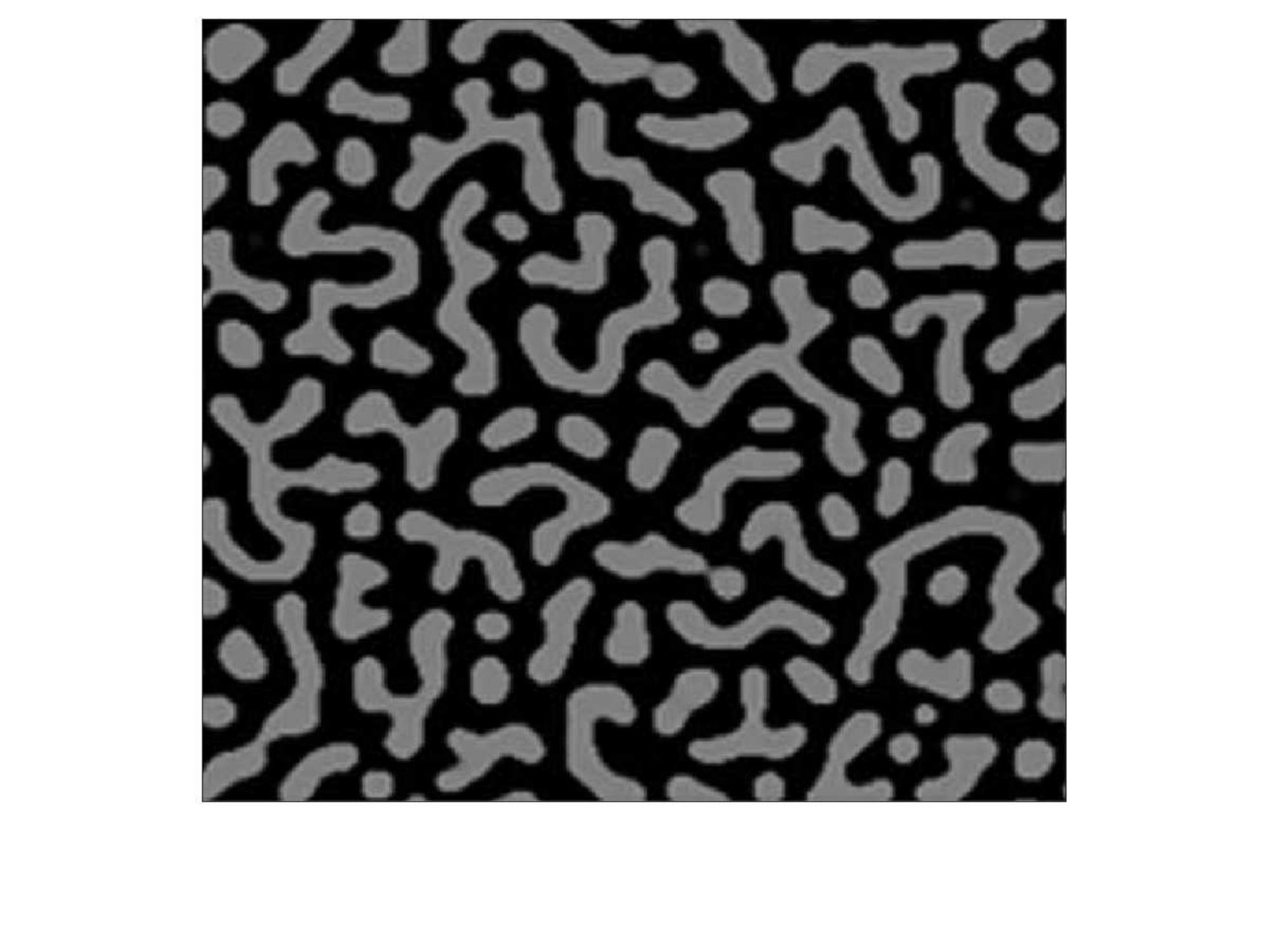}
			\label{fig:ex_c} } 
		\subfigure[][]{\includegraphics[ width=0.32\columnwidth]{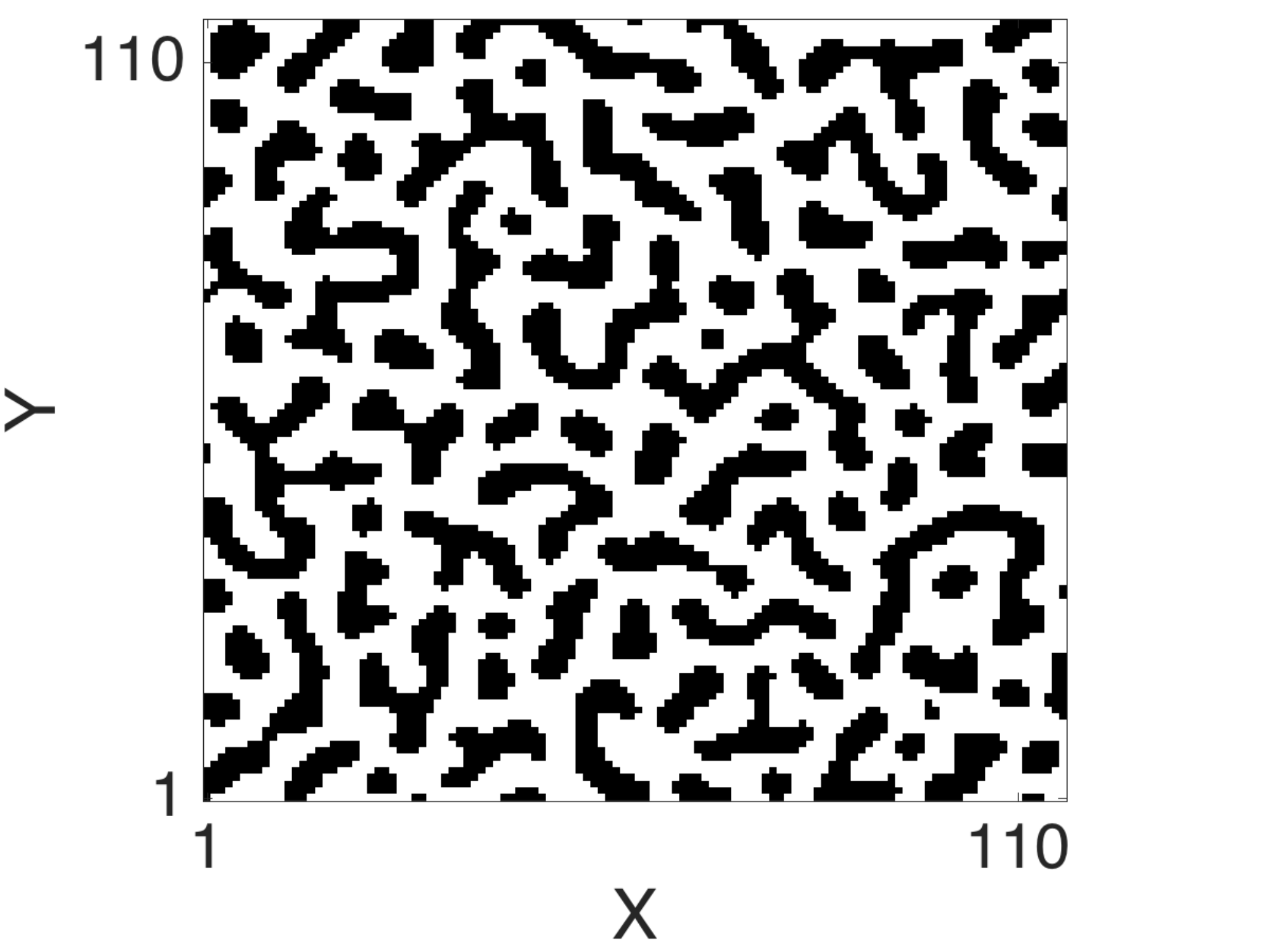}\label{fig:ex_c_matrix}}
		\subfigure[][]{\includegraphics[width=0.32 \columnwidth]{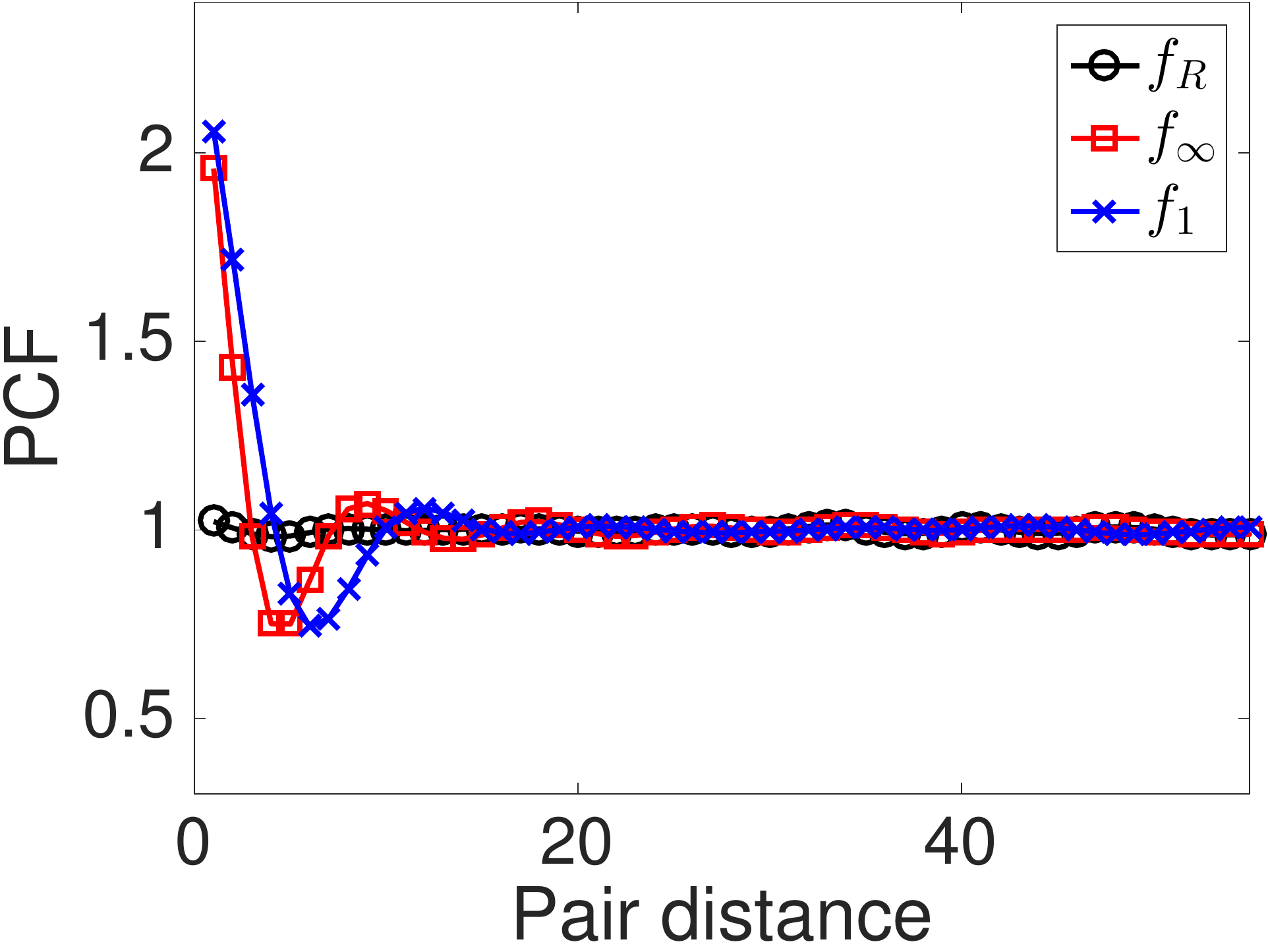}\label{fig:ex_c_PCF}} 
	\end{center}
	\caption{Spatial analysis of Turing patterns. Panels \subref{fig:ex_a}, \subref{fig:ex_b} and \subref{fig:ex_c} show original images representing the results of a reaction-diffusion mechanism between two chemical substances, \subref{fig:ex_a} is reprinted from \cite{ball2015fpm}, \subref{fig:ex_b} from \cite{ouyang1991tuh} and \subref{fig:ex_c} from \cite{kondo2017ukt}. Panels \subref{fig:ex_a_matrix}, \subref{fig:ex_b_matrix} and \subref{fig:ex_c_matrix} visualise the occupancy matrices corresponding to the original images (described in text). In the panels \subref{fig:ex_a_PCF}, \subref{fig:ex_b_PCF} and \subref{fig:ex_c_PCF} we compare the Square Taxicab, Square Uniform and Rectilinear PCFs for each of the examples.}
	\label{fig:ex_1} 
\end{figure}

The examples in Fig. \ref{fig:ex_scs} were constructed specifically to underline the main differences between the three PCFs. Nevertheless, similar patterns also arise in many biologically and mathematically relevant applications \cite{bhide1999dsc,ouyang1991tuh,murray2007mbi}. We conclude this section by comparing the three PCF approaches applied to some real-world examples taken from the literature. In Fig. \ref{fig:ex_1} we analyse three images representing examples of Turing patterns. A corresponding occupancy matrix is obtained by representing each pixel of the image as a value in a matrix which is 1 (i.e. occupied) if the three values of the RGB colourisation of the pixel are above a certain threshold ($80$) and 0 otherwise. In all cases the column and row densities are almost constant, hence the spatial structure again remains largely undetected by the Rectilinear PCF, whilst our Square Uniform and Square Taxicab PCF correctly identify the patterns. As already observed in the previous examples, we note that the estimated diameter of the aggregate, wavelength and amplitude of the oscillations differ according to the metric used.

\section{The Triangle, Hexagon and Cube PCFs} 
\label{sec:tri_hex}
Despite the square lattice being the most popular set up for spatially discrete models \cite{simpson2009mss,yates2012gfm,ross2015icc, baker2010fmm, deutsch2007cam}, in some situations other types of tessellation, either regular or irregular, can be more suitable \cite{browning2017bca, keeling1999els, deutsch2007cam}. 

In the following subsections we extend our definition of the PCFs in Section \ref{sec:PCF} to more general types of tessellations. 
We define the Triangle, Hexagon, Cube Uniform and Cube Taxicab PCFs under non-periodic and periodic BC for triangular, hexagonal and cuboidal tessellations respectively. 
the following subsections represent qualitative discussions of the different cases. We refer the reader to the SM Sections \ref{SUPP-sec:norm_tri} and \ref{SUPP-sec:norm_3D} for the full details of the derivation of the PCF formulae.

\subsection{Triangle and Hexagon PCF}
\label{sec:tria-hex}

\begin{figure}[h!!]
	\begin{center}
		\subfigure[][]{\includegraphics[height=3cm]{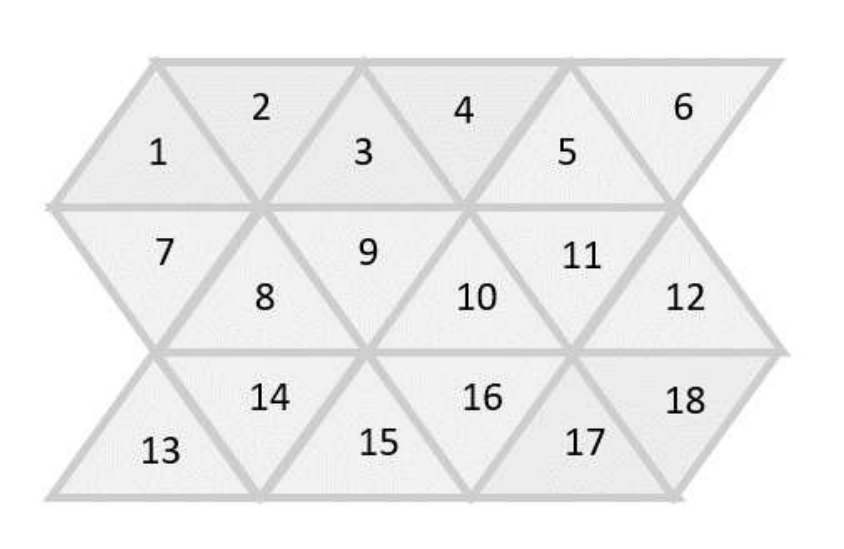}
			\label{fig:ex_a}} 
		\subfigure[][]{\includegraphics[height=3cm]{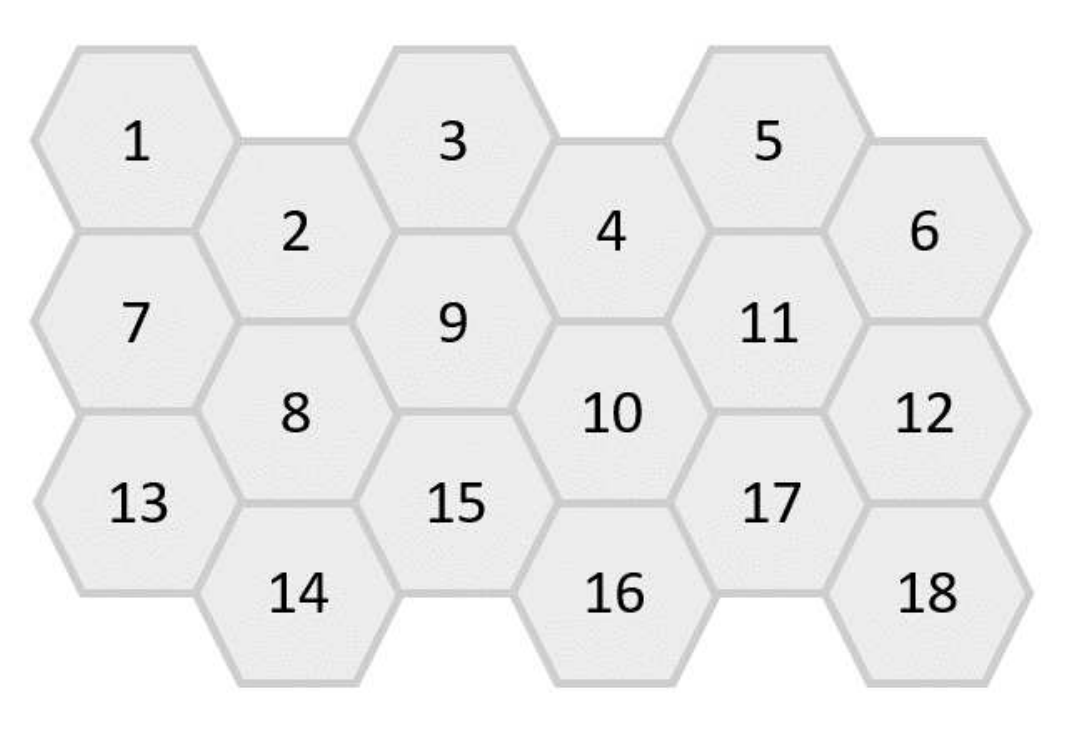}}
	\end{center}
	\caption{Example domains for (a) triangular and (b) hexagonal tessellations in which $L_x=6$ and $L_y=3$.}
	\label{fig:tex} 
\end{figure}

First, we define triangularly and hexagonally tessellated domains of size $L_x \times L_y$. These comprise is an array of $L_y$ rows of $L_x$ regular triangles or hexagons respectively. Examples for which $L_x = 6$ and $L_y = 3$ for each of the two cases are given in Fig. \ref{fig:tex} (a) and (b), respectively. Notice that, for periodic BC to be meaningful in these domain definitions $L_x$ must be even. Therefore, we enforce this as a condition in what follows.

In the context of triangular and hexagonal tessellations we focus our attention on the taxicab metric, both for simplicity and as the most natural metric on this domain type. Using the taxicab metric, the number of sites of distance $m$ from any given reference site is given by
\begin{subequations}
\label{eq:t_trihex}
\begin{align}
t_{tri}(m)&=3m,
\\t_{hex}(m)&=6m.
\end{align}
\end{subequations}
The proofs of the expressions \eqref{eq:t_trihex} are omitted, but they can be obtained easily by induction on $m$. Examples for $m=1,2,3$ are visualised in Fig. \ref{fig:tm_tex}. Using the same reasoning as in Section \ref{sec:PCF} under periodic BC:

\begin{figure}[h!!]
	\begin{center}
		\subfigure[][]{\includegraphics[height=4cm]{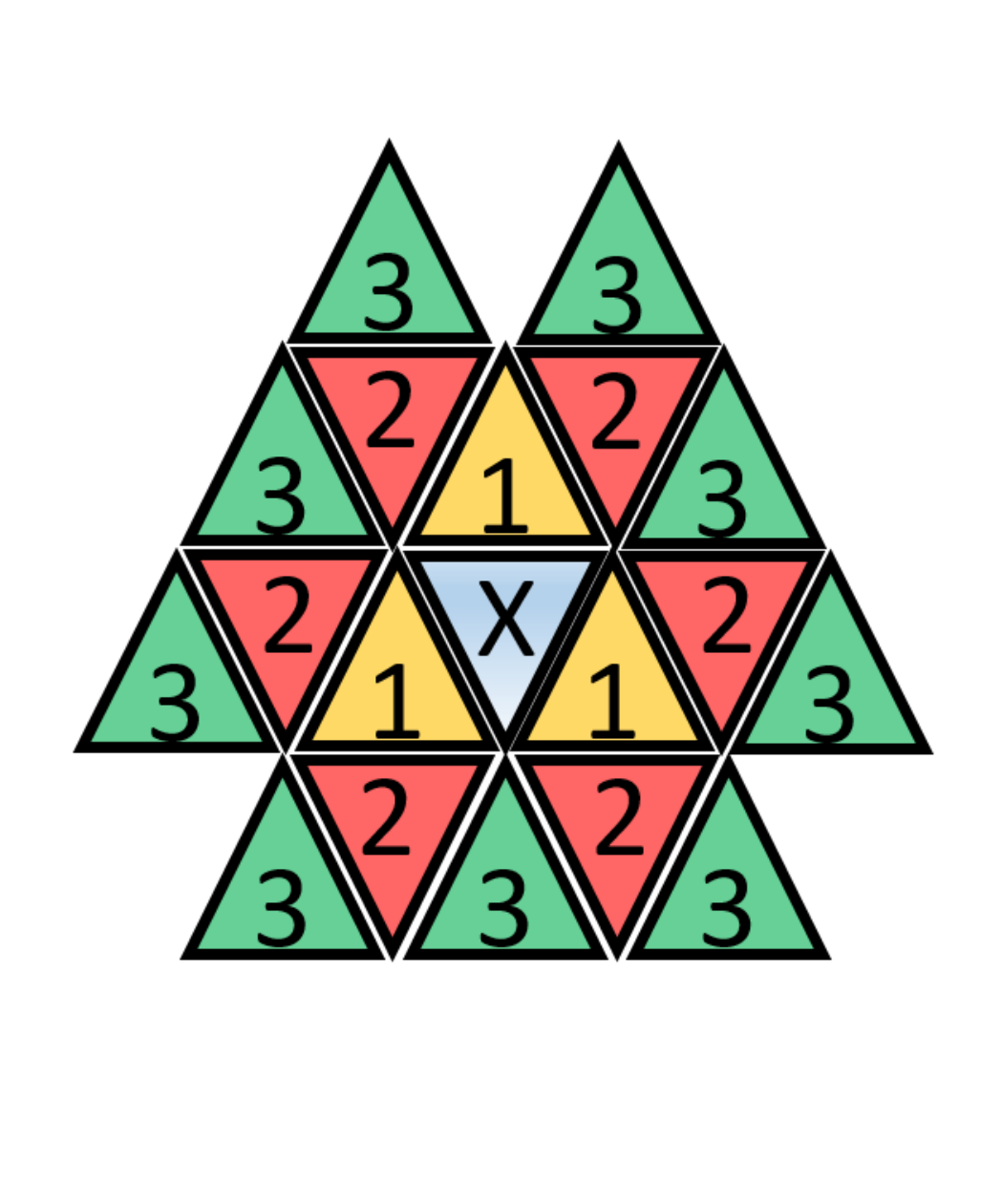}
			\label{fig:ex_a}} 
		\subfigure[][]{\includegraphics[height=4cm]{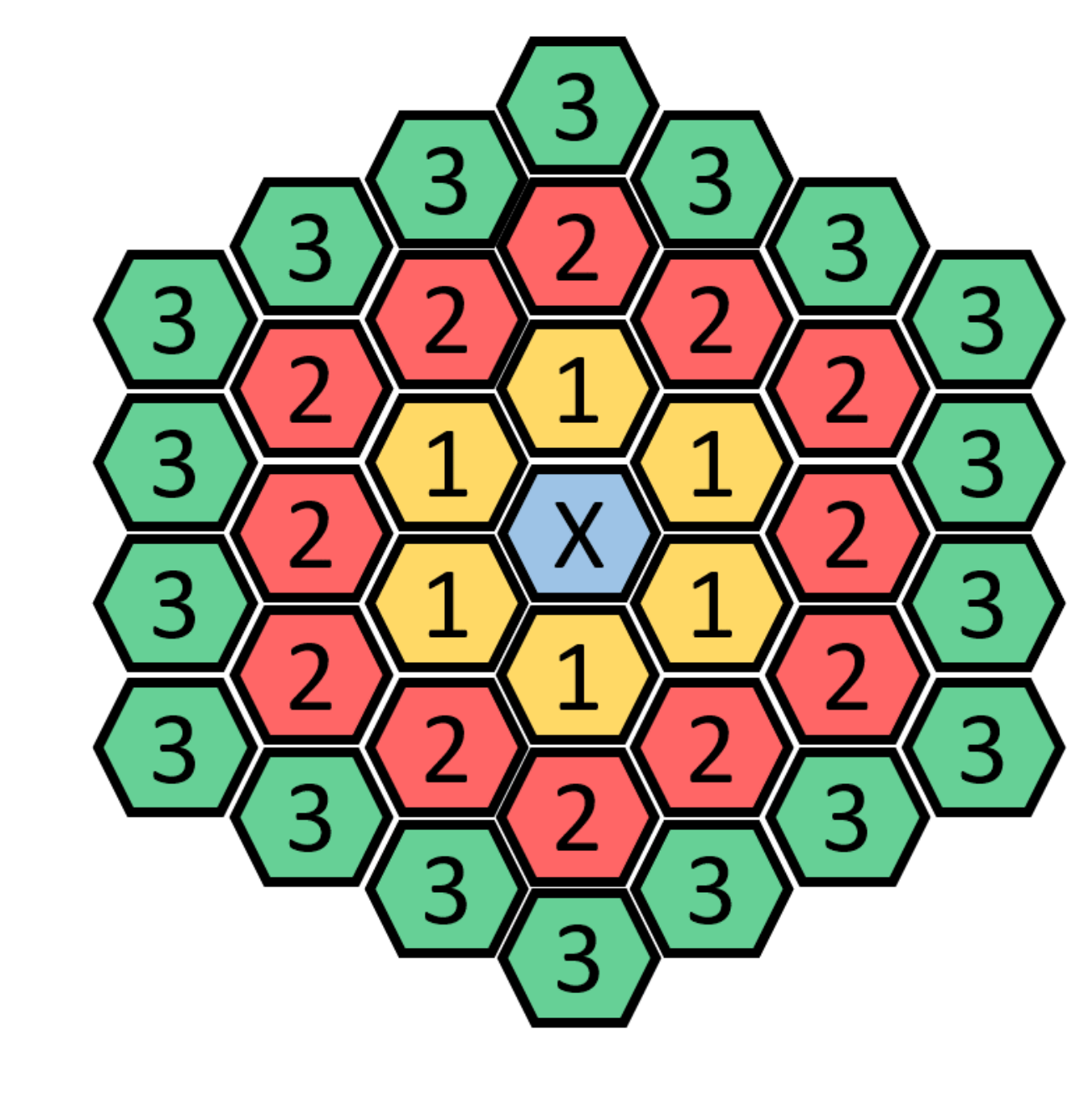}}
	\end{center}
	\caption{Schematic of agent pairs with the (a) triangular tessellation (b) hexagonal tessellation using the taxicab metric. Sites in yellow, red and green, labelled 1, 2 and 3 respectively, are distance-one, -two and -three neighbours from the blue site (labelled with X), respectively.}
	\label{fig:tm_tex} 
\end{figure}

\begin{subequations}
\label{eq:norm_tex}	
\begin{align}
s^{p}_{tri}(m)= 3 m \frac{L_x L_y}{2} \, ,\\ 
s^p_{hex}(m)= 3 m L_x L_y\, .
\end{align}	
\end{subequations}

Substituting equations \eqref{eq:norm_tex} into equation \eqref{normalisation_factor_1} we obtain the normalisations for the Triangle and Hexagon PCF, respectively, under periodic BC, namely 
\begin{subequations}
 \begin{align}
 \mathbb{E}\left[ \bar{c}^{p}_{tri}(m)\right]=\frac{3mN(N-1)}{2(L_xL_y-1)},
 \\\mathbb{E}\left[ \bar{c}^{p}_{hex}(m)\right]=\frac{3mN(N-1)}{L_xL_y-1}.
 \end{align}
 \end{subequations}
From these expressions one can obtain the formulae of $f_{tri}$ and $f_{hex}$ under periodic BC by using the definition \eqref{eq:PCF_definition}. The normalisations in the case of non-periodic BC are given in the SM Section \ref{SUPP-sec:norm_tri}. 
\subsection{Uniform Cube and Taxicab Cube PCF}
\label{sec:cube}

\begin{figure}[h!!]
	\begin{center}
		\subfigure[][]{\includegraphics[height=3cm]{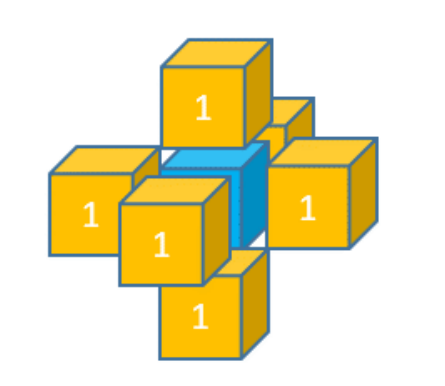}\label{fig:normsvis3D_TA} }
		\subfigure[][]{\includegraphics[height=3cm]{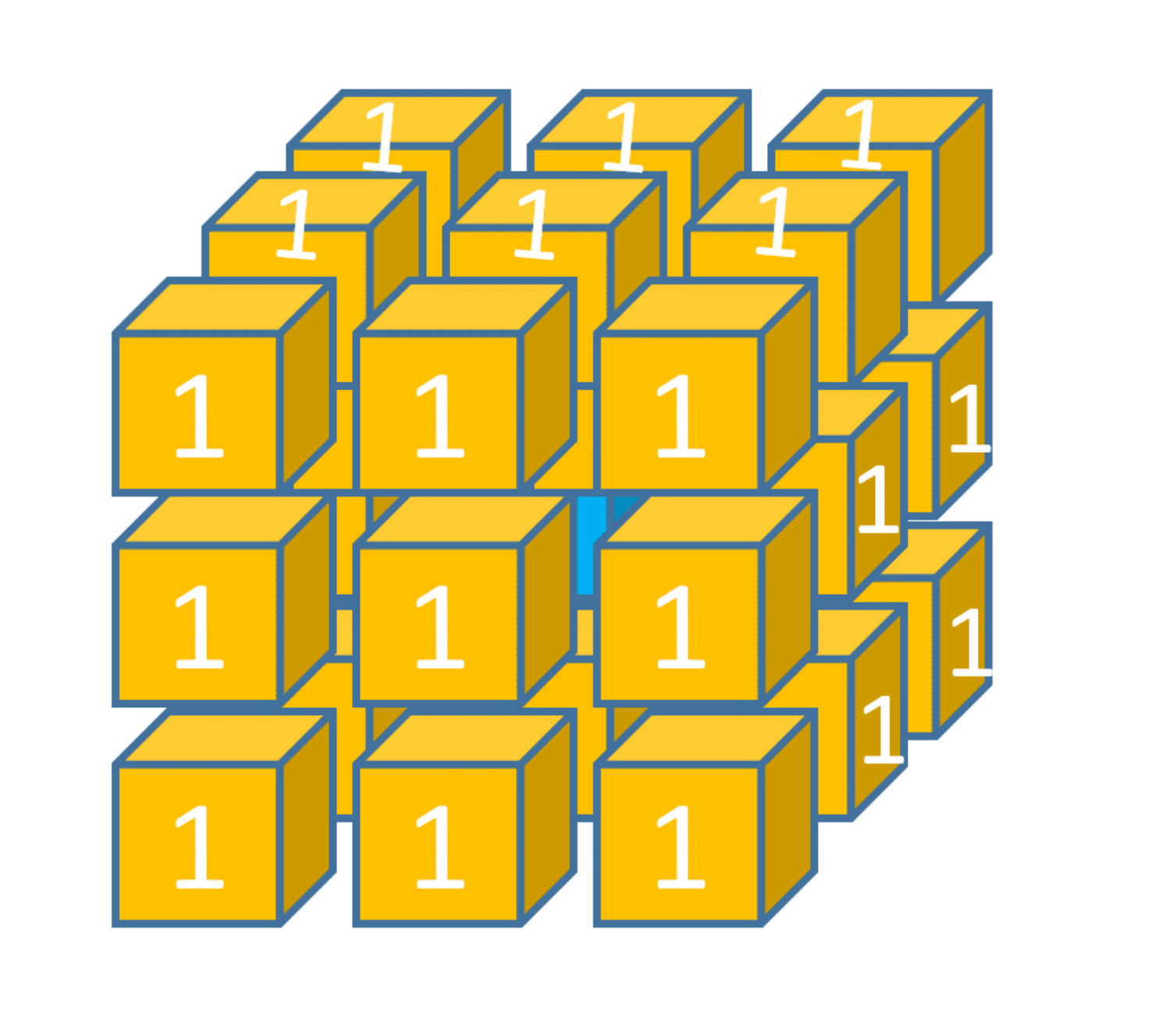}\label{fig:normsvis3D_UN} }
	\end{center}
	\caption{Schematic of agent pairs using (a) the taxicab metric (b) the uniform metric. Sites in yellow are defined to be distance-one neighbours from the site marked in blue.}
	\label{fig:normsvis3D} 
\end{figure}

We define a three dimensional $L_x \times L_y \times L_z$ cuboidal lattice with unit spacing. 
Using the taxicab and uniform metric respectively, the number of sites of distance $m$ from any given reference site is given by

\begin{subequations}
\label{eq:t_3d}
\begin{align}
t_{cube_1}(m)&=2(2m^2+1),\\
t_{cube_\infty}(m)&=2(12m^2+1) \, .
\end{align}
\end{subequations}
The proofs of the expressions \eqref{eq:t_3d} are omitted, but they can be obtained easily by induction on $m$. Examples of agent pairs for $m=1$ are given in Fig. \ref{fig:normsvis3D_TA} and \ref{fig:normsvis3D_UN} for the taxicab and uniform norms respectively. Using the same reasoning as in Section \ref{sec:PCF}, under periodic BC the normalisations for Taxicab and Uniform Cube PCFs, respectively, are as follows:
\begin{subequations}
\label{eq:norm_3d}
\begin{align}
s^{p}_{cube_1}(m)= &(2m^2+1) L_x L_y L_z \, ,
\\ s^{p}_{cube_\infty}(m)=& (12m^2+1)L_xL_yL_z \, .
\end{align}
\end{subequations}
For simplicity we refer the reader to Section \ref{SUPP-sec:norm_3D} of the SM for the normalisation factors for the cases with non-periodic BC.

\section{The General PCF}
\label{sec:ireg}
In this section we provide a comprehensive method for generating a PCF for any tessellation type, BC and metric but with the caveat of having a high computational cost.
This PCF is a valuable tool for irregular domain shapes and partitions although it can be used for any tessellation of any domain. 

\begin{figure}[h!!]
	\begin{center}
		\subfigure[][]{\includegraphics[height=5cm]{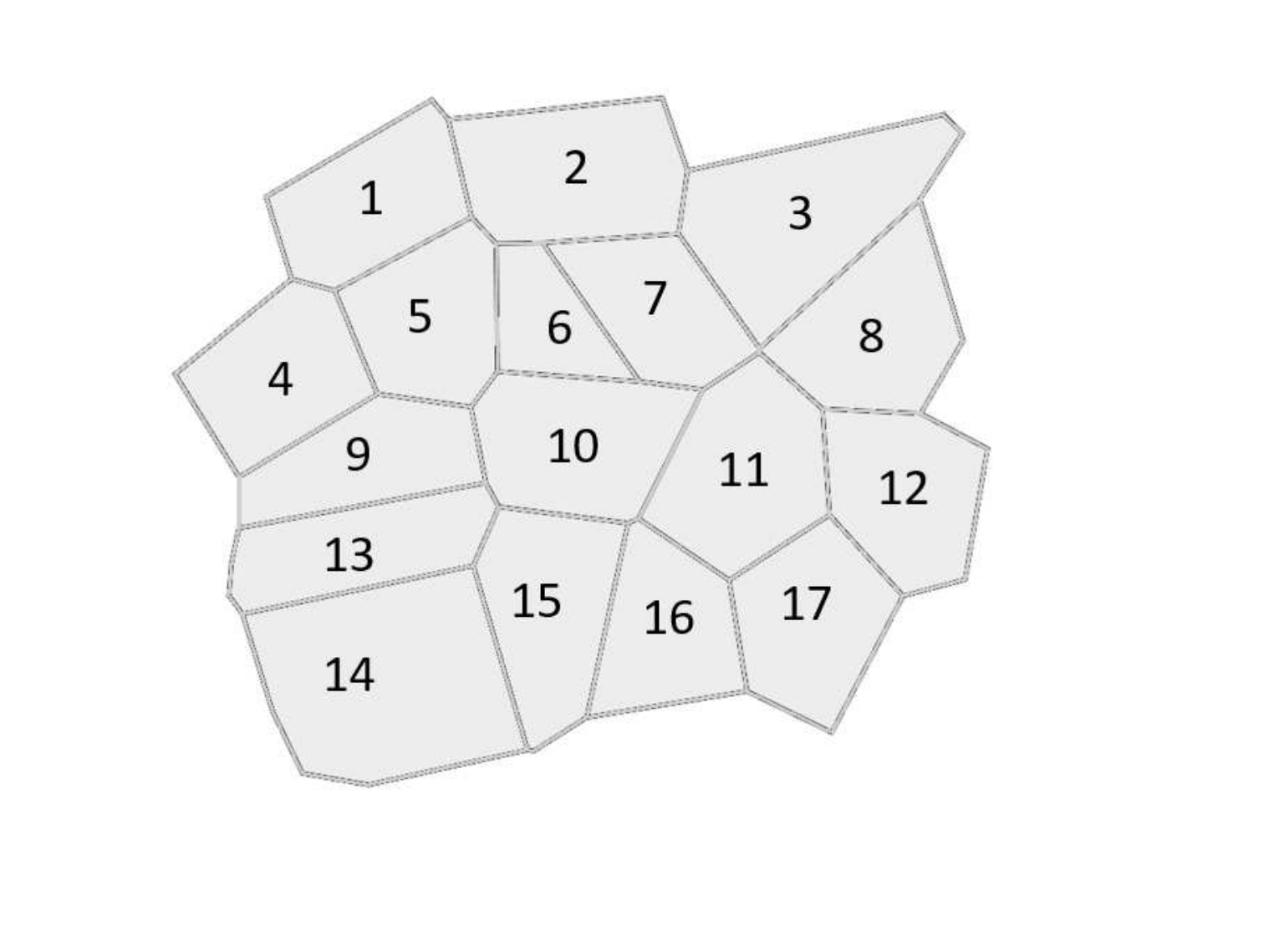}\label{fig:irregular_tesselation}}
		\subfigure[][]{\includegraphics[height=5cm]{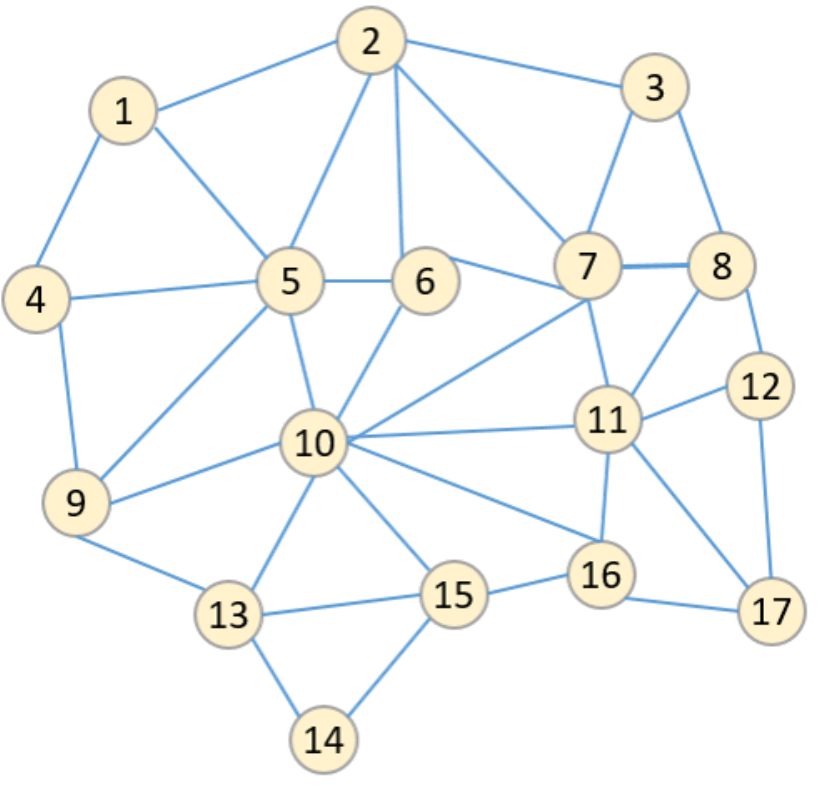}\label{fig:irregular_connectivity}}
	\end{center}
	\caption{An example of an irregular domain partition with its corresponding connectivity graph under the taxicab metric. Panel  \subref{fig:irregular_tesselation} shows A size 17 irregular lattice domain. Panel \subref{fig:irregular_connectivity} shows corresponding connectivity graph for the tessellation in \subref{fig:irregular_tesselation} under non-periodic BC using the taxicab metric for distances. 
}
	\label{tm_irr} 
\end{figure}

First, we consider a two-dimensional domain partitioned into $Z$ regions (or sites) with arbitrary shapes and sizes, each labelled with a number from 1 to $Z$. Fig. \ref{tm_irr} \subref{fig:irregular_tesselation} shows an example of an irregularly shaped domain partitioned in $Z=17$ regions. 
Given the domain, we choose a suitable metric. For the irregular lattice, which we consider in the following example, we consider the taxicab metric. This means that we define the distance, $m$, between two sites to be the minimum number of sites visited when starting at one site and moving consecutively through adjacent sites to the other. For example, in Fig. \ref{tm_irr} \subref{fig:irregular_tesselation} the sites 4 and 7 are at distance three. Similarly, adjacent sites are defined to be at distance one.

\label{sec:graph}

Having chosen and defined a suitable metric, we may now represent the connections between lattice sites as an undirected connectivity graph $G(V,E)$, where each vertex represents a lattice site and each edge connects vertices whose corresponding sites are distance-one neighbours. Using the taxicab metric, edges connect vertices whose corresponding sites are adjacent. Fig. \ref{tm_irr} \subref{fig:irregular_connectivity} shows an example of such an association (applied to the irregular lattice in Fig. \ref{tm_irr} \subref{fig:irregular_tesselation}) using the taxicab metric.

The corresponding adjacency matrix of graph $G$ is a $Z \times Z$ matrix defined as follows:

\begin{equation}
	A^G_{i,j} = \begin{cases} 1 & \text{for $(i,j) \in E$}, \\ 0  & \text{for $(i,j) \notin E$}. \end{cases}
\end{equation}
We use properties of the adjacency matrix to determine the number of sites at a given distance. In particular, we can compute $(A^G)^m$ whose entries $(A^G)^m_{i,j}$ are the number of walks of length $m$ from vertex $i$ to vertex $j$. To compute the minimum walk between two sites (and hence the distance between them) we produce the \textit{distance matrix} $D^G$. This is a $Z \times Z$ matrix defined as

\begin{equation}
	D^G_{i,j}= \begin{cases} \min \left\lbrace m \in \mathbb{N}^+ \, | \, (A^G)^m_{i,j} \neq 0 \right\rbrace & \text{for $i\neq j$}, \\ 
 0 &	\text{for $i=j$}. \end{cases}
\label{eq:distance_matrix}
\end{equation}
Notice that, each entry, $D^G_{i,j}$, denotes the distance between the vertices $i$ and $j$ on graph $G$ and hence on the original lattice.

Given the distance matrix of the domain, $D^G$, and the set of the occupied sites $ M \subseteq V $, the PCF of the system can be computed as follow. The number of pairs of agents at distance $m$ for a general metric $d$ is given by:
\begin{equation}
\label{eq:num_PCF_irr}
	c_{d}(m)=\frac{1}{2} \left| \lbrace (i,j)\in M \times M \, |\, D^G_{i,j}=m \rbrace\right| \, .
\end{equation}
Similarly, we can express the number of pairs of sites at distance $m$ as 
\begin{equation}
\label{eq:num_PCF_irr}
	s_{d}(m)=\frac{1}{2} \left| \lbrace (i,j)\in V \times V \, |\, D^G_{i,j}=m \rbrace\right| \, .
\end{equation}

To compute the normalisation factor, denote the total number of agents as $N=|M|$ and hence using the same argument as in Section \ref{sec:PCF} we can write
\begin{equation}
\mathbb{E}[\bar{c}_d (m)]=\left(\frac{N}{Z}\right)\left(\frac{N-1}{Z-1}\right) s_d (m).  
\label{eq:irr_normalisation_factor}
\end{equation} 
The General PCF is then defined by combining equations \eqref{eq:num_PCF_irr} and \eqref{eq:irr_normalisation_factor} in equation \eqref{eq:PCF_definition}.
Notice that the computation of the normalisation for the General PCF can be computationally expensive. This is because the computation of $D^G$ involves calculating powers of matrices of size $Z \times Z$, where $Z$ is often large. In particular, the cost of computing the normalisation of the General PCF is $\mathcal{O}(Z^3 m_{max})$ in which $m_{max}$ is the maximum value $m$ for which the PCF is computed. For this reason, the General PCF is better reserved for cases in which the expression for the normalisation factor cannot be computed analytically, unlike in Sections \ref{sec:PCF}, \ref{sec:tria-hex} and \ref{sec:cube}, although it can, of course, be used even if an analytical formula is available.  

\begin{figure}[h!!!!!!]
	\begin{center}
		\subfigure[][]{\includegraphics[ width=0.38 \columnwidth]{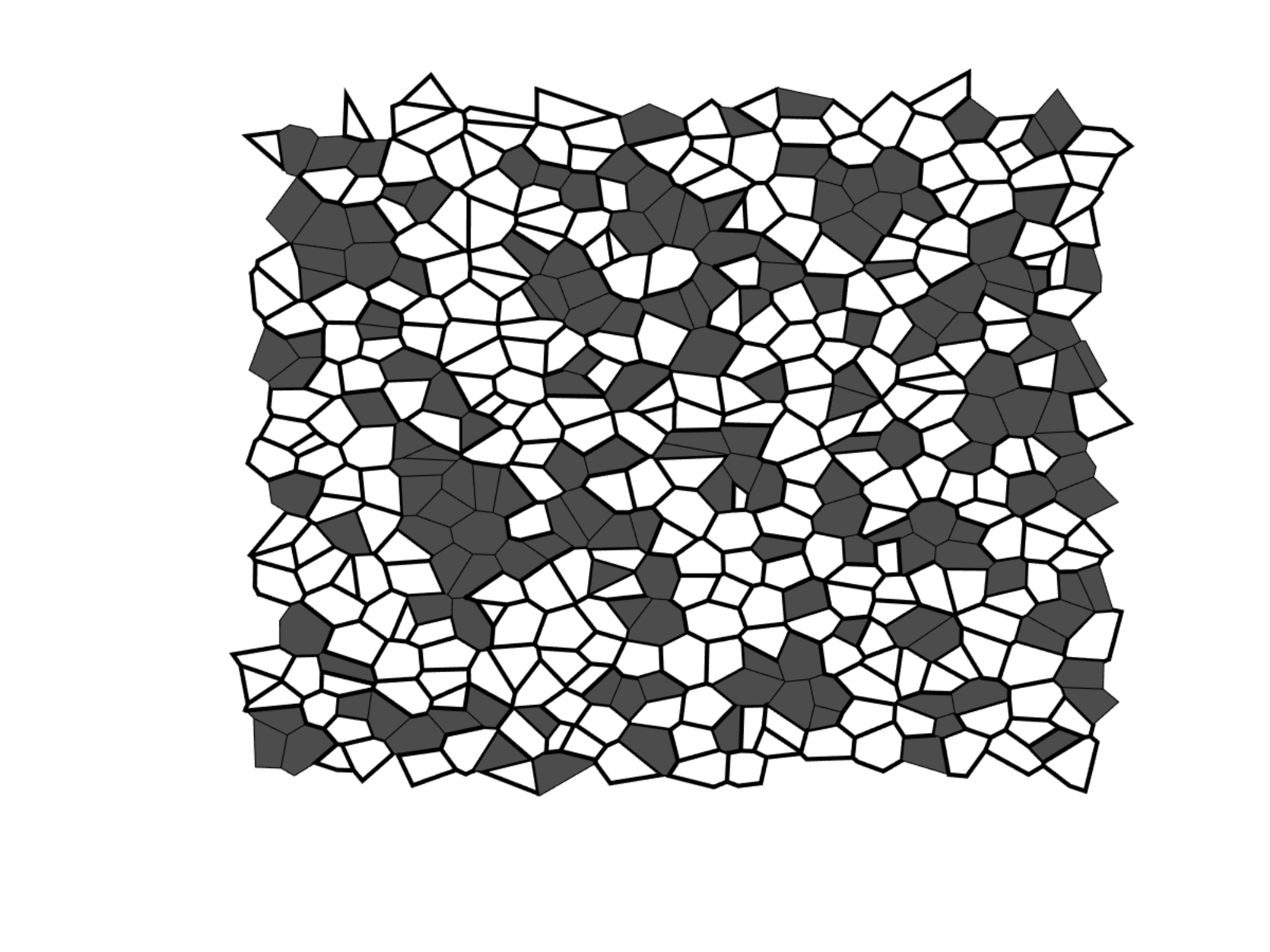}\label{fig:ex_irr_unif_lattice}}
		\subfigure[][]{\includegraphics[width=0.32 \columnwidth]{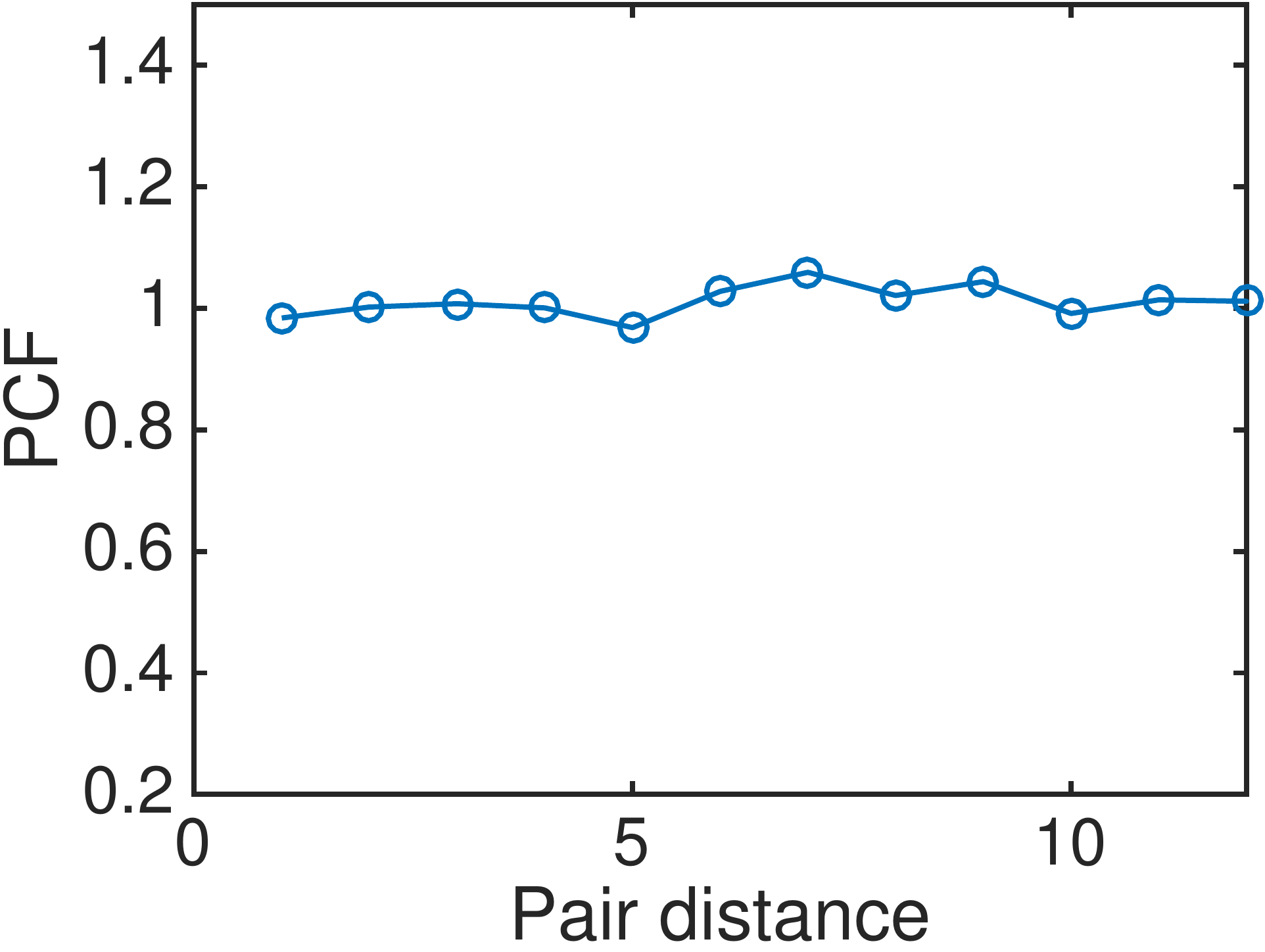}\label{fig:ex_irr_unif_PCF}}\\
		\subfigure[][]{\includegraphics[ width=0.38 \columnwidth]{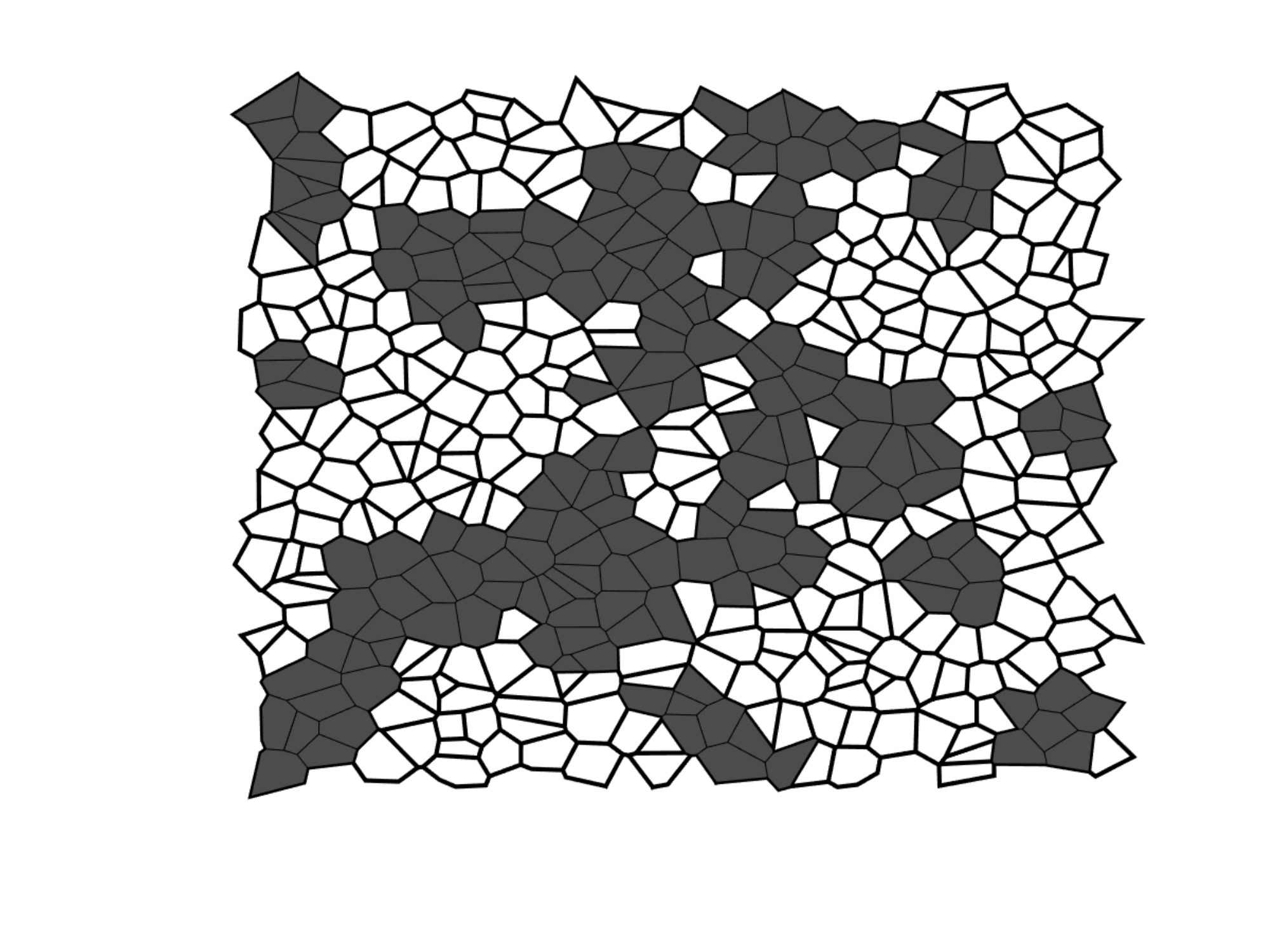}\label{fig:ex_irr_agg_lattice}}
		\subfigure[][]{\includegraphics[width=0.32 \columnwidth]{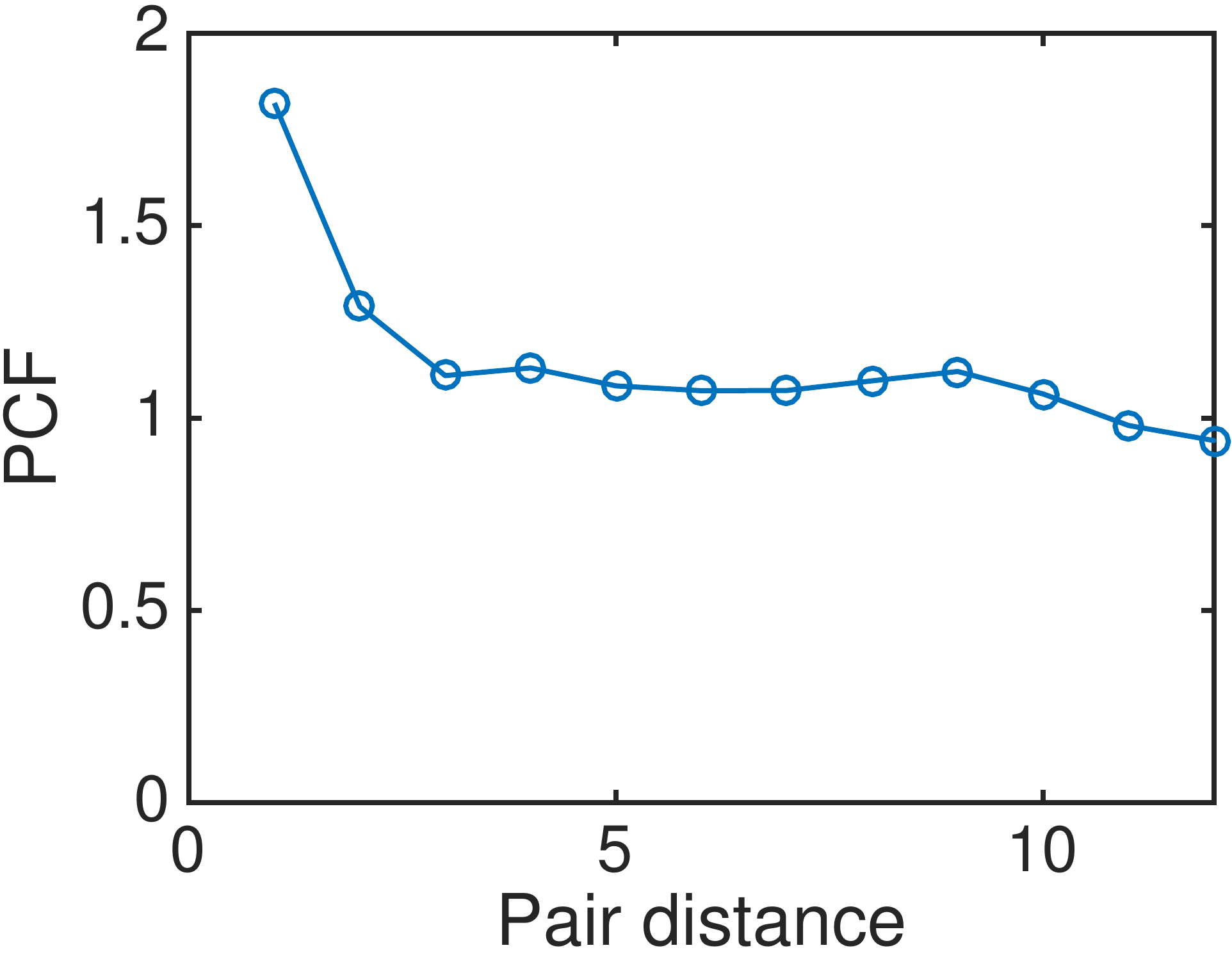}\label{fig:ex_irr_agg_PCF}}\\
		\subfigure[][]{\includegraphics[ width=0.38 \columnwidth]{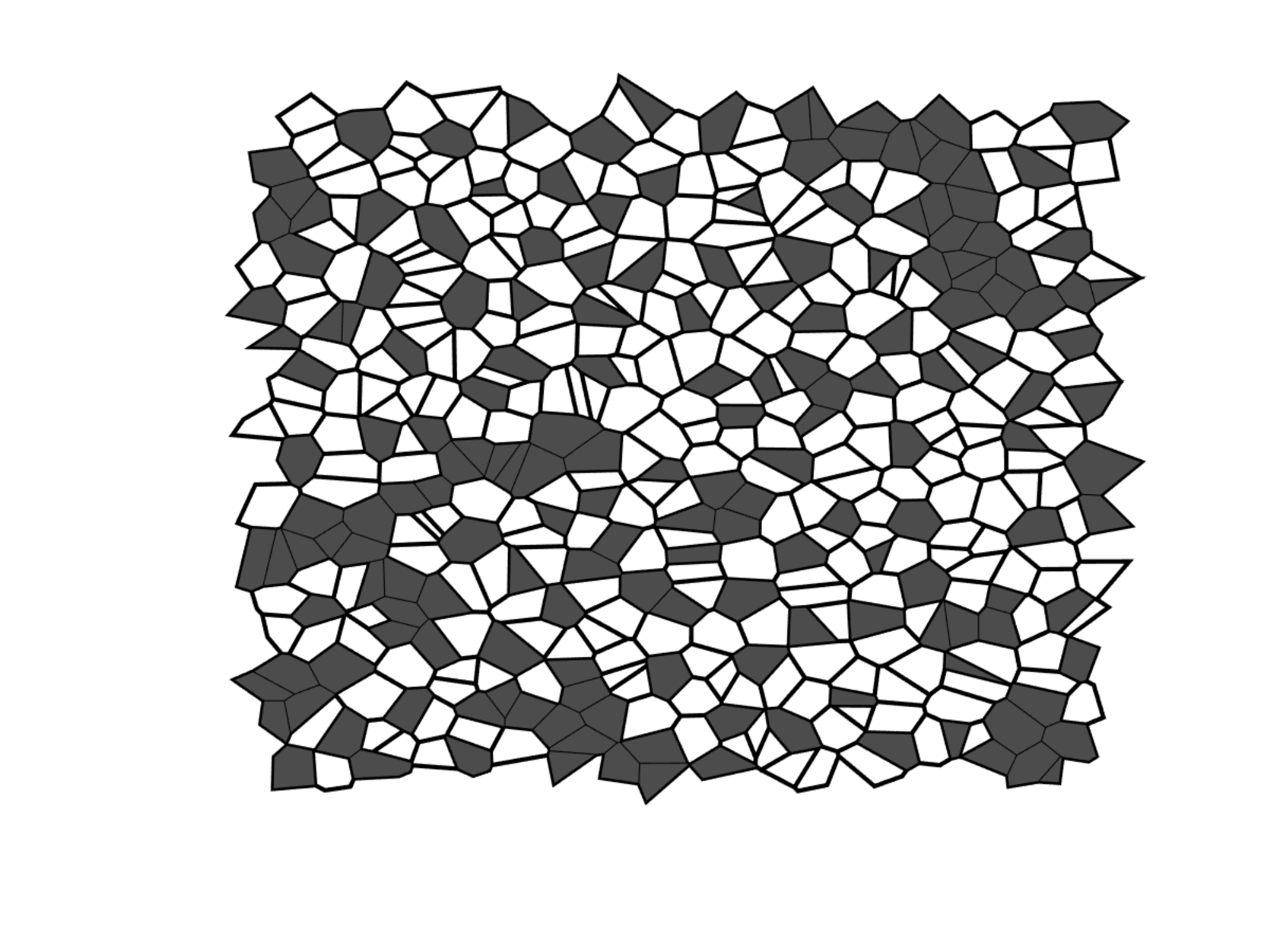}\label{fig:ex_irr_seg_lattice}}
		\subfigure[][]{\includegraphics[width=0.32 \columnwidth]{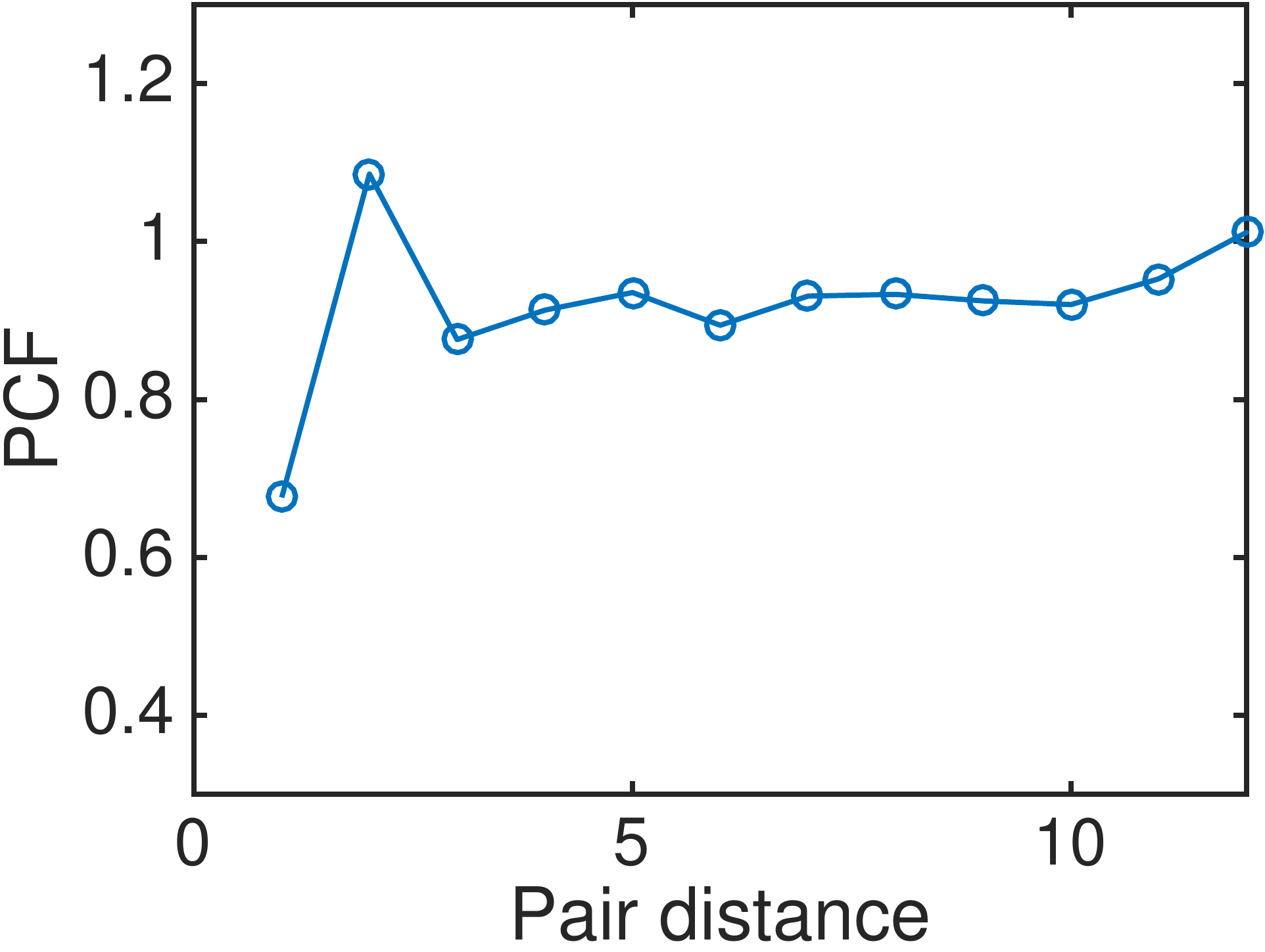}\label{fig:ex_irr_seg_PCF}}
	\end{center} 
	\caption{Examples of spatial correlation analysis on an irregular domain. Panels \subref{fig:ex_irr_unif_lattice}, \subref{fig:ex_irr_agg_lattice} and \subref{fig:ex_irr_seg_lattice} show three examples of irregular lattices populated with density $0.4$ with agents (grey sites). In panel \subref{fig:ex_irr_unif_lattice}, agents are displaced uniformly at random (no spatial correlation). In panel \subref{fig:ex_irr_agg_lattice}, agents are in a strong form of aggregation, while in panel \subref{fig:ex_irr_seg_lattice} agents are displaced in a segregate manner. Panels \subref{fig:ex_irr_unif_PCF}, \subref{fig:ex_irr_agg_PCF} and \subref{fig:ex_irr_seg_PCF} are the corresponding General PCF evaluations for panels \subref{fig:ex_irr_unif_lattice}, \subref{fig:ex_irr_agg_lattice} and \subref{fig:ex_irr_seg_lattice}, respectively.}
	\label{fig:ex_irr} 
\end{figure}

In Fig. \ref{fig:ex_irr} we apply the General PCF to three examples of agent-based systems on an irregular lattice.  In all three examples, the irregular tessellation is the Voronoi partition based on a set of randomly distributed points. The randomised points are obtained by starting with a square lattice, with lattice size $\Delta$, and perturbing the coordinates of each point $(x_i,y_i)$ to $(x_i+\delta^x_i, y +\delta^y_i)$ with each $\delta^x_i,\delta^y_i$ chosen uniformly at random in the interval $[-\frac{\Delta}{2},\frac{\Delta}{2}]$. 

In the first example, Fig. \ref{fig:ex_irr} \subref{fig:ex_irr_unif_lattice}, $N$ lattice sites selected uniformly at random to be occupied (grey sites). By eye, several larger clusters of occupied and unoccupied lattice sites are evident, indicating that there may be spatial correlation. Fig. \ref{fig:ex_irr} \subref{fig:ex_irr_unif_PCF} shows the corresponding General PCF. The values are close to unity, correctly identifying that there is, indeed, no spatial correlation in the system. This highlights the importance of accurate quantitative methods for determining spatial correlation rather than a reliance on \textit{ad hoc} judgements. 

In the other two cases we test our PCF with examples of strong spatial correlation. In Fig. \ref{fig:ex_irr} \subref{fig:ex_irr_agg_lattice} we consider a system in an aggregated state. To generate such a configuration, we start with an empty domain, and select an empty site uniformly at random. We then place an agent in this site and all of its neighbouring sites (if they are not already occupied). We repeat this process until we reach density $0.4$. The process leads to a strong form of aggregation which the General PCF, in Fig. \ref{fig:ex_irr} \subref{fig:ex_irr_agg_PCF}, correctly identifies. Finally, in Fig. \ref{fig:ex_irr} (e) we consider a system in a segregated state. To generate such a configuration, we start with a fully populated domain, and repeatedly select at random an occupied site. We remove all agents occupying adjacent sites but leave the initially selected site occupied. The process ends when density $0.4$ is reached. This mechanism generates a system which is unlikely to have adjacent sites occupied and more likely to have agents displaced at distance two. The corresponding PCF, shown in Fig. \ref{fig:ex_irr} \subref{fig:ex_irr_seg_PCF}, correctly captures both of these features: the PCF have value $0.68$ at pair distance one which implies negative correlation at the shortest distance, and value $1.1$ at pair distance two, highlighting the positive correlation at slightly larger distances.

 Note that, since the normalisation in the General PCF can give an exact value for the number of sites at certain pair distances, we used this method to check the normalisation factors for all of our previously proposed PCFs. In all cases the results confirmed the analytical expressions given in Sections \ref{sec:PCF}, \ref{sec:tria-hex}, \ref{sec:cube} and in the SM.

\section{Conclusions}
\label{sec:conclusions}

In this paper we have developed a set of new tools to study spatial correlation on discrete domains. We derived two discrete pair-correlation functions for an exclusion process on a square lattice: the Square Uniform and Square Taxicab PCF. We applied our PCFs to patterns observed in nature and to computational simulations and showed that our PCF can not only distinguish and quantify different types of correlation but also improved upon previous on-lattice PCFs. For example, we showed that our PCF was normalised correctly, unlike the Annular PCF, and was able to identify and quantify anisotropic patterns such as the chessboard or diagonal stripes that the Rectilinear PCF \cite{binder2013qss} missed. Furthermore, we highlighted how different measures of distance, taxicab and uniform, can lead to different quantifications of spatial correlation.

We extended the calculation of appropriate PCFs to deal with exclusion processes on the other regular spatial tessellations in two dimensions as well as the cubic lattice in three dimensions. We derived the Triangle PCF, Hexagon PCF, Cube Taxicab PCF and Cube Uniform PCF. These are the first PCFs defined specifically for these discrete lattice types. 
Finally, we derived a comprehensive PCF for any kind of discrete domain, BC and metric
which we referred to as the General PCF. The method can be computationally expensive, however, it allows complete freedom in defining a suitable PCF for more complex cases, including those for which recognising spatial correlation by eye becomes less intuitive. 

All of our PCFs are designed for a single species of agents.  However, in many applications the agents are divided into multiple species and it can be important to distinguish between different types of spatial correlations: either within agents of the same species (auto-correlation) or comparing the position of agents of different species (cross-correlation). \citet{dini2017uib} have recently investigated  correlation in multiple species by using the Rectilinear PCF. We believe a similar approach to that of \citet{dini2017uib} can be applied to our new isotropic PCFs to quantify heterogenous correlations, yet it lies beyond the scope of this paper and, as such, we will tackle it in a future publication. 
The new isotropic PCFs which we defined in this paper will be important in further studies and applications. In particular, our novel functions can improve previous studies (see \citep{johnston2014isa}, for example) which have used PCF as an efficient summary statistics to infer model parameters. 

\ack
The authors would like to thank the CMB-CNCB preprint club for constructive and helpful comments on a preprint of this paper. J.P.O. acknowledges support from the SWBio DTP.

\bibliography{database.bib}

\begin{thebibliography}{42}
\providecommand{\natexlab}[1]{#1}
\providecommand{\url}[1]{\texttt{#1}}
\expandafter\ifx\csname urlstyle\endcsname\relax
  \providecommand{\doi}[1]{doi: #1}\else
  \providecommand{\doi}{doi: \begingroup \urlstyle{rm}\Url}\fi

\bibitem[Agnew et~al.(2014)Agnew, Green, Brown, Simpson, and
  Binder]{agnew2014dmc}
D.J.G. Agnew, J.E.F. Green, T.M. Brown, M.J. Simpson, and B.J. Binder.
\newblock Distinguishing between mechanisms of cell aggregation using
  pair-correlation functions.
\newblock \emph{J. Theor. Biol.}, 352:\penalty0 16--23, 2014.

\bibitem[Bahcall and Soneira(1983)]{bahcall1983scf}
N.A. Bahcall and R.M. Soneira.
\newblock The spatial correlation function of rich clusters of galaxies.
\newblock \emph{Astrophys. J.}, 270:\penalty0 20--38, 1983.

\bibitem[Baker et~al.(2010)Baker, Yates, and Erban]{baker2010fmm}
R.E. Baker, C.A. Yates, and R.~Erban.
\newblock From microscopic to macroscopic descriptions of cell migration on
  growing domains.
\newblock \emph{Bull. Math. Biol.}, 72\penalty0 (3):\penalty0 719--762, 2010.

\bibitem[Ball(2015)]{ball2015fpm}
P.~Ball.
\newblock Forging patterns and making waves from biology to geology: a
  commentary on turing (1952) ``the chemical basis of morphogenesis''.
\newblock \emph{Phil. Trans. R. Soc. B}, 370\penalty0 (1666):\penalty0
  20140218, 2015.

\bibitem[Bhide and Yashonath(1999)]{bhide1999dsc}
S.~Y. Bhide and S.~Yashonath.
\newblock Dependence of the self-diffusion coefficient on the sorbate
  concentration: A two-dimensional lattice gas model with and without
  confinement.
\newblock \emph{J. Chem. Phys.}, 111\penalty0 (4):\penalty0 1658--1667, 1999.

\bibitem[Binder and Simpson(2013)]{binder2013qss}
B.J. Binder and M.J. Simpson.
\newblock Quantifying spatial structure in experimental observations and
  agent-based simulations using pair-correlation functions.
\newblock \emph{Phys. Rev. E}, 88:\penalty0 022705, 2013.

\bibitem[Binder and Simpson(2015)]{binder2015sap}
B.J. Binder and M.J. Simpson.
\newblock Spectral analysis of pair-correlation bandwidth: application to cell
  biology images.
\newblock \emph{R. Soc. Open Sci.}, 2\penalty0 (2):\penalty0 140494, 2015.

\bibitem[Binder et~al.(2015)Binder, Sundstrom, Gardner, Jiranek, and
  Oliver]{binder2015qtf}
B.J. Binder, J.F. Sundstrom, J.M. Gardner, V.~Jiranek, and S.G. Oliver.
\newblock Quantifying two-dimensional filamentous and invasive growth spatial
  patterns in yeast colonies.
\newblock \emph{Plos Comput. Biol.}, 11\penalty0 (2):\penalty0 e1004070, 2015.

\bibitem[Binny et~al.(2016)Binny, James, and Plank]{binny2016ccb}
R.N. Binny, A.~James, and M.J. Plank.
\newblock Collective cell behaviour with neighbour-dependent proliferation,
  death and directional bias.
\newblock \emph{Bull. Math. Biol.}, 2016.

\bibitem[Bone et~al.(1986)Bone, Bachor, and Sandeman]{bone1986fpa}
D.J. Bone, H.A. Bachor, and R.J. Sandeman.
\newblock Fringe-pattern analysis using a 2-d fourier transform.
\newblock \emph{Appl. Optics.}, 25\penalty0 (10):\penalty0 1653--1660, 1986.

\bibitem[Bourke and Franks(1995)]{bourke1995sea}
A.F.G. Bourke and N.R. Franks.
\newblock \emph{Social evolution in ants}.
\newblock Princeton University Press, 1995.

\bibitem[Browning et~al.(2017)Browning, McCue, and Simpson]{browning2017bca}
A.P. Browning, S.~McCue, and M.J. Simpson.
\newblock A bayesian computational approach to explore the optimal duration of
  a cell proliferation assay.
\newblock \emph{bioRxiv}, page 147678, 2017.

\bibitem[Cavagna et~al.(2010)Cavagna, Cimarelli, Giardina, Parisi, Santagati,
  Stefanini, and Viale]{cavagna2010sfc}
A.~Cavagna, A.~Cimarelli, I.~Giardina, G.~Parisi, R.~Santagati, F.~Stefanini,
  and M.~Viale.
\newblock Scale-free correlations in starling flocks.
\newblock \emph{{Proc. Natl. Acad. Sci. USA.}}, 107\penalty0 (26):\penalty0
  11865--11870, 2010.

\bibitem[Chalub et~al.(2006)Chalub, Dolak-Struss, Markowich, Oelz, Schmeiser,
  and Soreff]{chalub2006mhc}
F.~Chalub, Y.~Dolak-Struss, P.~Markowich, D.~Oelz, C.~Schmeiser, and A.~Soreff.
\newblock Model hierarchies for cell aggregation by chemotaxis.
\newblock \emph{Math. Models Methods Appl. Sci.}, 16:\penalty0 1173--1197,
  2006.

\bibitem[Deutsch and Dormann(2007)]{deutsch2007cam}
A.~Deutsch and S.~Dormann.
\newblock \emph{Cellular automaton modeling of biological pattern formation:
  characterization, applications, and analysis}.
\newblock Springer Science \& Business Media, 2007.

\bibitem[Diggle et~al.(1976)Diggle, Besag, and Gleaves]{diggle1976sas}
P.J. Diggle, J.~Besag, and T.J. Gleaves.
\newblock Statistical analysis of spatial point patterns by means of distance
  methods.
\newblock \emph{Biometrics}, 1976.

\bibitem[Donev et~al.(2005)Donev, Torquato, and Stillinger]{donev2005pcf}
A.~Donev, S.~Torquato, and F.H. Stillinger.
\newblock Pair correlation function characteristics of nearly jammed disordered
  and ordered hard-sphere packings.
\newblock \emph{Phys. Rev. E}, 71\penalty0 (1):\penalty0 011105, 2005.

\bibitem[Friedman et~al.(1985)Friedman, Rigel, and Kopf]{friedman1985edm}
R.J. Friedman, D.S. Rigel, and A.W. Kopf.
\newblock Early detection of malignant melanoma: The role of physician
  examination and self-examination of the skin.
\newblock \emph{CA. Cancer J. Clin.}, 35\penalty0 (3):\penalty0 130--151, 1985.

\bibitem[Gavagnin et~al.(2018)Gavagnin, Owen, and Yates]{gavagnin2018pcf}
E.~Gavagnin, J.P. Owen, and C.A. Yates.
\newblock Pair correlation functions for identifying spatial correlation in
  discrete domains.
\newblock \emph{Phys. Rev. E}, 97:\penalty0 062104, 2018.

\bibitem[Green et~al.(2010)Green, Waters, Whiteley, Edelstein-Keshet,
  Shakesheff, and Byrne]{green2010nlm}
J.E.F. Green, S.L. Waters, J.P. Whiteley, L.~Edelstein-Keshet, K.M. Shakesheff,
  and H.M. Byrne.
\newblock Non-local models for the formation of hepatocyte--stellate cell
  aggregates.
\newblock \emph{J. Theor. Biol.}, 267\penalty0 (1):\penalty0 106--120, 2010.

\bibitem[Hackett-Jones et~al.(2012)Hackett-Jones, Davies, Binder, and
  Landman]{hackett2012gis}
E.J. Hackett-Jones, K.J. Davies, B.J. Binder, and K.A. Landman.
\newblock Generalized index for spatial data sets as a measure of complete
  spatial randomness.
\newblock \emph{Phys. Rev. E}, 85\penalty0 (6):\penalty0 061908, 2012.

\bibitem[Hinde et~al.(2010)Hinde, Cardarelli, Digman, and
  Gratton]{hinde2010vpc}
E.~Hinde, F.~Cardarelli, M.A. Digman, and E.~Gratton.
\newblock In vivo pair correlation analysis of egfp intranuclear diffusion
  reveals dna-dependent molecular flow.
\newblock \emph{Proc. Natl. Acad. Sci. U.S.A.}, 107\penalty0 (38):\penalty0
  16560--16565, 2010.

\bibitem[Illian et~al.(2008)Illian, Penttinen, Stoyan, and
  Stoyan]{illian2008sam}
J.~Illian, A.~Penttinen, H.~Stoyan, and D.~Stoyan.
\newblock \emph{Statistical analysis and modelling of spatial point patterns}.
\newblock John Wiley \& Sons, 2008.

\bibitem[Johnston et~al.(2014)Johnston, Simpson, McElwain, Binder, and
  Ross]{johnston2014isa}
S.T. Johnston, M.J. Simpson, D.L.S. McElwain, B.J. Binder, and J.V. Ross.
\newblock Interpreting scratch assays using pair density dynamics and
  approximate bayesian computation.
\newblock \emph{Op. bio.}, 4\penalty0 (9):\penalty0 140097, 2014.

\bibitem[Keeling(1999)]{keeling1999els}
M.J. Keeling.
\newblock The effects of local spatial structure on epidemiological invasions.
\newblock \emph{Proc. R. Soc., Ser. B, Biol. Sc., Lond.}, 266\penalty0
  (1421):\penalty0 859--867, 1999.

\bibitem[Kondo(2017)]{kondo2017ukt}
S.~Kondo.
\newblock An updated kernel-based turing model for studying the mechanisms of
  biological pattern formation.
\newblock \emph{J. Theor. Biol.}, 414:\penalty0 120--127, 2017.

\bibitem[Murakawa and Togashi(2015)]{murakawa2015cmc}
H.~Murakawa and H.~Togashi.
\newblock Continuous models for cell-cell adhesion.
\newblock \emph{J. Theor. Biol.}, 374:\penalty0 1--12, 2015.

\bibitem[Murray(2007)]{murray2007mbi}
J.D. Murray.
\newblock \emph{Mathematical biology: I. An introduction}, volume~17.
\newblock Springer Science \& Business Media, 2007.

\bibitem[Othmer and Stevens(1997)]{othmer1997abc}
H.G. Othmer and A.~Stevens.
\newblock Aggregation, blowup, and collapse: the {A}{B}{C}'s of taxis in
  reinforced random walks.
\newblock \emph{{SIAM} J. Appl. Math.}, 57\penalty0 (4):\penalty0 1044--1081,
  1997.

\bibitem[Ouyang and Swinney(1991)]{ouyang1991tuh}
Q.I. Ouyang and H.L. Swinney.
\newblock Transition from a uniform state to hexagonal and striped turing
  patterns.
\newblock \emph{Nature}, 352\penalty0 (6336):\penalty0 610--612, 1991.

\bibitem[Raghib et~al.(2011)Raghib, Hill, and Dieckmann]{raghib2011mme}
M.~Raghib, N.A. Hill, and U.~Dieckmann.
\newblock A multiscale maximum entropy moment closure for locally regulated
  space--time point process models of population dynamics.
\newblock \emph{J. Math. Biol.}, 62\penalty0 (5):\penalty0 605--653, 2011.

\bibitem[Ross et~al.(2015)Ross, Yates, and Baker]{ross2015icc}
R.J.H. Ross, C.A. Yates, and R.E. Baker.
\newblock Inference of cell--cell interactions from population density
  characteristics and cell trajectories on static and growing domains.
\newblock \emph{Math. Biosci.}, 264:\penalty0 108--118, 2015.

\bibitem[Simpson et~al.(2009)Simpson, Landman, and Hughes]{simpson2009mss}
M.J. Simpson, K.A. Landman, and B.D. Hughes.
\newblock {Multi-species simple exclusion processes}.
\newblock \emph{Phys. A}, 388\penalty0 (4):\penalty0 399--406, 2009.

\bibitem[Simpson et~al.(2013)Simpson, Binder, Haridas, Wood, Treloar, McElwain,
  and Baker]{simpson2013emi}
M.J. Simpson, B.J. Binder, P.~Haridas, B.~K Wood, K.K. Treloar, D.~McElwain,
  and R.E. Baker.
\newblock Experimental and modelling investigation of monolayer development
  with clustering.
\newblock \emph{Bull. Math. Biol.}, 75\penalty0 (5):\penalty0 871--889, 2013.

\bibitem[Steinberg(1996)]{steinberg1996adh}
M.~S. Steinberg.
\newblock Adhesion in development: an historical overview.
\newblock \emph{Dev. Biol.}, 180\penalty0 (2):\penalty0 377--388, 1996.

\bibitem[Stevens(2000)]{stevens2000sca}
A.~Stevens.
\newblock A stochastic cellular automaton modeling gliding and aggregation of
  myxobacteria.
\newblock \emph{SIAM J. Appl. Math.}, 61\penalty0 (1):\penalty0 172--182, 2000.

\bibitem[Takeda et~al.(1982)Takeda, Ina, and Kobayashi]{takeda1982ftm}
M.~Takeda, H.~Ina, and S.~Kobayashi.
\newblock Fourier-transform method of fringe-pattern analysis for
  computer-based topography and interferometry.
\newblock \emph{JosA}, 72\penalty0 (1):\penalty0 156--160, 1982.

\bibitem[Thomas et~al.(2005)Thomas, Bennett, Thomson, and
  Shakesheff]{thomas2005hsc}
R.J. Thomas, A.~Bennett, B.~Thomson, and K.M. Shakesheff.
\newblock Hepatic stellate cells on poly (dl-lactic acid) surfaces control the
  formation of 3d hepatocyte co-culture aggregates in vitro.
\newblock \emph{ecells \& materials}, 11:\penalty0 16--26, 2005.

\bibitem[Treloar et~al.(2013)Treloar, Simpson, Haridas, Manton, Leavesley,
  McElwain, and Baker]{treloar2013mtd}
K.K. Treloar, M.J. Simpson, P.~Haridas, K.J. Manton, D.I. Leavesley, D.S.
  McElwain, and R.E. Baker.
\newblock Multiple types of data are required to identify the mechanisms
  influencing the spatial expansion of melanoma cell colonies.
\newblock \emph{BMC Syst. Biol.}, 7\penalty0 (1):\penalty0 137, 2013.

\bibitem[Weinstock(2000)]{weinstock2000edm}
M.A. Weinstock.
\newblock Early detection of melanoma.
\newblock \emph{JAMA}, 284\penalty0 (7):\penalty0 886--889, 2000.

\bibitem[Yates et~al.(2012)Yates, Baker, Erban, and Maini]{yates2012gfm}
C.A. Yates, R.E. Baker, R.~Erban, and P.K. Maini.
\newblock Going from microscopic to macroscopic on nonuniform growing domains.
\newblock \emph{Phys. Rev. E}, 86\penalty0 (2):\penalty0 021921, 2012.

\bibitem[Young et~al.(2001)Young, Roberts, and Stuhne]{young2001rpc}
W.R. Young, A.J. Roberts, and G.~Stuhne.
\newblock Reproductive pair correlations and the clustering of organisms.
\newblock \emph{Nature}, 412\penalty0 (6844):\penalty0 328--331, 2001.

\end{thebibliography}


\begin{thebibliography}{1}
\providecommand{\natexlab}[1]{#1}
\providecommand{\url}[1]{\texttt{#1}}
\expandafter\ifx\csname urlstyle\endcsname\relax
  \providecommand{\doi}[1]{doi: #1}\else
  \providecommand{\doi}{doi: \begingroup \urlstyle{rm}\Url}\fi

\bibitem[Gavagnin et~al.(2018)Gavagnin, Owen, and Yates]{gavagnin2018pcf}
E.~Gavagnin, J.P. Owen, and C.A. Yates.
\newblock Pair correlation functions for identifying spatial correlation in
  discrete domains.
\newblock \emph{Phys. Rev. E}, 97:\penalty0 062104, 2018.

\end{thebibliography}
\bibliographystyle{unsrtnat}

\end{document}


\title{Supplementary materials to accompany: ``Pair correlation functions for identifying spatial correlation in discrete domains''}
\author{\textbf{Enrico Gavagnin\footnote{These authors contributed equally to this work. Corresponding authors: EG: e.gavagnin@bath.ac.uk, JO: j.owen@bath.ac.uk}   \,, Jennifer P. Owen$^*$ and Christian A. Yates} \\ \small{\textit{Department of Mathematical Sciences},}\\ \small{\textit{University of Bath, Claverton Down, Bath, BA2 7AY,  UK}}}%
\date{}

\maketitle

\vfill
This document contains the supplementary materials which accompany the paper of \citet*{gavagnin2018pcf}. In Section \ref{sec:norm_unif} a full derivation of the normalisation factor for the Square Uniform PCF using non-periodic boundary conditions (BC) is provided. In Section \ref{sec:norm_tri} we report the expressions of the normalisations for the Triangle PCF and Hexagon PCF, using non-periodic BC. Section \ref{sec:norm_3D} contains the derivations of the normalisations for the Cube Taxicab PCF and the Cube Uniform PCF under non-periodic BC. This supplementary document concludes with a summary of the normalisation factors for all of our PCFs in Section \ref{sec:summary_norm}. We refer the reader to the main text for   the definitions of all of our PCFs.

\section{Derivation of the normalisation factor of the Square Uniform PCF under non-periodic BC}
\label{sec:norm_unif}
\begin{figure}[h!!!!!!]
	\begin{center}
		\subfigure[][]{\includegraphics[height=0.3 \columnwidth]{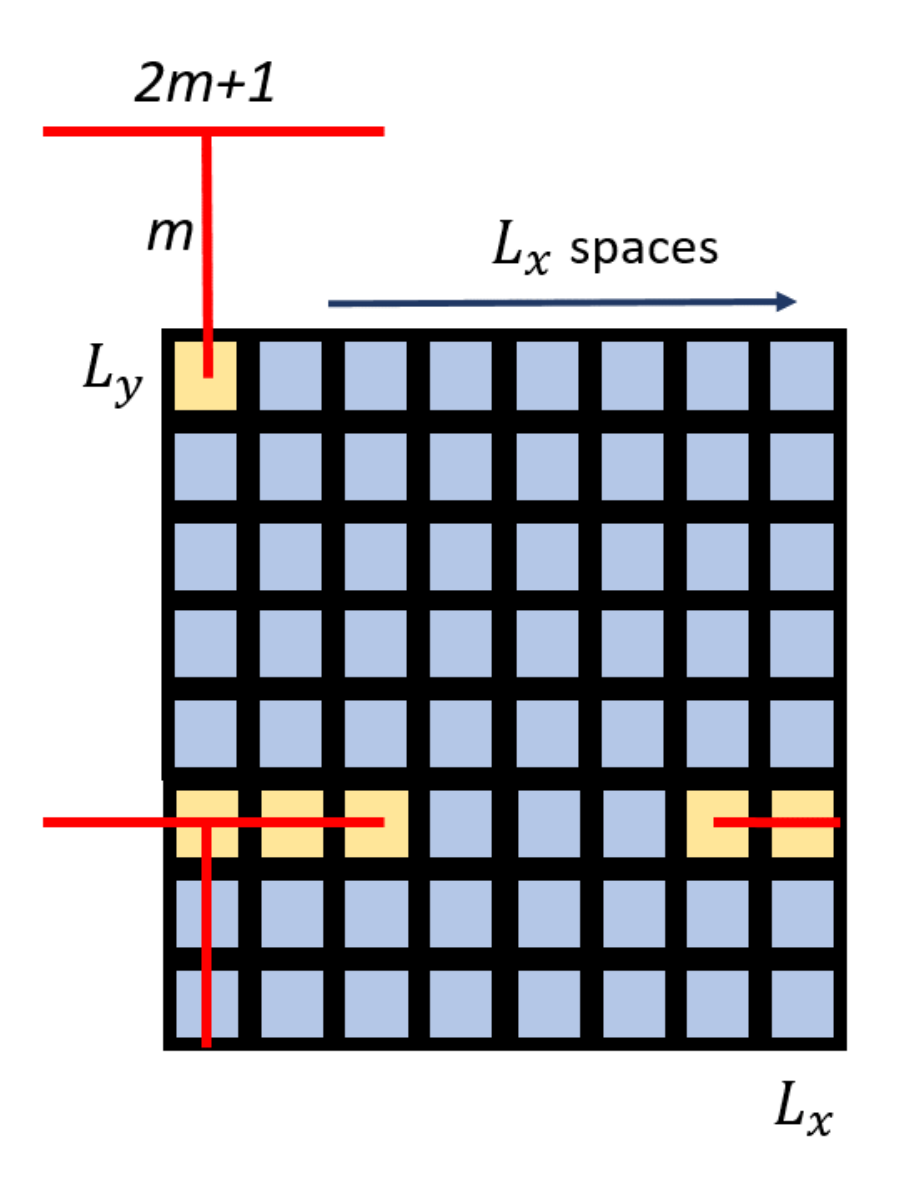}
			\label{fig:ex_unif_a} }
		\subfigure[][]{\includegraphics[height=0.3 \columnwidth]{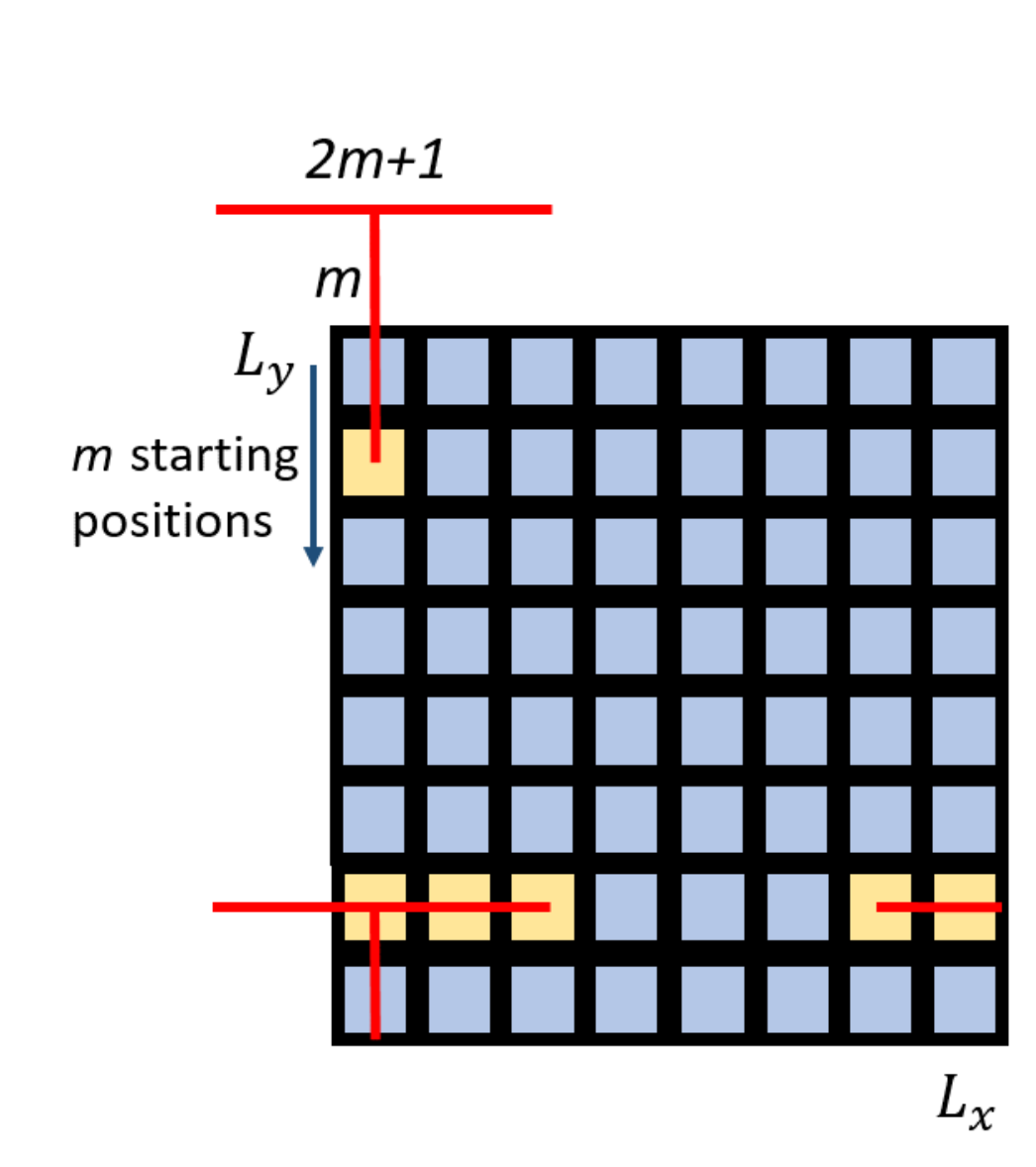}
			\label{fig:ex_unif_b} }\\
		\subfigure[][]{\includegraphics[ width=0.8 \columnwidth]{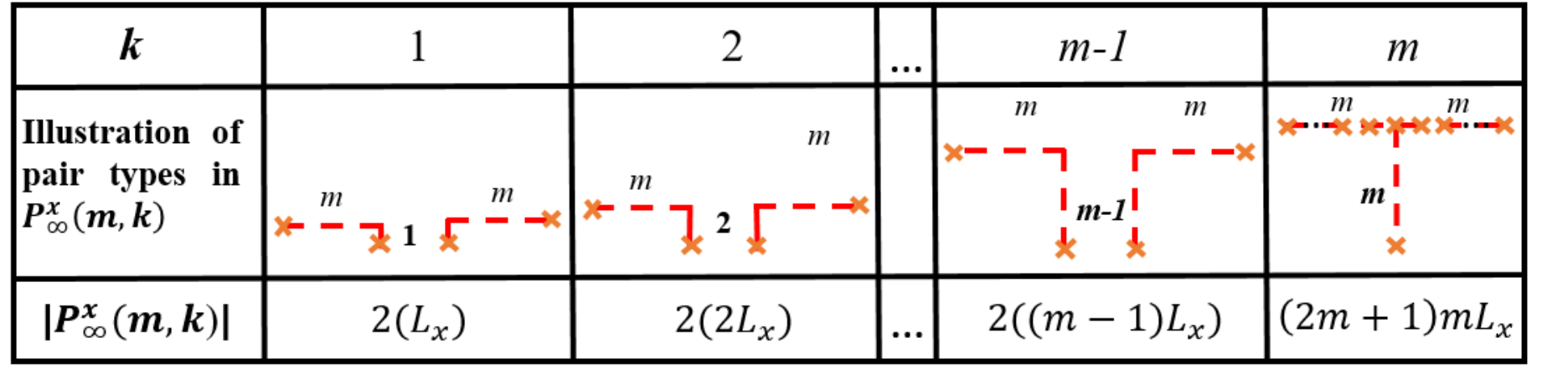}
			\label{fig:ex_binder_a} }
	\end{center}
	\caption{A visualisation of examples of the pairs of sites in $P^x_\infty(m)$. Panels (a) and (b) show different site pairs in $P^x_\infty(m,m)$. Panel (c) shows all the different types of pairs in $P^x_\infty(m,k)$ for $k=1, \dots , m$ with the corresponding value of $|P^x_\infty(m,k)|$.For example, for each of the $L_x$ columns, each site in the rows $\lbrace L_y-m+1, L_y-m+2, \dots , L_y \rbrace$ has $2m+1$ corresponding sites at distance $m$ separated by $m$ rows and reached by crossing the horizontal boundary. Therefore $|P^x_\infty(m,m)|=(2m+1)mL_x$. }
	\label{fig:ex_unif} 
\end{figure}

In this section we derive the normalisation factor for the Square Uniform PCF, defined in Section \ref{MAIN-sec:PCF} in the main text, under non-periodic BC. This corresponds to computing $\mathbb{E} [\bar{c}^{n}_\infty(m)]$ of equation\,\eqref{MAIN-eq:PCF_definition} of the main text. The approach that we use is similar to the one in Section \ref{MAIN-sec:norm_non_per}. In particular, by using the equation\,\eqref{MAIN-normalisation_factor_1} of the main text and the expression of $s^p_\infty(m)$ given by equation\,\eqref{MAIN-eq:def_s_unif} we only need to determine the expression of the reminder defined as
\begin{equation}
\label{eq:def_rem_unif}
	r_\infty(m)=s^p_\infty(m)-s^n_\infty(m) \; .
\end{equation}
We define the set of pairs of sites separated by distance $m \in \mathcal{D}^n_\infty$ that cross the $x$ boundary (horizontal axis) or $y$ boundary (vertical axis), respectively, as
\begin{subequations}
\begin{align}
P_\infty^x(m)=&\big\lbrace (\boldsymbol{a}, \boldsymbol{b}) \in S_\infty^p(m)\, \  \big| \ |y_a-y_b|>L_y-|y_a-y_b| \big\rbrace, 
\\ P_\infty^y(m)=&\big\lbrace (\boldsymbol{a}, \boldsymbol{b}) \in S_\infty^p(m)\, \  \big| \ |x_a-x_b|>L_x-|x_a-x_b| \big\rbrace, 
\end{align}
\end{subequations}
where $S_\infty^p(m)$ is defined as in equation\,\eqref{MAIN-eq:S_def_2} of the main text. Within these sets, $P^x_\infty$ and $P^y_\infty$, we select those separated by  $k \in \lbrace 1,... m \rbrace$ rows or columns, respectively. We define these subsets as:
\begin{subequations}
\begin{align}
P^x_\infty(m,k)=&\big\lbrace (\boldsymbol{a}, \boldsymbol{b}) \in \, P_\infty^x(m) \, \,\big| \quad L_y-|y_a-y_b|=k \big\rbrace,
\\P^y_\infty(m,k)=&\big\lbrace (\boldsymbol{a}, \boldsymbol{b}) \in \, P_\infty^y(m) \, \, \big| \quad L_x-|x_a-x_b|=k \big\rbrace . 
\end{align}
\end{subequations}
Fig.\,\ref{fig:ex_unif} (a) and (b) show an example of pairs of sites within $P_\infty^x(m,m)$.

By following the same steps as in Section \ref{MAIN-sec:norm_non_per}, we reduce to the expression for the remainder given by
\begin{align}
r_\infty(m)&=|P^{x}_\infty(m) \cup P^{y}_\infty(m)| \nonumber
\\&=\sum_{k=1}^m | P_\infty^x(m,k) | + \sum_{k=1}^m | P_\infty^y(m,k) | - |P_\infty^x(m) \cap P_\infty^y(m)| \label{eq:rem_unif}	\, .
\end{align}
To conclude the computation we derive an expression for the two sums in equation\,\eqref{eq:rem_unif} and the corresponding equation for the size of the intersection. 
By counting the contribution of each type of pair (see Fig.\,\ref{fig:ex_unif}\,(c) for a visualisation), one can write down the following expressions for the two sums in equation\,\eqref{eq:rem}: 
\begin{subequations}
\begin{align}
	\sum_{k=1}^m | P_\infty^x(m,k) | &= 2(L_x + 2L_x +  \dots L_x(m-1)) + m(2m+1)L_x \;,
	\\\sum_{k=1}^m | P_\infty^y(m,k) | &= 2(L_y + 2L_y +  \dots L_y(m-1)) + m(2m+1)L_y \; .
\end{align}
\end{subequations}
Hence
 \begin{align}
\sum_{k=1}^m | P_\infty^x(m,k) | + \sum_{k=1}^m | P_\infty^y(m,k) |&=2(L_x + 2L_x + , \dots, L_x(m-1)) + m(2m+1)L_x \nonumber \\&+2(L_y + 2L_y + , \dots, L_y(m-1)) + m(2m+1)L_y \nonumber
\\&=(L_x+L_y)\bigg(m(2m+1)+2\sum_{i=1}^{m-1}i\bigg) \nonumber
\\&=(L_x+L_y)\bigg(m(2m+1)+2\frac{(m-1)m}{2}\bigg) \nonumber
\\&=3(L_x+L_y)m^2 \, .
\label{eq:diff_unif}
\end{align}

\begin{figure}
		\includegraphics[width=0.35 \columnwidth]{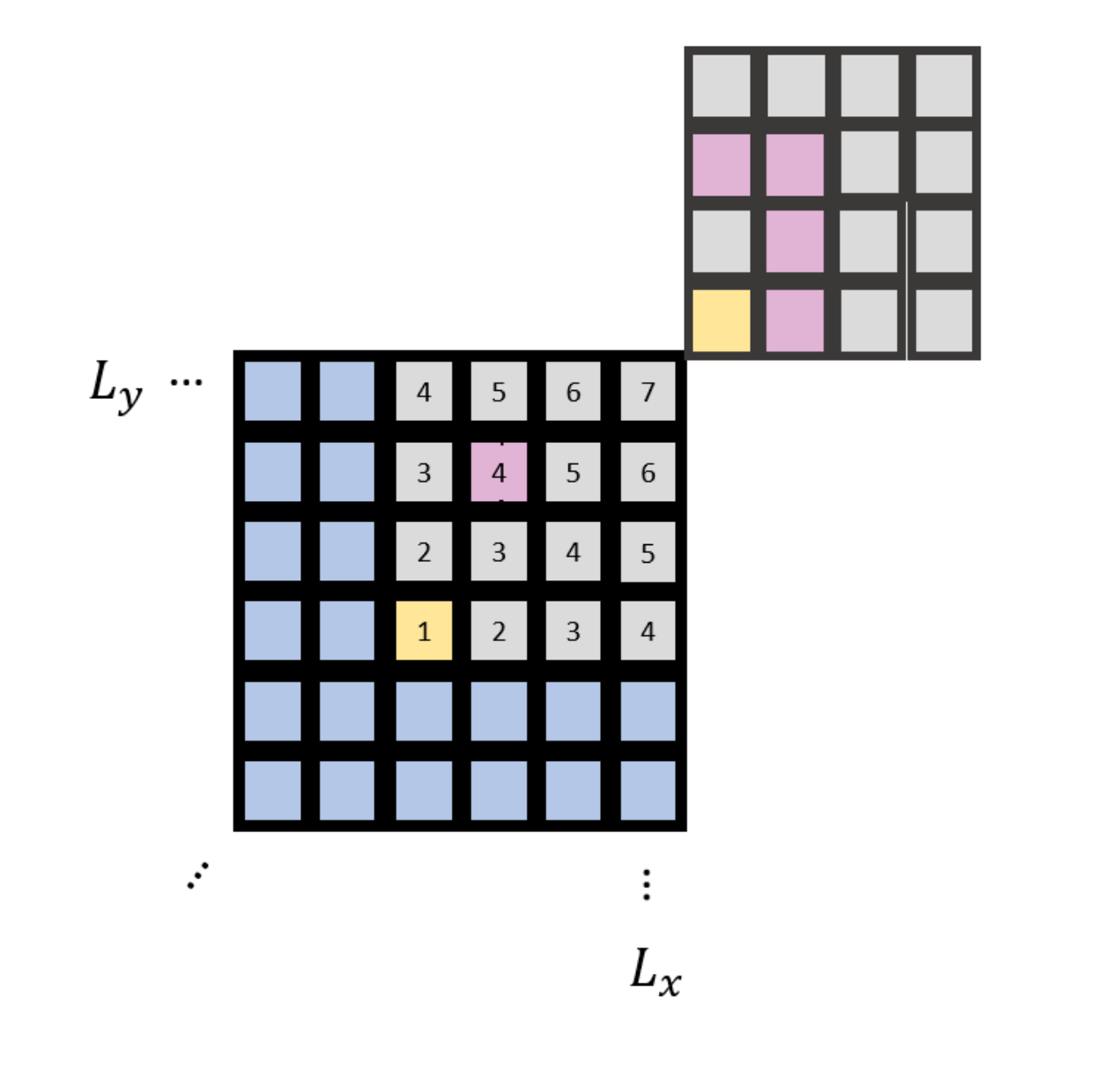}
		\caption{Examples of pairs in $P_\infty^x(4) \cap P_\infty^y(4)$ on a zoomed-in corner of a larger domain. Sites in pink and yellow are distance $m=4$ away from other sites in pink and yellow respectively. The number in each site correspond to the number of sites at distance $m=4$, reached by crossing the $x$ and $y$ boundaries.}
		\label{fig:cross_two_axis_unif}
\end{figure}
In order to calculate of the residue, $r_\infty (m)$, we need to compute the expression for the size of intersection, $|P_\infty^x(m) \cap P_\infty^y(m)| $, in equation\,\eqref{eq:rem_unif}. This counts the pairs of sites separated by distance $m$ by crossing both the $x$ and $y$ boundaries. Fig.\,\ref{fig:cross_two_axis_unif}  shows an example of pairs of sites within $P_\infty^x(4) \cap P_\infty^y(4)$. The number in each site corresponds to the number of sites at distance $m=4$, reached by crossing both the $x$ and $y$ boundaries. So $|P_\infty^x(4) \cap P_\infty^y(4)|$ is the sum of all the numbers in the coloured sites multiplied by two to account for the second corner region.

Extrapolating, for any value of $m$, the number of pairs of sites that cross the two boundaries is exactly
\begin{align}
|P_\infty^x(m) \cap P_\infty^y(m)|&=2\sum_{i=1}^m \sum_{j=i}^{m+i-1} j \nonumber
\\&=\sum_{i=1}^m [(m+i)(m+i-1)-i(i-1)] \nonumber
\\&=\sum_{i=1}^m [m^2+2m i -m] \nonumber
\\&=m^3-m^2+m^2(m+1)\nonumber
\\&=2m^3 \; .
\label{eq:unif_rem}
\end{align}
By substituting equations \eqref{eq:unif_rem} and \eqref{eq:diff_unif} into equation\,\eqref{eq:rem_unif} we obtain the expression for the remainder $r_\infty(m)$. By rearranging equation\,\eqref{eq:def_rem_unif} we determine $s_1^n(m)$, which we then substitute into equation\,\eqref{MAIN-normalisation_factor_1} of the main text to obtain the exact expression for the normalisation in the non-periodic case. This is given by
\begin{align}
\mathbb{E}\left[c^{U,n}_\infty(m)\right] =&\bigg(\frac{N}{L_xL_y}\bigg)\bigg(\frac{N-1}{L_xL_y-1}\bigg)\bigg(4mL_xL_y-3(L_x+L_y)m^2+2m^3\bigg). \label{eq:unif_non}
\end{align}

\section{Normalisation factor for the non-periodic Triangle PCF and Hexagon PCF}
\label{sec:norm_tri}
Here we provide the expressions of the normalisation factors for the Triangle PCF and the Hexagon PCF defined in Section \ref{MAIN-sec:tri_hex}, under non-periodic BC.

In both cases, the normalisation can be computed by using equation\,\eqref{MAIN-normalisation_factor_1} of the main text. Hence, we only need to provide the expression of the terms $s^n_{tri}(m)$, for the Triangle PCF, and $s^n_{hex}(m)$, for the Hexagon PCF. These correspond to the exact numbers of site pairs at distance $m$, assuming non-periodic BC, in the two types of tessellations (see Section \ref{MAIN-sec:tri_hex} of the main text for more details on the definition of the tessellated domain). We omit the derivations of the following expressions which can be obtained with similar steps to the square lattice cases. All the following expressions have been confirmed for a wide range of values of $L_x$, $L_y$ and $m$ by applying the General PCF (see Section \ref{MAIN-sec:ireg} of the main text).

\subsubsection*{Number of site pairs on a triangular tessellation with non-periodic BC}
\begin{subequations}
\begin{align}
s^n_{tri}(1)=&3 \frac{L_x L_y}{2}-\frac{L_x}{2}-L_y \, ,\\[15pt]
s^n_{tri}(2)=&3 L_x L_y-2L_x-4L_y+2 \, ,\\[15pt]
\begin{split}
s^n_{tri}(m)=&3 m \frac{L_x L_y}{2}-\frac{L_xm^2}{2} \qquad \qquad\text{(for $m\ge3$)}\\
&+Ly\left[2k_6(k_6-2k_3+1)+k_3(m-6)-m^2+m-2 \right]\\
& +\frac{1}{3} (m-1)(m^2-2m+6)\\&-\frac{1}{3}(k_7+1)(20k_7^2+37k_7+12)\\ 
&-(m-7-4k_7)(k_7+1)(m+k_7-2)\,,
\end{split}
\end{align}
\end{subequations}
where $k_j=\left\lfloor \frac{m-j}{4} \right\rfloor$.

\subsubsection*{Number of site pairs on a hexagonal tessellation with non-periodic BC}
\begin{align}
	s^{n}_{hex}(m)=3mL_x L_y -\frac{1}{4} \left( 7 m^2+k\right)L_x-2m^2 L_y+\frac{11m^3}{12}-\frac{(2-3k)m}{12} \, ,
\end{align}
where $k=m \pmod{2}$.


\section{Derivation of the normalisation factor of the Cube Taxicab PCF and Cube Uniform PCF under non-periodic BC}
\label{sec:norm_3D}
Consider a system of agents on a three-dimensional square lattice of size $L_x \times L_y \times L_z$, with lattice step $\Delta$ and with the exclusion property that, at any given time, each lattice site can be occupied by at most one agent.  
If $N$ agents occupy the domain, then the occupancy of the lattice can be represented by a matrix $M$: 
\begin{align}
M_{xyz}=
\begin{cases}
&  0 \text{  if $(x,y,z)$ is vacant,}\\
&1 \text{  if $(x,y,z)$ is occupied.}
\end{cases}
\end{align}
where
\begin{align}N=\sum_{x=1}^{L_x}\sum_{y=1}^{L_y}\sum_{z=1}^{L_z}M_{xyz} \leq L_xL_yL_z .\end{align}
Let $\psi$ be the set of all agent pairs in the lattice i.e.
\begin{align}
\begin{split}
\psi=&\{(\boldsymbol{a},\boldsymbol{b})\in \mathbb{L}\times\mathbb{L}\hspace{0.6em} |\hspace{0.6em} \boldsymbol{a}=(x_a,y_a,z_a),\   \boldsymbol{b}=(x_b,y_b,z_b),\  \boldsymbol{a}\neq \boldsymbol{b},\  M_{x_a,y_a,z_a}=M_{x_b,y_b,z_b}=1 \} 
\end{split} ,
\end{align}
where $\mathbb{L}=\{1,\dots, L_x\}\times \{1, \dots , L_y\} \times \{1, \dots , L_z\}$ 
is the set of all sites in the lattice. 

In this Section we derive the normalisations for the Cube Taxicab PCF and the Cube Uniform PCF under non-periodic BC.
Using the same notation as in Section \ref{MAIN-sec:previous} of the main text, we define the subsets of agent pairs separated by distance $m$ under non-periodic BC as
\begin{subequations}
\label{eq:def_bigS_nonper}
	\begin{align}
	C^{n}_1(m)=&\{(\boldsymbol{a},\boldsymbol{b}) \in \psi \ \big|\  \norm{\boldsymbol{a}-\boldsymbol{b}}_{1}=m\}, \quad  \, m \in \mathcal{D}^n_1,
	\\C_{\infty}^{n}(m)=&\{(\boldsymbol{a},\boldsymbol{b}) \in \psi \  \big|\  \norm{\boldsymbol{a}-\boldsymbol{b}}_\infty=m \}, \quad m \in \mathcal{D}^n_\infty,
	\end{align}
\end{subequations}
for the taxicab and uniform metric respectively, where $\mathcal{D}^n_1 = \mathcal{D}^n_\infty= \big\lbrace 1,2,...\max \lbrace L_x,L_y,L_z \rbrace -1 \big\rbrace$. Using the definitions of the uniform and taxicab metrics, we can express these sets as: 
\begin{subequations}
	\begin{align}
	C^{n}_1(m)=&\{(\boldsymbol{a},\boldsymbol{b}) \in \psi\  \big|\  |x_a-x_b|+|y_a-y_b|+|z_a-z_b|=m\} ,
	\\C_{\infty}^{n}(m)=&\{(\boldsymbol{a},\boldsymbol{b}) \in \psi\  \big |\  \max \lbrace |x_a-x_b|\text{, }|y_a-y_b|\text{, } |z_a-z_b| \rbrace=m  \} .
	\end{align}
\end{subequations}
Similarly we define the subsets of agent pairs separated by distance $m$ under periodic BC as

\begin{subequations}
\label{eq:def_bigS_per}	
	\begin{align}	
	&\begin{aligned}
	C^{p}_1(m)= \Big\lbrace(\boldsymbol{a},\boldsymbol{b})\in  \psi\mid \min \lbrace & |x_a-x_b|, L_x-|x_a-x_b| \rbrace \\&+ \min \lbrace |y_a-y_b|, L_y-|y_a-y_b|\rbrace \\&+ \min \lbrace |z_a-z_b|, L_y-|y_a-y_b|\rbrace=m\Big\rbrace , \quad  \, m \in \mathcal{D}^p_1,
	\end{aligned}\\
	&\begin{aligned}
	C^{p}_\infty(m)=\Big\lbrace (\boldsymbol{a},\boldsymbol{b})\in \psi\mid \max \big\lbrace & \min \lbrace |x_a-x_b|, L_x-|x_a-x_b| \rbrace, \\& \min\lbrace |y_a-y_b|, L_y-|y_a-y_b|\rbrace , \\& \min\lbrace |z_a-z_b|, L_z-|z_a-z_b|\rbrace \big\rbrace =m \Big\rbrace , \quad  \, m \in \mathcal{D}^p_\infty,
	\end{aligned}
	\end{align}
\end{subequations}
where $\mathcal{D}^p_1 = \mathcal{D}^p_\infty= \big\lbrace 1,2,...\min \big\lbrace \big\lfloor\frac{L_x}{2}\big\rfloor,\big\lfloor\frac{L_y}{2}\big\rfloor, \big\lfloor\frac{L_z}{2}\big\rfloor \big\rbrace \big\rbrace$. 
The schematics in Fig.  \ref{MAIN-fig:normsvis3D} of the main text represent examples of sites separated by distance $m=1$ using the taxicab (a) and uniform (b) metrics in three dimensions.
The formulae for the number of pairs of sites separated by distance $m$ under periodic BC are given in Section \ref{MAIN-sec:cube} of the main text. Here we compute the numbers of pairs of sites separated by  distance $m$ under non-periodic BC, \textit{i.e} $s^{n}_1(m)$ and $s^{n}_\infty (m)$, which we can then substitute into the equation\,\eqref{MAIN-normalisation_factor_1} in the main text to obtain the normalisation factors of the corresponding PCFs.

Similarly to the two-dimensional cases, we focus on computing an expression the remainder given by
	\begin{equation}
	\label{eq:def_rem}
	r_d(m)=s^{p}_d(m)-s^{n}_d(m),
	\end{equation} 	
where $d=1,\infty$. The expression of $s^n_d(m)$ will then follow by rearrangement. 

\begin{figure}[h!!]
	\begin{center}
		\subfigure[][]{\includegraphics[height=0.4 \columnwidth]{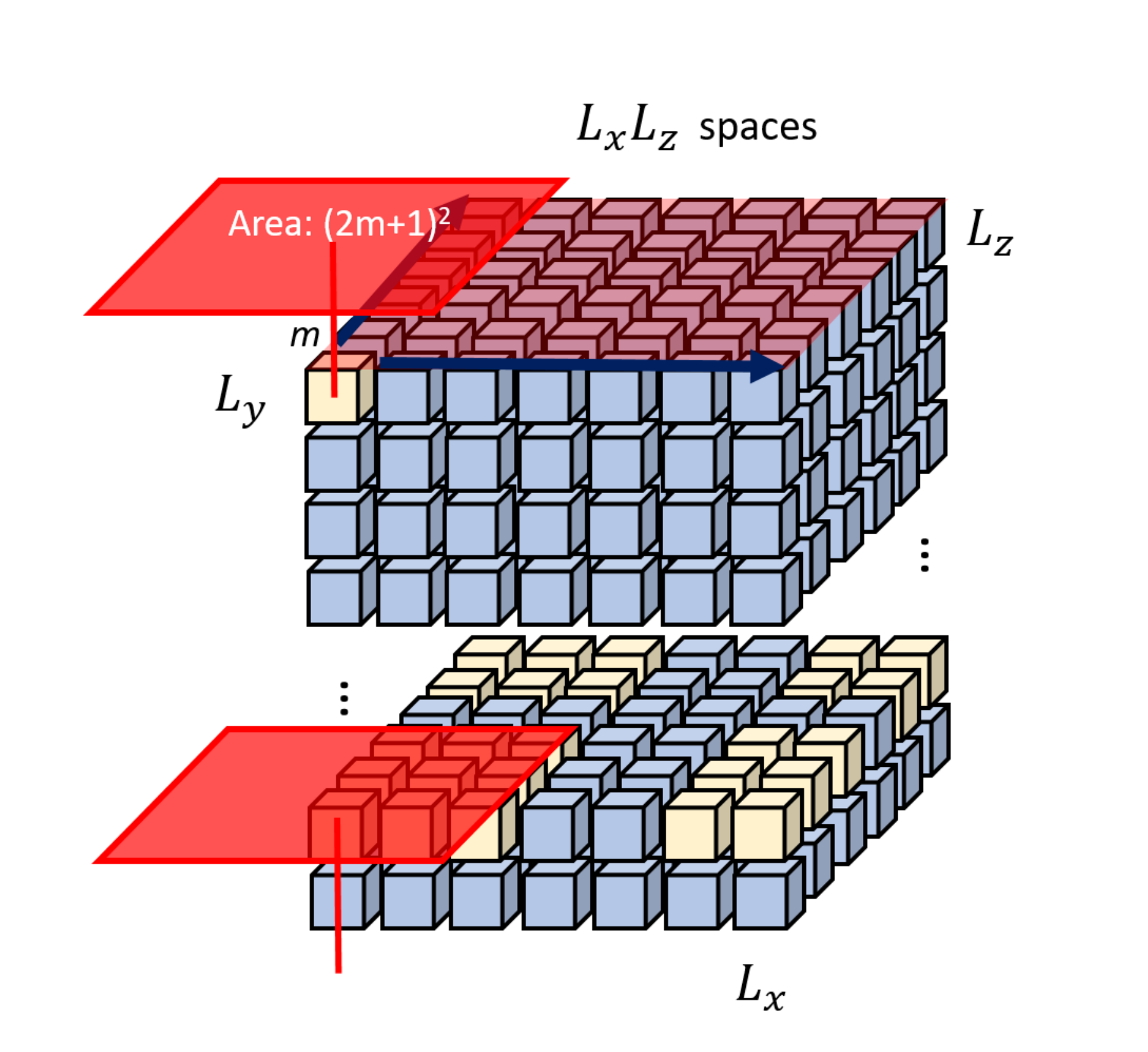}
			\label{fig:ex_unif_a} }
		\subfigure[][]{\includegraphics[height=0.4 \columnwidth]{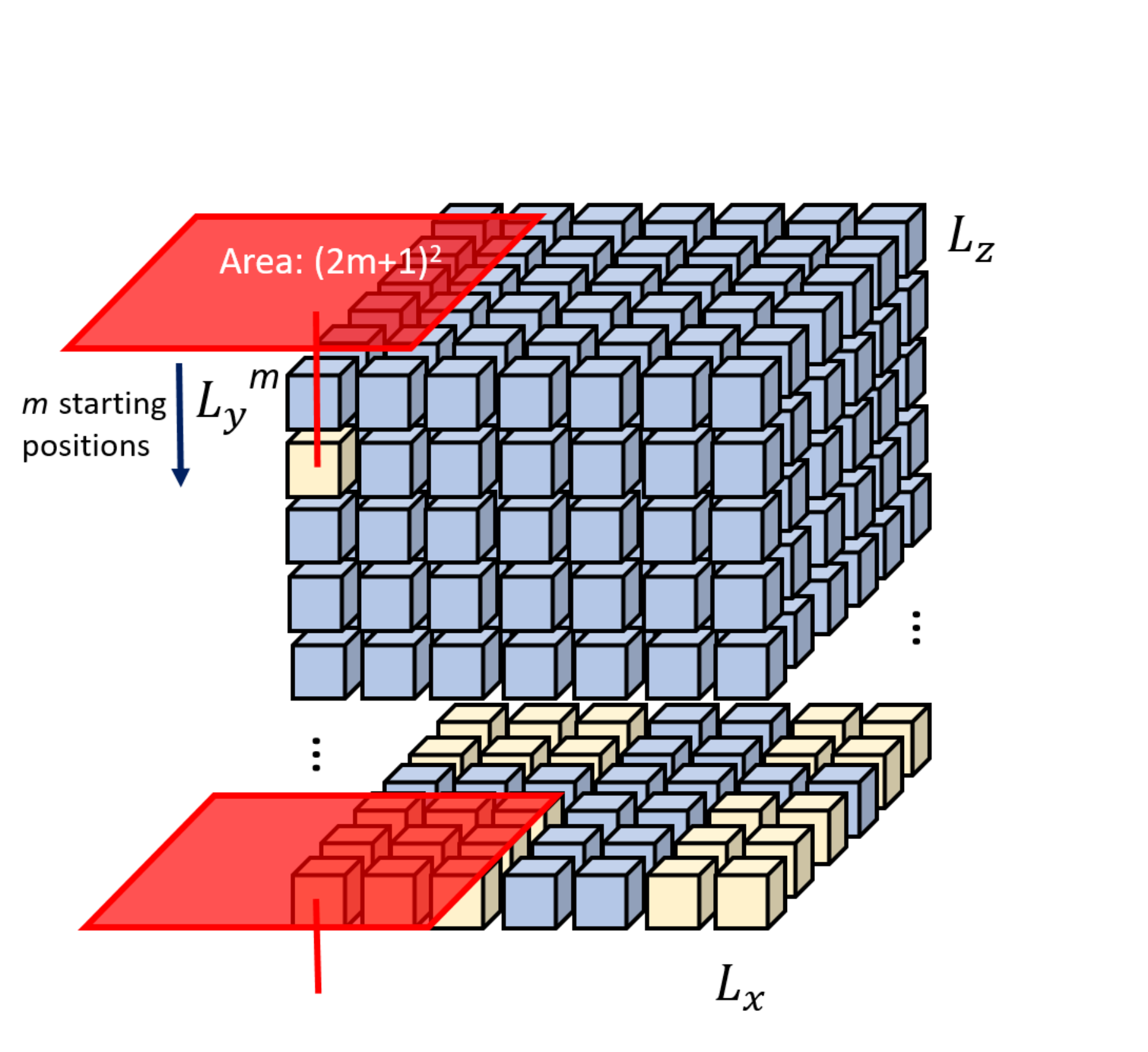}
			\label{fig:ex_unif_b} }
		\subfigure[][]{\includegraphics[height=0.2 \columnwidth]{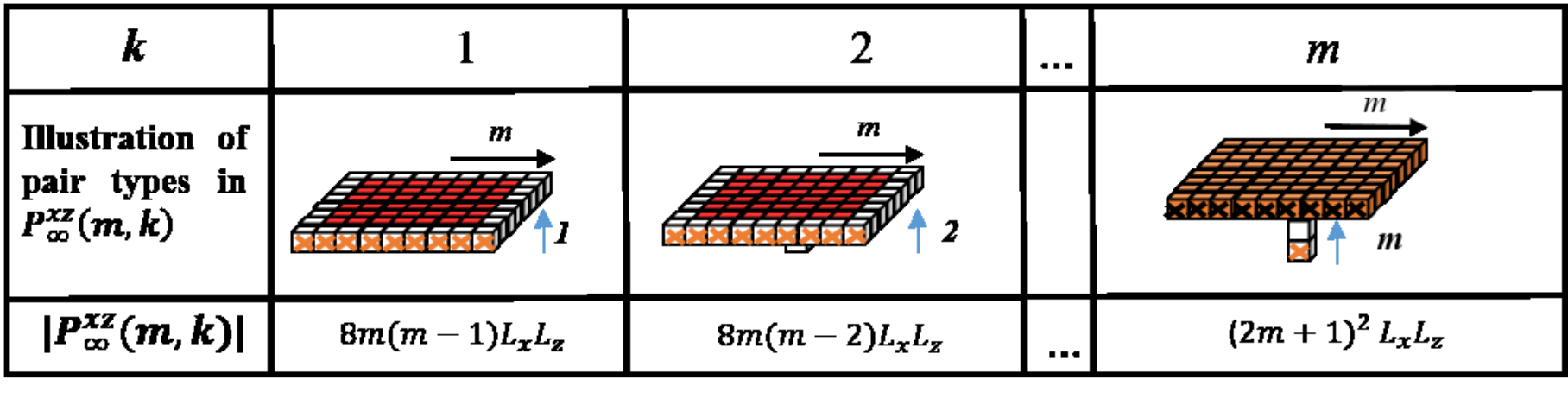}
		} \label{eq:example}
	\end{center}
	\caption{A visualisation of the pairs of sites in $P^{xz}_\infty(m)$. Panels (a) and (b) show two different site pairs in $P^{xz}_\infty(m,m)$. Panel (c) shows all the different types of pairs in $P^{xz}_\infty(m,k)$ for $k=1, \dots , m$ with the corresponding value of $|P^{xz}_\infty(m,k)|$. For example, each site in the top $m$ horizontal planes, \textit{i.e.} $\lbrace y=L_y-m+1, y=L_y-m+2, \dots , y=L_y \rbrace$, has $(2m+1)^2$ sites at distance $m$ separated by $m$ horizontal planes which are reached by crossing the $\xz$ boundary. Therefore $|P^{xz}_\infty(m,m)|=(2m+1)^2L_xL_z$.}
	\label{fig:3D_unif} 
\end{figure}

Let us define the set of pairs of sites separated by distance $m \in \mathcal{D}^n_1$ that cross the $\yz$ plane boundary (\textit{i.e.} $x=0$), $\xz$  plane boundary (\textit{i.e.} $y=0$) and $\xy$ plane boundary (\textit{i.e.} $z=0$), respectively as
\begin{subequations}
	\begin{align}
P_1^{yz}(m)=&\big\lbrace (\boldsymbol{a}, \boldsymbol{b}) \in S_1^p(m)\, \  \big| \ |x_a-x_b|>L_x-|x_a-x_b| \big\rbrace	, 
	\\ P_1^{xz}(m)=&\big\lbrace (\boldsymbol{a}, \boldsymbol{b}) \in S_1^p(m)\, \  \big| \ |y_a-y_b|>L_y-|y_a-y_b| \big\rbrace,
		\\ P_1^{xy}(m)=&\big\lbrace (\boldsymbol{a}, \boldsymbol{b}) \in S_1^p(m)\, \  \big| \ |z_a-z_b|>L_z-|z_a-z_b| \big\rbrace. 
	\end{align}
\end{subequations}
For $P_1^{yz}(m)$, $P_1^{xz}(m)$ and $P_1^{xy}(m)$ let us consider those pairs of sites separated by $k \in \lbrace 1,... m \rbrace$ sites in the corresponding orthogonal direction, \textit{i.e.} $x-$, $y-$ and $z-$direction, respectively. We define these subsets as:
\begin{subequations}
	\begin{align}
	P^{yz}_1(m,k)=&\big\lbrace (\boldsymbol{a}, \boldsymbol{b}) \in \, P_1^{yz}(m) \big| \quad L_x-|x_a-x_b|=k \big\rbrace,
	\\P^{xz}_1(m,k)=&\big\lbrace (\boldsymbol{a}, \boldsymbol{b}) \in \, P_1^{xz}(m) \big| \quad L_y-|y_a-y_b|=k \big\rbrace
		\\P^{xy}_1(m,k)=&\big\lbrace (\boldsymbol{a}, \boldsymbol{b}) \in \, P_1^{xy}(m) \big| \quad L_z-|z_a-z_b|=k \big\rbrace  . 
	\end{align}
\end{subequations}
Notice that we have $P_d^{yz}(m)=\bigcup_{k =1}^m P_d^{yz}(m,k)$, $P_d^{xz}(m)=\bigcup_{k =1}^m P_d^{xz}(m,k)$ and $P_d^{xy}(m)=\bigcup_{k =1}^m P_d^{xy}(m,k)$.
Figs.\,\ref{eq:example}\,(a)-(b) provide a visualisation of pairs of sites within $P_1^{xy}(m,m)$. In Fig.\,\ref{eq:example}\,(c) we show the possible types of pairs of sites in $P_1^x(m,k)$, for $k=1,\dots,m$.
By definitions \eqref{eq:def_bigS_nonper} and \eqref{eq:def_bigS_per} we have that
\begin{align}
S_d^p(m) \, \backslash \, S_d^n(m)=P^{yz}_d(m) \cup P^{xz}_d(m) \cup P^{xy}_d(m)\,,
\end{align}
where $d=1,\infty$. Hence, we obtain
\begin{align}
r_d(m)=&|P^{yz}_d(m) \cup P^{xz}_d(m) \cup P^{xy}_d(m)| \nonumber
\\\begin{split}=&\sum_{k=1}^m | P_d^{yz}(m,k) | + \sum_{k=1}^m | P_d^{xz}(m,k) | +\sum_{k=1}^m | P_d^{xy}(m,k) | \\&-|P_d^{yz}(m) \cap P_d^{xy}(m)|-|P_d^{xz}(m) \cap P_d^{yz}(m)|-|P_d^{xz}(m) \cap P_d^{xy}(m)|\\&+ |P_d^{yz}(m) \cap P_d^{xz}(m) \cap P_d^{xy}(m)| \label{eq:rem}	\, ,
\end{split}	
\end{align}
for $d=1,\infty$.

\begin{figure}[h!!]
\label{fig:3D_taxicab} 
	\begin{center}
		\subfigure[][]{\includegraphics[height=0.4 \columnwidth]{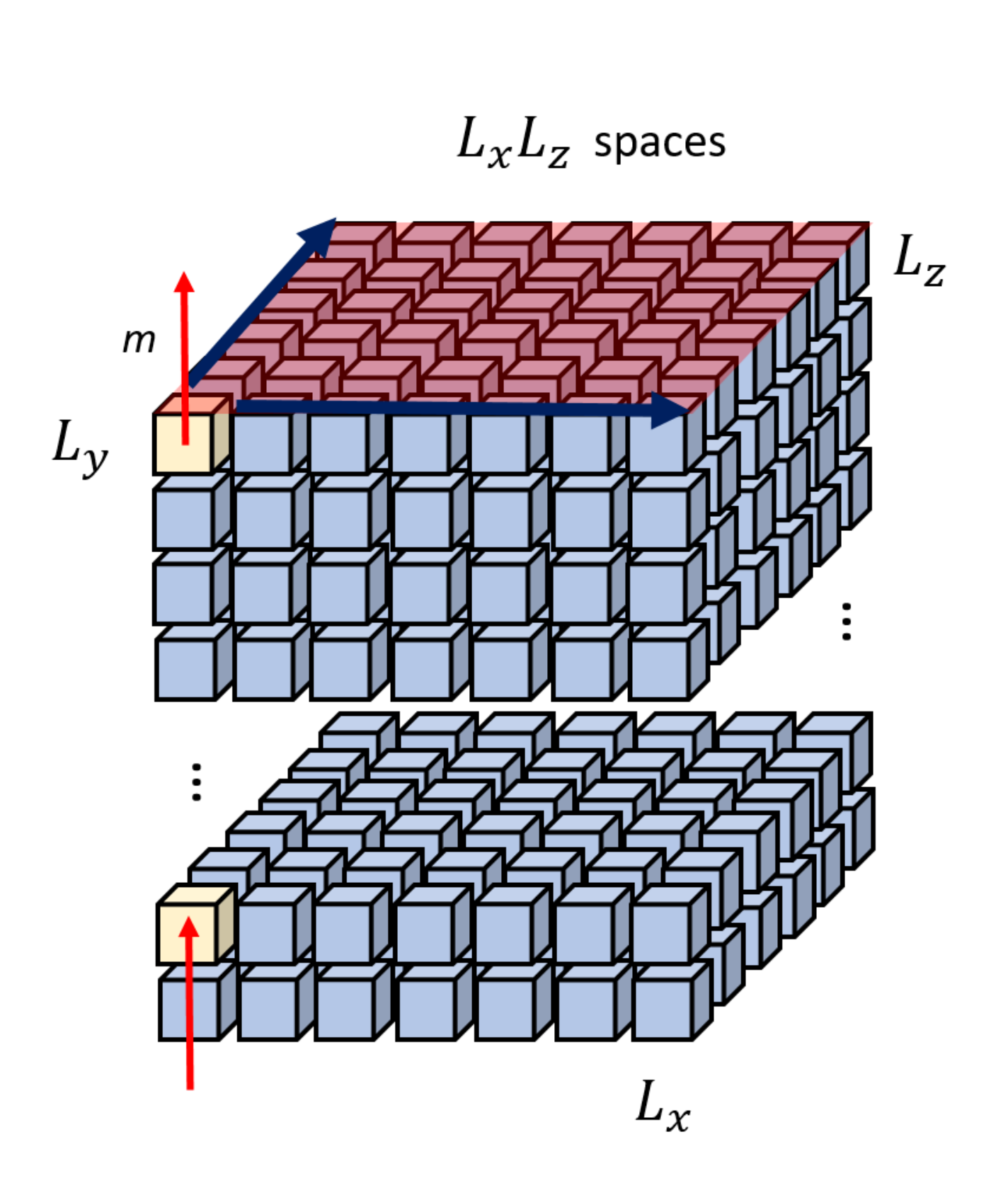}
			\label{fig:ex_unif_a} }
		\subfigure[][]{\includegraphics[height=0.4 \columnwidth]{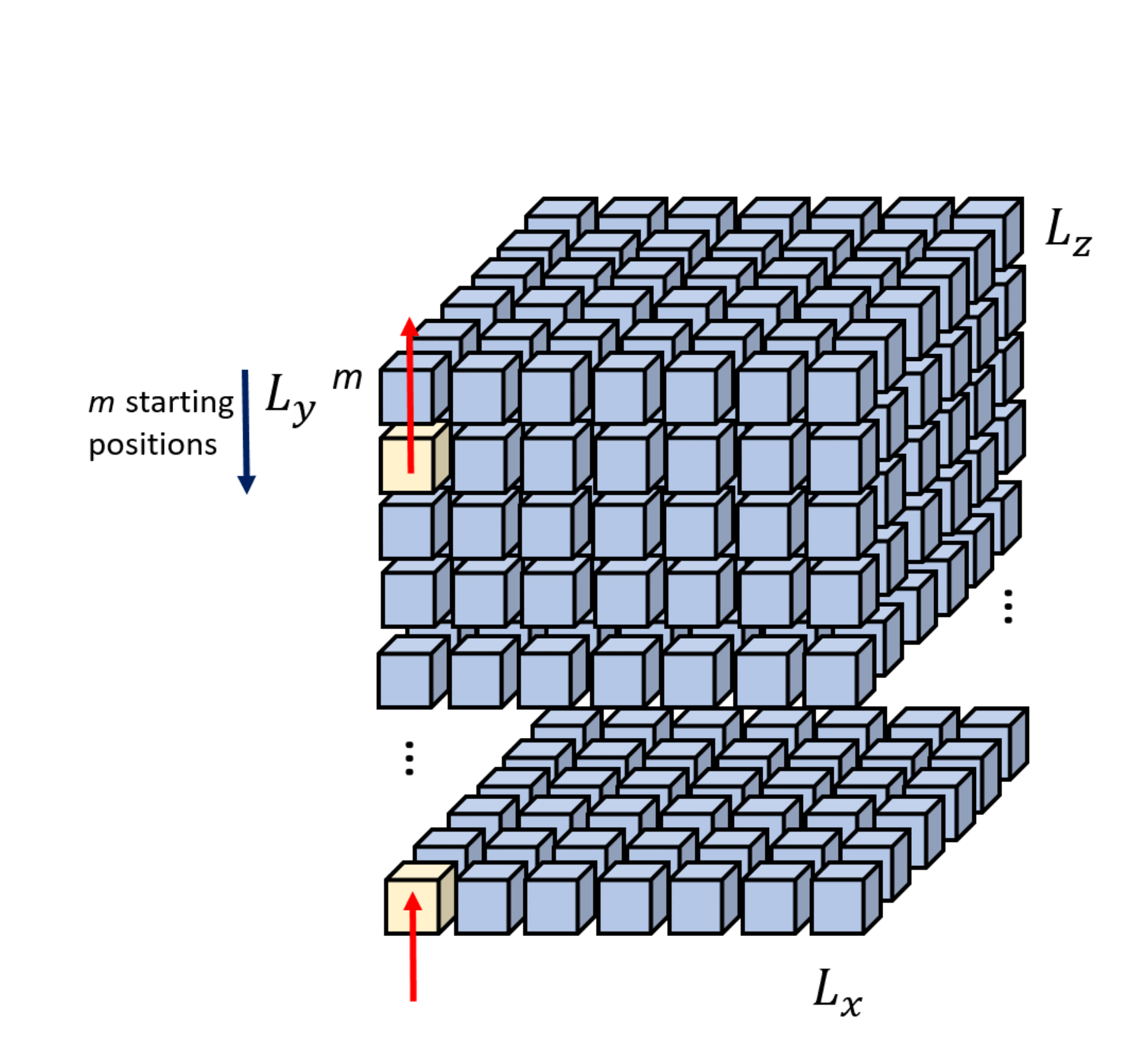}
			\label{fig:ex_unif_b} }
				\subfigure[][]{\includegraphics[height=0.2 \columnwidth]{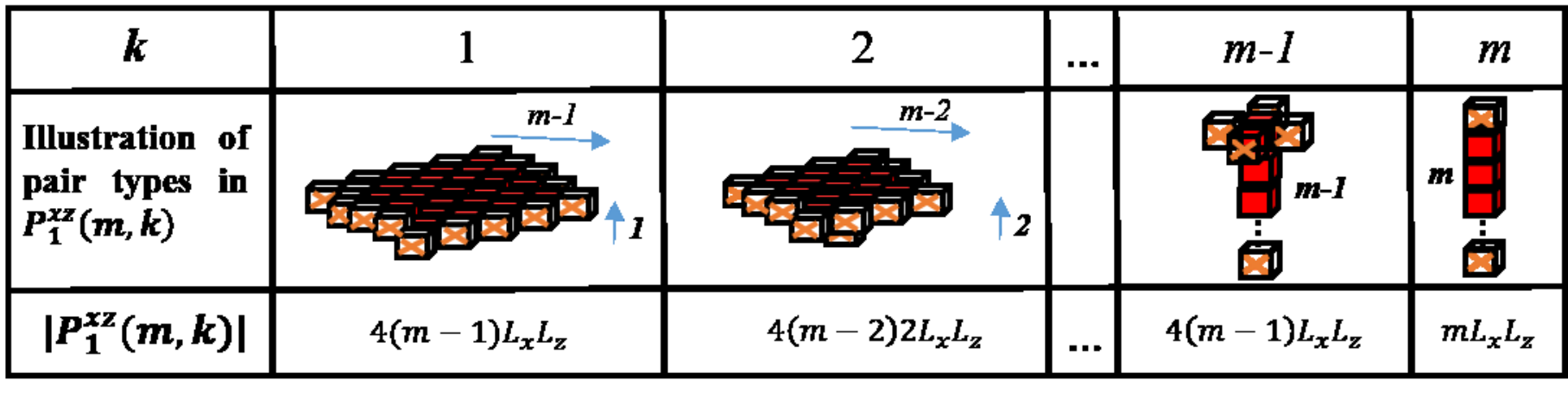} \label{fig:ex_taxicab_c} }
	\end{center}
	\caption{A visualisation of the pairs of sites in $P^{xz}_1(m)$. Panels (a) and (b) show two different site pairs in $P^{xz}_1(m,m)$. Panel (c) shows all the different types of pairs in $P^{xz}_1(m,k)$ for $k=1, \dots , m$ with the corresponding value of $|P^{xz}_1(m,k)|$. For example, each site in the top $m$ horizontal planes, \textit{i.e.} $\lbrace y=L_y-m+1, y=L_y-m+2, \dots , y=L_y \rbrace$, has a single corresponding site at distance $m$ separated by $m$ horizontal planes and 0 vertical planes, which is reached by crossing the $\xz$ boundary. Therefore $|P^{xz}_1(m,m)|=mL_xL_y$. The axis orientation is chosen to be consistent with the two-dimensional case. }
\end{figure}

\begin{figure}[h!!]
\begin{center}
{		\subfigure[][]{\includegraphics[height=0.28 \columnwidth]{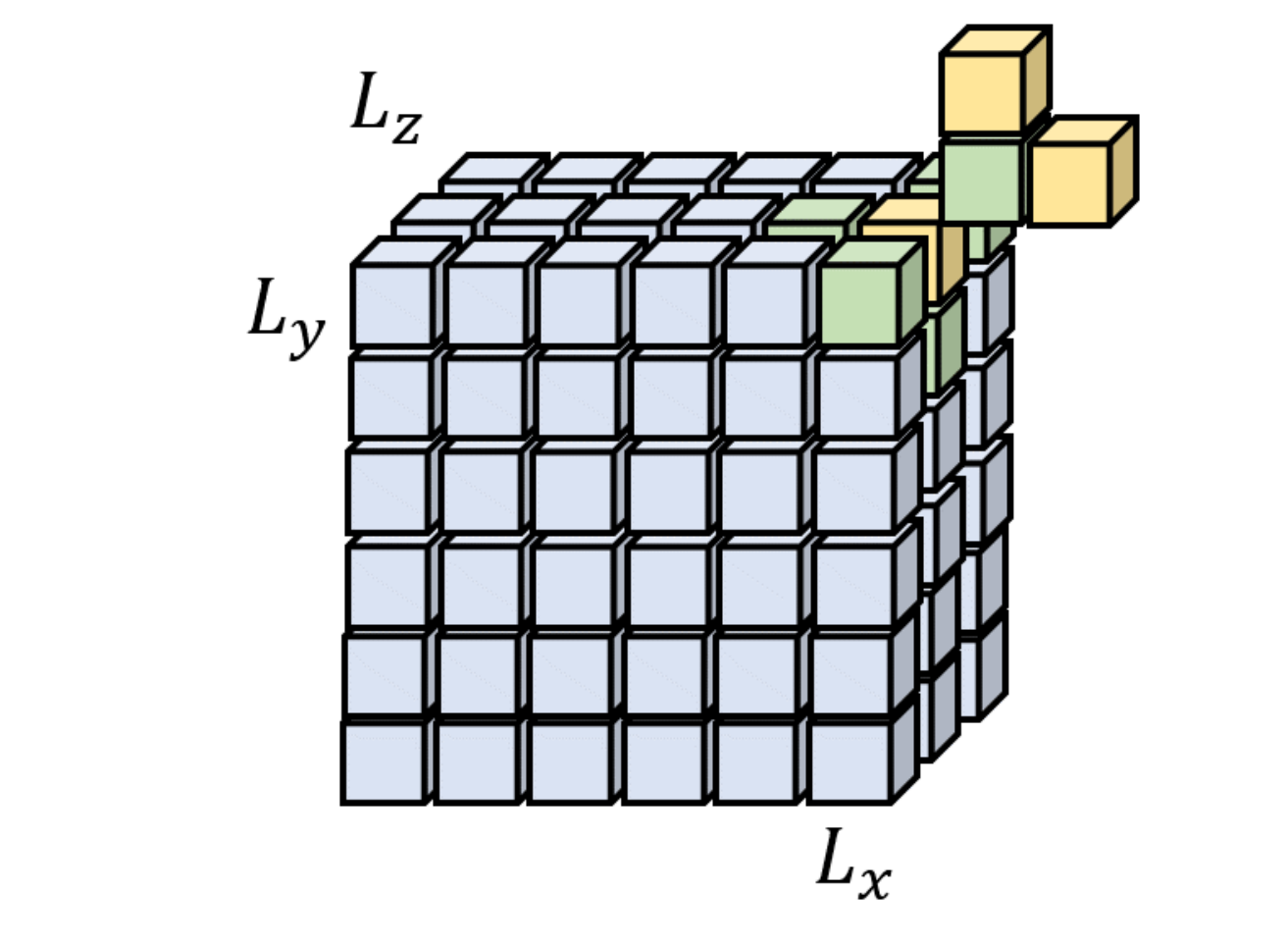}}
		\subfigure[][]{\includegraphics[height=0.32 \columnwidth]{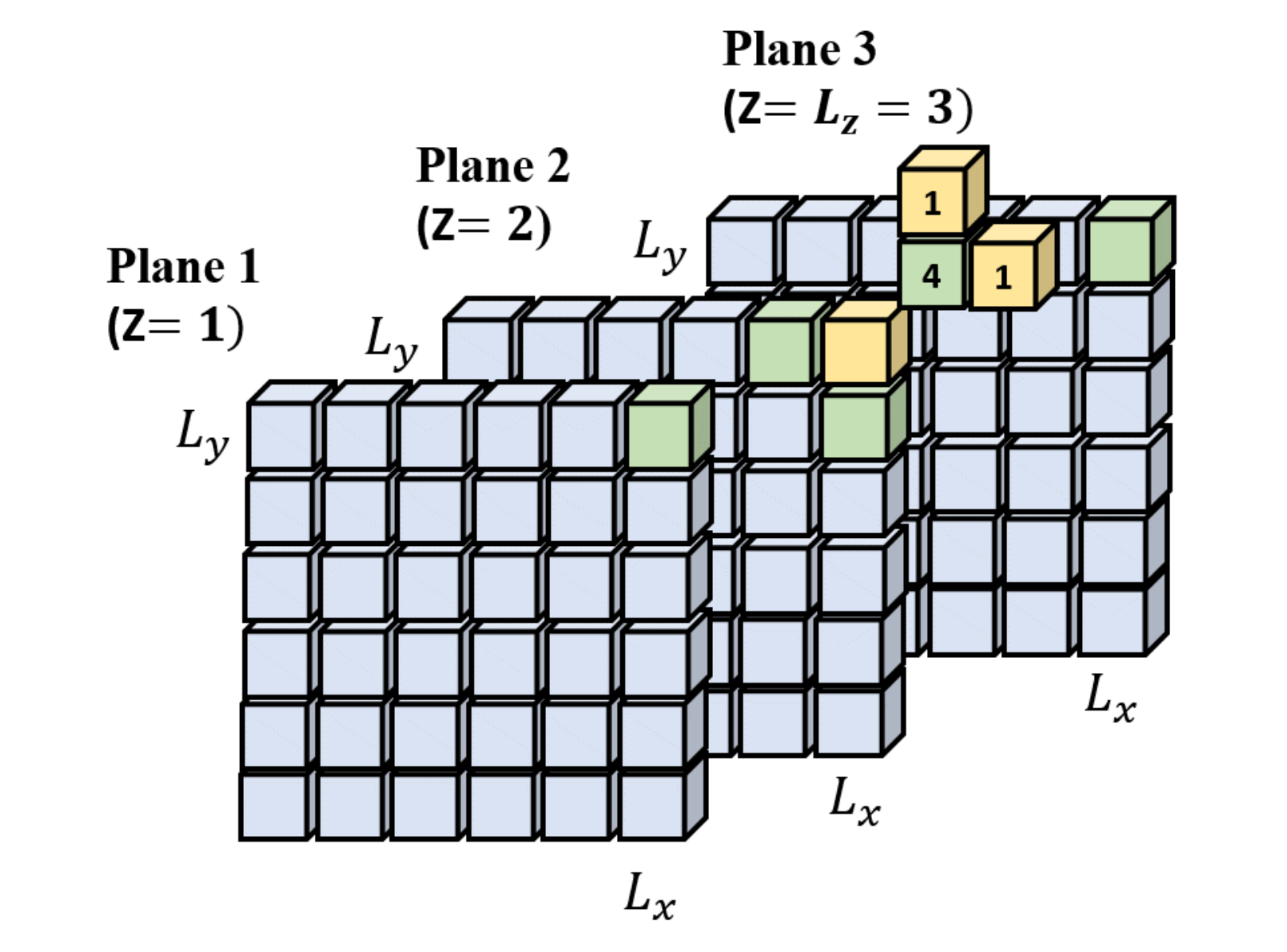}
			 }
\caption{Examples of pairs of sites separated by distance $m=3$ that cross both the $\xz$ and $\yz$ plane boundaries, \textit{i.e.} pairs of sites in $P_1^{xz}(3) \cap P_1^{yz}(3)$. Green and yellow sites inside of the cube domain are distance $m=3$ from green and yellow sites outside of the domain respectively. Panel (a) illustrates an example three-dimensional cube that can be split into $L_z$ vertical planes, panel (b) illustrates how the $\xy$ plane $z=2$ contributes sites that are distance $m$ from other sites by crossing both the $\xz$ and $\yz$ boundaries. }
\label{fig:taxcab_edge}}
\end{center}
\end{figure}

By counting the contribution of each type of pair (see Fig.\,\ref{fig:3D_taxicab} (c) for a visualisation), one can write down the following expressions for the first three sums in equation\,\eqref{eq:rem}: 
\begin{subequations}
	\begin{align}
	\sum_{k=1}^m | P_1^{ij}(m,k) | &= L_iL_j\bigg(4\left(m-1+2(m-2)+3(m-3)+...+m-1\right) +m\bigg)\;,
	\\ 	\sum_{k=1}^m | P_\infty^{ij}(m,k) | &= L_iL_j\bigg(8m\left(m-1+(m-2)+(m-3)+...+1\right) +(2m+1)^2\bigg)\;,
	\end{align}
\end{subequations}
where $ij={yz,xz,xy}$. Hence for the case of taxicab metric we have
\begin{align}
\sum_{k=1}^m | P_1^{yz}(m,k) | + \sum_{k=1}^m | P_1^{xz}(m,k) |+& \sum_{k=1}^m | P_1^{xy}(m,k) |\nonumber \\=&\bigg(4 \sum_{i=1}^{m-1} i(m-i)\bigg)(L_xL_y+L_yL_z+L_zL_x) \nonumber
\\=&\frac{2}{3}m(m-1)(m+1)(L_xL_y+L_yL_z+L_zL_x)\nonumber \\=&\frac{1}{3}(2m^3+m)(L_xL_y+L_yL_z+L_zL_x) \, ;
\label{eq:diff_tc}
\end{align}
and similarly for the uniform metric:
\begin{align}
\sum_{k=1}^m | P_\infty^{yz}(m,k) | + \sum_{k=1}^m | &P_\infty^{xz}(m,k) |+ \sum_{k=1}^m | P_\infty^{xy}(m,k) | \nonumber \\=&\bigg((2m+1)^2+8m \sum_{i=1}^{m-1} i\bigg)(L_xL_y+L_yL_z+L_zL_x) \nonumber
\\=&m(8m^2+1)(L_xL_y+L_yL_z+L_zL_x) \,.
\label{eq:diff_tc}
\end{align}
We now focus on deriving an expression for the unknown expressions in \eqref{eq:rem}, that is the size of intersections, $P_d^{xy}(m) \cap P_d^{yz}(m)$, $P_d^{xy} \cap P_d^{xz}(m) $ and $P_d^{xy} \cap P_d^{xz}(m)$. These sets consist of pairs of sites separated by distance $m$ that cross \textit{at least} two of the three $\yz,\xz,\xy$ boundaries. Examples of pairs of sites in $P_1^{xz}(m) \cap P_1^{yz}(m)$ are illustrated in Fig.\,\ref{fig:taxcab_edge}\,(a) and Fig.\,\ref{fig:taxcab_edge}\,(b). We show that $P_1^{xz}(m) \cap P_1^{yz}(m)$ consists of all of the pairs of sites located in the top left and right edges (right demonstrated only) which cross the $\xz$ and $\yz$ boundaries of each of the $L_z$ vertical planes making up the cube (see Fig.\,\ref{fig:taxcab_edge}\,(b)). The sites that comprise pairs in $P_1^{xz}(m) \cap P_1^{yz}(m)$ can both be in the same $\xy$ plane or they can be in neighbouring (up to $m-2$) planes away. The number of pairs in these sets can then be calculated by addition, enumerating possible pairs as in Fig.\,\ref{MAIN-fig:cross_two_axis} in the main text for taxicab case and Fig.\,\ref{fig:cross_two_axis_unif} for the uniform case respectively to give the following. For any $s,t,u \in \lbrace x,y,z \rbrace$, such that $(s,t) \neq (t,u)$ 
\begin{align}
|P_1^{st}(m) \cap P_1^{tu}(m)|=&2L_t\sum_{i=1}^m i^2 (m-i) \nonumber
\\=&\frac{1}{6}(m-1)m^2(m+1)L_t,
\\|P_\infty^{st}(m) \cap P_\infty^{tu}(m)|=&2L_t\bigg((2m+1)m^3 + 2\sum_{i=1}^{m-1}\sum_{j=1}^{m-1}ij\bigg)\ \nonumber
\\=&m^2(5m^2+1)L_t.
\end{align}

\begin{figure}[h!!!!!!!!]
	\begin{center} {
		\includegraphics[height=0.35 \columnwidth]{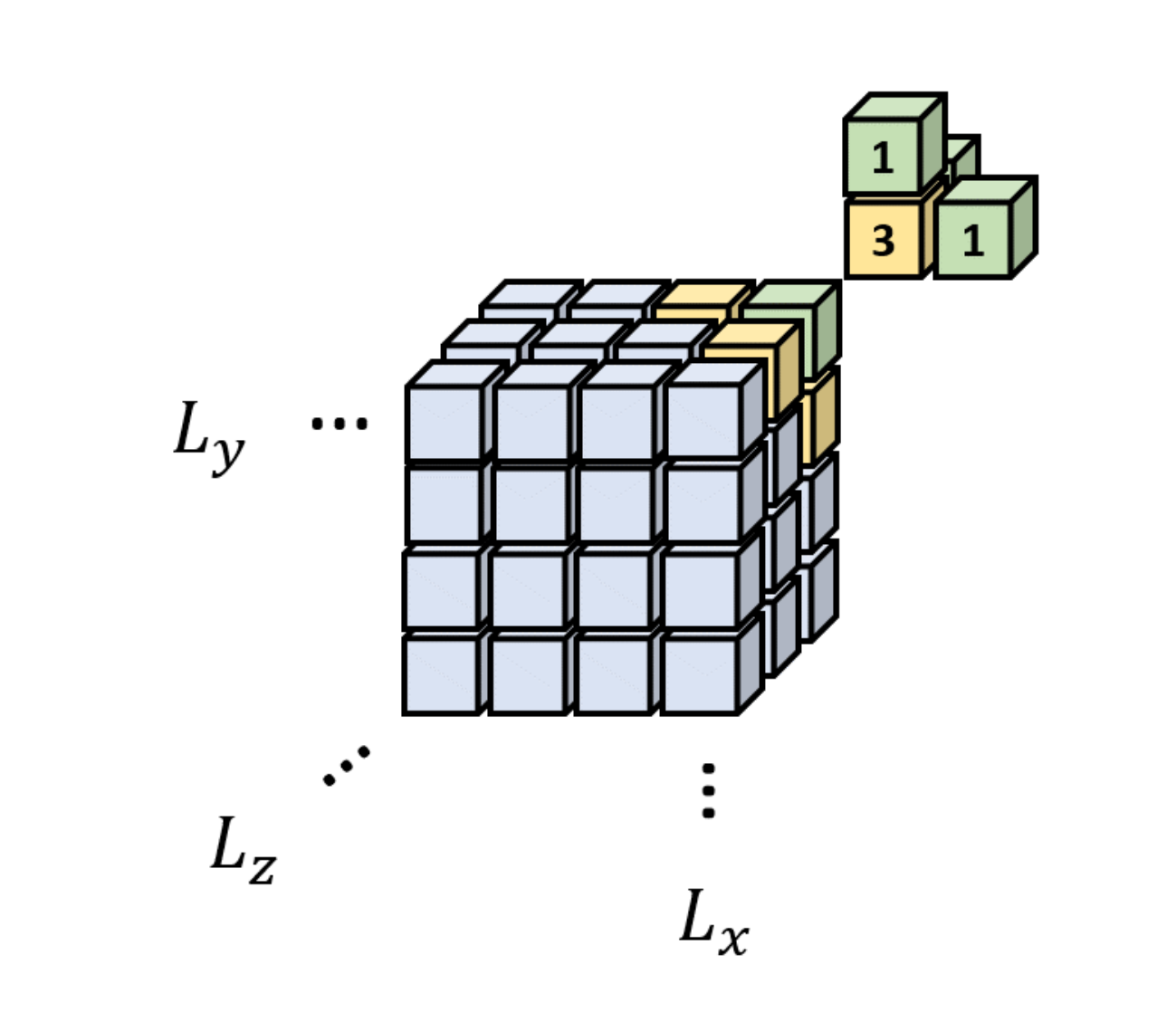}
				\caption{Examples of pairs of sites separated by distance $m=4$ that cross the $\xy$, $\yz$ and $\xz$ boundaries, \textit{i.e.} pairs of sites in $P_1^{xz}(4) \cap P_1^{yz}(4) \cap P_1^{xy}(4)$. Green and yellows sites inside of the cube domain are distance $m=3$ from green and yellow sites outside of the domain. The numbers correspond to the number of distance $m=4$ neighbours each site has which lie in $P_1^{xz}(4) \cap P_1^{yz}(4) \cap P_1^{xy}(4)$.}
				\label{fig:taxcabcorner}	}
	\end{center}
\end{figure}

Finally, we compute $|P_d^{yz}(m) \cap P_d^{xz}(m) \cap P_d^{xy}(m)|$, the number of pairs of sites of distance $m$ that cross the $\yz$, $\xz$ and $\xy$ plane boundaries. Notice that if $m<3$ no pairs can be connected by crossing all three boundaries, which means that $P_1^{yz}(m) \cap P_1^{xz}(m) \cap P_1^{xy}(m)$ is the empty set. 
Fig.\,\ref{fig:taxcabcorner} is an illustration of pairs of sites in $P_1^{yz}(4) \cap P_1^{xz}(4) \cap P_1^{xy}(4)$. These consist of the pairs of sites in four of the eight corners of the cube. In the case of the taxicab metric, the number of pairs of sites in $P_1^{yz}(m) \cap P_1^{xz}(m) \cap P_1^{xy}(m)$ can be calculated by first considering the site at the very corner of the cube (marked in green in Fig.\,\ref{fig:taxcabcorner}). There is exactly $Tr(1)=1$ site at the very corner of the cube where $Tr(x)$ refers to the $x$th triangular number. This site is distance $m$ from $Tr(m-2)$ sites (yellow sites outside of the domain) by crossing the $\yz$, $\xz$ and $\xy$ plane boundaries. Next we consider the nearest neighbouring sites to the very corner of the cube. There are $Tr(2)=3$ of these sites (marked in yellow). Each of the these sites are distance $m$ from $Tr(m-3)$ sites by crossing the $\yz$, $\xz$ and $\xy$ plane boundaries (marked in yellow outside of the domain). We continue this method to count site pairs until we consider sites that separated by $m-3$ sites from the corner site (green site in Fig.\,\ref{fig:taxcabcorner}). Hence we deduce that:
\begin{align}
\label{eq:res_taxicabcorner}
|P_1^{yz}(m) \cap P_1^{xz}(m) \cap P_1^{xy}(m)|=&4\sum_{i=1}^{m-2}Tr(i)Tr(m-i-1) \nonumber
\\=&4\sum_{i=1}^{m-2}\frac{i(i+1)}{2}\cdot \frac{(m-i-1)(m-i)}{2}\, , \nonumber
\\=&\frac{m^5}{30}-\frac{m^3}{6}+\frac{2m}{15}.
\end{align}

With a similar argument, one can obtain the corresponding expression by using the uniform metric which reads
\begin{align}
\label{eq:res_uniformcorner}
|P_\infty^{yz}(m) \cap P_\infty^{xz}(m) \cap P_\infty^{xy}(m)|=&m^3 (3 m^2+1) \, .
\end{align}

The formulae \eqref{eq:res_taxicabcorner} and \eqref{eq:res_uniformcorner} complete the formulae needed for the computation of the remainder of equation\,\eqref{eq:rem}. This can be used together with Eqs.\,\eqref{MAIN-eq:norm_3d} in the main text to rearrange equation\,\eqref{eq:def_rem}. The resulting expressions for the numbers of pairs at distance $m$ under non-periodic BC in three dimensions, for the taxicab and uniform metrics, respectively, read


\begin{subequations}
\label{eq:norm_3d_full}
\begin{align}
\begin{split}
s^n_{{cube}_1}(m)=&(2m^2+1)L_xL_yL_z-\frac{1}{3}(2m^3+m)(L_xL_y+L_yL_z+L_zL_x)\\&+\frac{m^2(m^2-1)}{6}(L_x+L_y+L_z)-\frac{m^5}{30}+\frac{m^3}{6}-\frac{2m}{15}\,,
\end{split}
\end{align}
and
\begin{align}
\begin{split}
s^n_{{cube}_\infty}(m)=&(12m^2+1)L_xL_yL_z-m(8m^2+1)(L_xL_y+L_yL_z+L_zL_x)\\&+m^2(5m^2+1)(L_x+L_y+L_z)-m^3 (3 m^2+1)  \, .
\end{split}
\end{align}
\end{subequations}
The normalisations for the Cube Taxicab PCF and Cube Uniform PCF under non-periodic BC are obtained by substituting expressions \eqref{eq:norm_3d_full} into equation\,\eqref{MAIN-normalisation_factor_1} in the main text.
\section{Summary of normalisation factors}
\label{sec:summary_norm}
In this section we summarise all of the normalisation factors which have been computed in the paper. For all cases considered, the normalisation, $\mathbb{E}[\bar{c}_d(m)]$, is obtained by equation\,\eqref{MAIN-normalisation_factor_1} in the main text. In Table \ref{table:summary} we report only the counts of pairs of sites separated by distance $m$, $s_d(m)$, depending on the metric and the BC used.

\begin{table}
\begin{center}
\begin{adjustbox}{angle=90}
\def\arraystretch{2}
\begin{tabular}{ |c|c|c|c| } 
\hline
PCF ($d$) &  Periodic ($s^p_d$)  & Non-periodic ($s^n_d$) \\[2ex]
\hline \hline 
Square Taxicab ($1$) & $2mL_xL_y$ & $2mL_xL_y-(L_x+L_y)m^2+\frac{m^3-m}{3}$ \\[2ex] 
\hline
Square Uniform ($\infty$) & $4mL_xL_y$ & $4mL_xL_y-3(L_x+L_y)m^2+2m^3$ \\[2ex] 
\hline
Triangle ($tri$) & $\frac{3}{2} mL_x L_y$ & 
\begin{tabular}{@{}c@{}}
for $m=1$: \quad  $3 \frac{L_x L_y}{2}-\frac{L_x}{2}-L_y$ \, ,
\quad for $m=2$: \quad  $3 L_x L_y-2L_x-4L_y+2$  \, , \\
for $m\ge3$: \quad  $3 m \frac{L_x L_y}{2}-\frac{m^2}{2}L_x+\left[2k_6(k_6-2k_3+1)+k_3(m-6)-m^2+m-2 \right]L_y$  
\\$+\frac{1}{3} (m-1)(m^2-2m+6)-\frac{1}{3}(k_7+1)(20k_7^2+37k_7+12)$\\$-(m-7-4k_7)(k_7+1)(m+k_7-2)$\,, \qquad where $k_j=\left\lfloor \frac{m-j}{4} \right\rfloor$
\end{tabular}
\\[2ex] 
\hline
Hexagon ($hex$) & $3mL_xL_y$ & 
\begin{tabular}{@{}c@{}}
$3mL_x L_y -\frac{1}{4} \left( 7 m^2+k\right)L_x-2m^2 L_y+\frac{11m^3}{12}-\frac{(2-3k)m}{12}$	 
\\ where $k=m \pmod{2}$
\end{tabular}
\\ [2ex] 
\hline
Cube Taxicab (${cube}_1$) & $(2m^2+1) L_x L_y L_z$ & 
\begin{tabular}{@{}c@{}}
$(2m^2+1)L_xL_yL_z-\frac{1}{3}(2m^3+m)(L_xL_y+L_yL_z+L_zL_x)$\\ 
$+\frac{m^2(m^2-1)}{6}(L_x+L_y+L_z)-\frac{m^5}{30}+\frac{m^3}{6}-\frac{2m}{15}$
\end{tabular}
\\[2ex] 
\hline
Cube Uniform (${cube}_\infty$) & $(12m^2+1)L_xL_yL_z$ & 
\begin{tabular}{@{}c@{}}
$(12m^2+1)L_xL_yL_z-m(8m^2+1)(L_xL_y+L_yL_z+L_zL_x)$\\ 
$+m^2(5m^2+1)(L_x+L_y+L_z)-m^3 (3 m^2+1) $
\end{tabular}
\\[2ex] 
\hline
\end{tabular}
\end{adjustbox}
\end{center}
\caption{Summary of the number of sites pairs, $s_d(m)$, according to the metric and BC used.}
\label{table:summary}
\end{table}
%
%
\newpage
\bibliography{database.bib}
\bibliographystyle{unsrtnat}